\documentclass[final,3p,times]{elsarticle}

\usepackage{amssymb}
\usepackage{multicol}
\usepackage[svgnames]{xcolor}

\usepackage{booktabs}
\usepackage{array}
\newcommand{\PreserveBackslash}[1]{\let\temp=\\#1\let\\=\temp}
\newcolumntype{C}[1]{>{\PreserveBackslash\centering}p{#1}}
\newcolumntype{R}[1]{>{\PreserveBackslash\raggedleft}p{#1}}
\newcolumntype{L}[1]{>{\PreserveBackslash\raggedright}p{#1}}

\journal{Computer Methods in Applied Mechanics and Engineering (CMAME)}

\usepackage{amsmath}
\usepackage{amsfonts}
\usepackage{amsthm}
\usepackage{mathtools}
\usepackage{subfig}
\usepackage{lineno}
\usepackage[hidelinks]{hyperref}
\usepackage{textcomp}
\usepackage{algorithm,algpseudocode,algorithmicx}

\usepackage{float}
\usepackage{array,multirow}
\usepackage{color}
\usepackage{tikz}
\usetikzlibrary{patterns}
\usetikzlibrary{spy}
\usetikzlibrary{spy,decorations.fractals,shapes.misc}

\usepackage{pgfplots}
\pgfplotsset{compat=1.16}
\newsavebox\Axis

\definecolor{lightgray}{gray}{0.80}

\usetikzlibrary{decorations.pathreplacing}
\usepackage[listings,theorems,skins,breakable]{tcolorbox}
\tcbset{skin=enhanced}
\newtcolorbox{lbracebox}[1][Word]{%
   frame hidden,enlarge left by=2cm,width=\linewidth-2cm,%
  overlay unbroken = {\draw [decorate,decoration={brace,amplitude=10pt},]%
                     (frame.south west)-- (frame.north west)
                    node [black,midway,left,xshift=-.6cm] {#1};},%
}

\definecolor{green1}{rgb}{0.4660, 0.6740, 0.1880} 
\definecolor{blue1}{rgb}{0, 0.4470, 0.7410} 
\definecolor{red1}{rgb}{0.8500, 0.3250, 0.0980}
\definecolor{yellow1}{rgb}{0.9290, 0.6940, 0.1250}
\definecolor{purple1}{rgb}{0.4940, 0.1840, 0.5560}
\definecolor{lightblue1}{rgb}{0.3010, 0.7450, 0.9330}
\definecolor{bordeaux1}{rgb}{0.6350, 0.0780, 0.1840}
\definecolor{gray1}{rgb}{0.7, 0.7, 0.7}
\definecolor{gray2}{rgb}{0.5, 0.5, 0.5}

\definecolor{burntorange}{rgb}{0.74902,0.341176,0}

\theoremstyle{plain}
\newtheorem{theorem}{Theorem}[section]

\theoremstyle{definition}
\newtheorem{definition}[theorem]{Definition}

\newcommand\tikzcircle[1][2.5]{\tikz[baseline=-#1]{\draw[thick](0,2pt) circle [radius=#1mm];}}

\newcommand{\vect}[1]{\boldsymbol{#1}} 									
\newcommand{\mat}[1]{\mathsf{#1}} 											
\newcommand{\mmat}[1]{\mathbf{#1}} 											

\newcommand{\domain}{\tikzcircle[1.5]}														
\newcommand{\N}{N}																				

\newcommand{\p}{p}																	

\newcommand{\sspace}[2]{\mathbb{S}^{#1}_{#2}}					

\newcommand{\bilinearform}[3]{#1(#2,\,#3)}				
\renewcommand{\aa}[2]{\bilinearform{a}{#1}{#2}}						
\newcommand{\bb}[2]{\bilinearform{b}{#1}{#2}}						
\newcommand{\cc}[2]{\bilinearform{c}{#1}{#2}}						

\modulolinenumbers[5]

\begin{document}

\begin{frontmatter}

\title{Fourier analysis of membrane locking and unlocking}

\corref{cor1}
\author[mech]{Ren\'e R. Hiemstra}
\ead{rene.hiemstra@tu-darmstadt.de}

\author[imc]{Federico Fuentes}
\ead{federico.fuentes@uc.cl}

\author[mech]{Dominik Schillinger}
\ead{dominik.schillinger@tu-darmstadt.de}

\cortext[cor1]{Corresponding author}

\address[mech]{Institute for Mechanics, Computational Mechanics Group, Technical University of Darmstadt, Germany}
\address[imc]{Institute for Mathematical and Computational Engineering (IMC), School of Engineering and Faculty of Mathematics, Pontificia Universidad Cat\'olica de Chile, Chile
\\[1.2cm]
\begin{large} 
Dedicated to Thomas J.R.\ Hughes on the occasion of his 80th birthday.
\end{large}
}

\begin{abstract}
Membrane locking in finite element approximations of thin beams and shells has remained an unresolved topic despite four decades of research. In this article, we utilize Fourier analysis of the complete spectrum of natural vibrations and propose a criterion to identify and evaluate the severity of membrane locking. To demonstrate our approach, we utilize standard and mixed Galerkin formulations applied to a circular Euler-Bernoulli ring discretized using uniform, periodic B-splines. By analytically computing the discrete Fourier operators, we obtain an exact representation of the normalized error across the entire spectrum of eigenvalues. Our investigation addresses key questions related to membrane locking, including mode susceptibility, the influence of polynomial order, and the impact of shell/beam thickness and radius of curvature. Furthermore, we compare the effectiveness of mixed and standard Galerkin methods in mitigating locking. By providing insights into the parameters affecting locking and introducing a criterion to evaluate its severity, this research contributes to the development of improved numerical methods for thin beams and shells.
\end{abstract}


\end{frontmatter}

\newpage

\section{Introduction \label{sec:introduction}}

\begin{em}
It was late 2009 when I first heard about isogeometric analysis. As a fresh M.Sc. student with a keen interest in research, particularly in finite element analysis, I delved into the first few papers by Tom and his former students. It became evident that this technology was a game-changer. I promptly acquired a copy of the IGA book released earlier that year and absorbed its contents in a matter of weeks. The prospect of working with such celebrated researchers ignited my imagination.

Our paths converged two years later at the World Congress on Computational Mechanics in Vienna. Tom, a prominent figure in the field, delivered the opening talk to a crowd of roughly 4000 people. After his presentation, I seized the opportunity to approach him, and to my delight, he generously engaged in a conversation. It was during this interaction that I discovered Tom's genuine interest in my research topic. Two days later, as I presented my work, I was pleasantly surprised to find Tom among the audience. Despite my initial nerves, I managed to deliver a solid presentation. This pivotal moment led to a remarkable turn of events—Tom extended an invitation to visit his research group in Austin, an invitation that effectively offered me a Ph.D. opportunity.

I had the incredible fortune of working alongside Tom, and it was a truly remarkable experience. Our initial encounter turned out to be just the beginning of a series of memorable interactions. Tom possesses a unique quality—he always makes time for his students and others, engaging with them not only on scientific matters but also on personal and amusing topics. He is quite the collector of anecdotes, regaling us with tales from the ``good old days" of finite elements in the late sixties and seventies.

Allow me to share two of Tom's captivating stories, a testament to his sense of humor and ability to find joy in unexpected situations. Picture a young Tom in an elevator with the late Oleg Zienkiewicz and a youthful Robert Taylor. In this particular episode, Robert earnestly tries to explain to Oleg why his penalty formulation is fundamentally flawed, emphasizing that ``the solution is zero to all problems". However, Zienkiewicz interrupts him with a playful remark, urging Robert not to bother him with details and to just implement the idea. In another story from a conference in western Germany, Tom recounts how a group of young scientists, including Ted Belytschko, Jurgen Bathe, and Wing Kam Liu, found themselves scribbling equations in a caf\'e. Unwittingly, their activity raised suspicion, and they were mistaken for terrorists, leading to an encounter with the police. It's these glimpses into Tom's past that reveal not only his expertise but also his wit and capacity to navigate scientific collaborations with a light-hearted approach.

Dear Tom, on the momentous occasion of your 80th birthday, I extend my warmest congratulations to you. As you gracefully embrace this respectable age, it's worth noting that you possess a vibrant spirit that defies the passing years. I dedicate this contribution to you, considering your unwavering interest in shell finite elements, spectral analysis, and the intriguing subject of membrane locking. I believe you will thoroughly enjoy reading it and find it enriching.
With warm regards, Rene
\end{em}

\

Membrane locking in finite element formulations of curved beam and shell models refers to a manifestation of artificial bending stiffness resulting from the coupling between the bending and membrane response induced by the curvature \cite{bischoff2004models, stolarski1982membrane}. In an extreme case, a curved beam, which would be expected to exhibit a bending-dominated response, instead displays an in-plane membrane-dominated discrete solution. Membrane locking is only one of several sources of locking, the most well-known being transverse shear locking in beam, plate and shell elements \cite{bischoff2004models,suri1995locking} and volumetric locking due to incompressibility in solid elements \cite{babuvska1992locking2, hughes2012finite}. 
All sources of locking may severely affect accuracy and robustness of finite element discretizations. 
The development of locking-preventing discretization technology
has a history of more than 40 years, first within classical finite elements and more recently in isogeometric analysis. Common techniques include, but are not limited to, reduced and selective integration techniques, strain projection methods, and certain mixed formulations based on alternative variational principles. An extensive discussion on several forms of locking in relation to finite element design can be found in \cite[Chapter 5.4]{bischoff2004models} and \cite[Chapter 9.4, and Chapter 10]{wriggers2008nonlinear}. An analysis of locking effects from the mathematical point of view can be found in \cite[Chapter 6]{braess2007finite}. We also refer to \cite{ambroziak2013locking} for a concise overview of the relevant literature until 2013 and to \cite[Table 1]{nguyen2022leveraging} for an overview of more recent approaches, in particular in the context of isogeometric analysis.

From a numerical analysis viewpoint, locking has been described as a loss of robustness (or uniform approximation) as a problem dependent parameter $t$ approaches a limiting value $t_0$ \cite{babuvska1992locking2,arnold1981discretization}. For example, Poisson's ratio approaching the incompressibility limit of one half and plate or shell thickness approaching zero. Most a priori error estimates for elliptic boundary value problems yield optimal asymptotic convergence rates when the parameter $t \neq t_0$ is fixed. 
This fact does not necessarily aid the study of locking phenomena, because these often manifest at coarse discretizations and vanish with mesh refinement. A series of early works \cite{suri1995locking,babuvska1992locking2, arnold1981discretization, babuvska1992locking1} are devoted to the mathematical analysis of finite element methods in the context of parameter dependency and locking. It has been shown how solid elements lose uniform convergence in the incompressibility limit \cite{babuvska1992locking2} and how uniform approximation is lost due to shear locking as the beam or plate thickness approaches zero \cite{suri1995locking,arnold1981discretization}.
However, membrane locking as a standalone phenomenon has not received the same level of mathematical scrutiny as the two other main sources of locking. The limited attention to membrane locking can be attributed, in part, to its physical rather than mathematical definition. It refers to the propensity for a higher presence of membrane stresses or displacements, as opposed to bending stresses and displacements, in comparison to the exact solution. Consequently, its distinct characteristics and mathematical analysis have not been thoroughly explored in the literature.

Benchmark problems designed to evaluate the performance and correctness of beam and shell elements \cite{belytschko1985stress, macneal1985proposed} have proven indispensable in the development of robust shell elements, but alone are not adequate to assess membrane locking and unlocking. These benchmarks cannot possibly cover all conditions and can therefore not positively rule out membrane locking. In addition, due to the complexity in loading and boundary conditions, these tests cannot isolate and separate the effect of membrane locking from other parasitic behavior due to, e.g., boundary layers.

Recently, we presented a novel approach to enhance the understanding of  membrane locking and unlocking in thin beams and shells through spectral analysis \cite{nguyen2022leveraging}. The eigenvalue decomposition is a valuable tool because it offers a comprehensive understanding of a linear operator\footnote{Consider a well posed linear elliptic boundary value problem: find $u \in V$ that satisfies $a(u,v) = l(v)$ for any $v \in V$. Suppose $u_n \in V$ are given and satisfy $a(u_m, u_n) = \lambda_n \delta_{mn}$, then the solution allows the decomposition: $u = \sum \nolimits_{n} \alpha_n \cdot u_n$, where $\alpha_n = l(u_n) / \lambda_n$ whenever $\lambda_n \neq 0$ and free otherwise. The functions $u_n \in V$ are called the eigenfunctions or eigenmodes of $a$ and $\lambda_n = a(u_n,u_n)$ its eigenvalues.}. By analyzing the eigenvalues and eigenfunctions, we can gain valuable insights into the properties and behavior of the operator, without being confined by the specific details of a boundary value problem, such as boundary conditions and the right-hand side. Eigenvalue and modal analysis provides a holistic view of the discretization error associated with a method. Specifically, the discretization error, such as the one occurring in the energy norm of an elliptic boundary value problem, can be expressed solely in terms of the errors made in eigenvalues and modes  \cite{hughes_finite_2014,strang_analysis_2008}.

In \cite{nguyen2022leveraging} we examined five different techniques to prevent or mitigate membrane locking in the context of a curved Euler-Bernoulli circular ring. The simplicity and periodic nature of the model problem made it possible to isolate the effect of membrane locking from shear locking and spurious boundary phenomena. Locking was measured, at fixed normalized mode number, as the discrepancy between the spectrum error and the spectrum error in the limit of refinement. The severity of locking and the performance of locking mitigation or removal could thereby be measured in a quantitative manner, which enabled an objective comparison. We note that the use of spectral or Fourier analysis in relation to locking is not a completely new idea. Inspection of the low eigenvalues to investigate the presence and effect of locking has been considered in the design of efficient shell elements \cite{lee2015modal,zou2021galerkin} or multigrid methods \cite{braess1988multigrid}. Uncoupled discrete Fourier operators and corresponding characteristic equations have been used to identify spurious effects and locking in early work \cite{park1984fourier}.

In this article, we provide a comprehensive analysis of the discrete error through Fourier (spectral) analysis, revealing important insights into locking phenomena. We analytically show that the numerical solution of the isolated membrane or bending response exhibits normalized spectra that remain invariant regardless of mesh refinement. However, our model displays a coupled response for which the normalized spectra no longer possess this invariance, except in the limit of mesh refinement. This observation motivates the definition of a concrete criterion that characterizes the presence and severity of membrane locking as a function of the normalized spectrum. To illustrate the issue, we investigate a standard Galerkin formulation employing quadratic, cubic, and quartic spline function spaces, and demonstrate its susceptibility to membrane locking. In contrast, a mixed formulation based on the Hellinger-Reissner variational principle proves to be free from locking for all polynomial degrees and offers significantly improved absolute accuracy across all levels of mesh refinement.

\begin{figure}[t]
	\centering
	\includegraphics[width=0.8\textwidth]{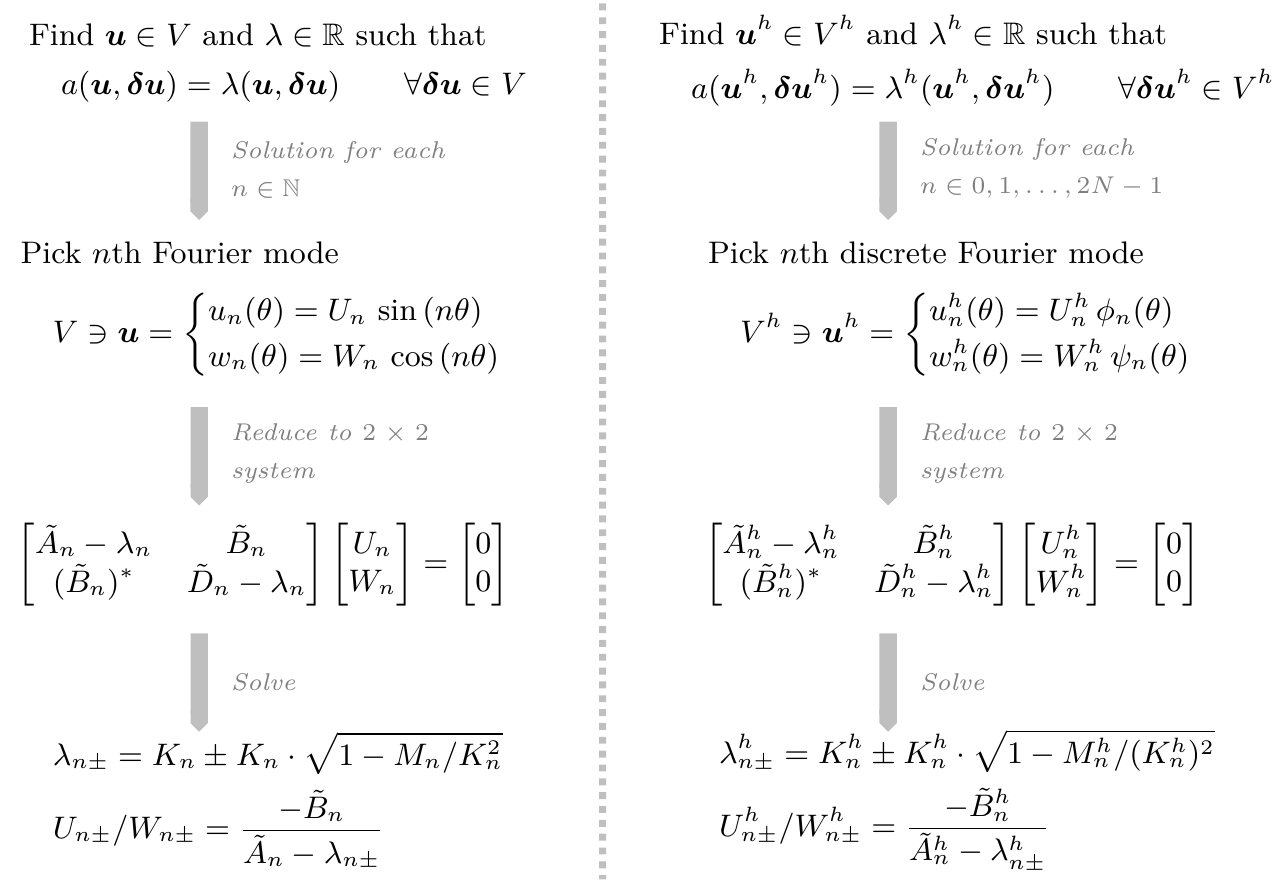}
	\caption{The solution approach of the discrete eigenvalue problem mirrors the solution approach to the exact solution.}
	\label{fig:overview}
\end{figure}

Our solution approach revolves around an analytical technique to solve the discrete eigenvalue problem, leveraging properties of circulant matrices that arise during discretization with periodic uniform spline basis functions. We establish the existence of a discrete analogue to the Fourier basis, effectively decoupling the $2N \times 2N$ system of equations into $N$ systems of $2 \times 2$ equations. The solution approach, depicted in Figure \ref{fig:overview}, closely resembles the exact solution obtained using a Fourier basis. This analytical approach enables us to compute exact expressions of the normalized error across the complete eigenvalue spectrum, incorporating parameters such as element size, thickness, and radius of curvature. Our analysis reveals how these relevant physical and discretization parameters contribute to membrane locking within the standard formulation. Notably, while our focus lies on quadratic discretizations of the standard formulation, our techniques apply to general polynomial degrees and can be adapted to the mixed formulation as well.

Although the dynamics of a circular Euler-Bernoulli ring are well-documented, as detailed in \cite[Chapter 5, Section 5.3, page 82]{soedel_vibrations_2004}, a comparison between the analytical and discrete spectra, for both membrane-dominated and bending-dominated eigenmodes, offers fresh insights into locking phenomena. Significantly, our findings establish that membrane locking solely affects eigenvalues corresponding to modes with a bending-dominated response, effectively overdampening the contributions of these eigenmodes from the discrete solution, while leaving the eigenvalues associated to a membrane-dominated response mostly unaffected. We also find that the phenomenon of membrane locking depends nontrivially on an interplay between a normalized thickness and the discretization size, and in a particular regime, is mathematically analogous to the phenomenon of shear locking (despite the fact that there are no shear strains in our model).

The paper is organized as follows: in Section \ref{sec:2}, we present the standard and mixed formulation of the eigenvalue model problem and investigate the behavior of the analytical solution for our model problem of a circular ring. In Section \ref{sec:3}, we present the standard and mixed Galerkin discretization using periodic uniform splines, and discuss an analytical solution technique based on the Discrete Fourier Transform (DFT). Section \ref{sec:4} introduces a criterion to identify and measure the severity of membrane locking in thin beams and shells. It is based on the discrepancy of the normalized spectrum of eigenvalues in relation to the spectrum in the limit of mesh refinement. These notions are applied in Section \ref{sec:5} to the numerical results stemming from the computation of the spectra of the circular Euler-Bernoulli ring. Here we juxtapose the standard formulation with the mixed formulation and discuss performance with respect to accuracy and robustness. In Section \ref{sec:6}, we perform a mathematical analysis of the spectrum in the standard formulation and derive exact expressions for the relative error and its dependence on different parameters. In Section \ref{sec:7} we further discuss some observations on the phenomenon of membrane locking in the context of our study. Finally, in Section \ref{sec:8} we draw conclusions, and propose potential future work.
\section{The continuous eigenvalue problem \label{sec:2}}
In the scope of this paper, we study the phenomenon of membrane locking in the special case of an Euler-Bernoulli ring. We highlight that this problem rules out the presence of shear strains, so is ideal to isolate membrane-related phenomena. In this section, we introduce the continuous mixed and standard eigenvalue problem. We then discuss properties of the analytical solution.

\begin{figure}[!h]
	\centering
	\captionsetup[subfigure]{labelformat=empty}
	\includegraphics[width=0.75\textwidth]{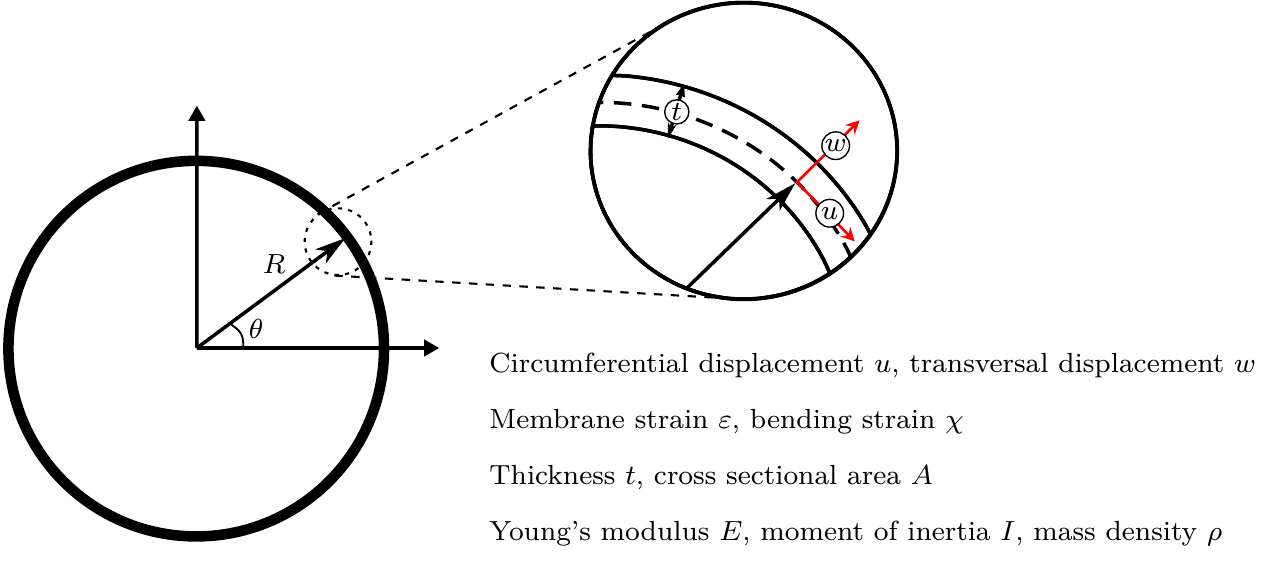}
	\caption{Closed circular ring modeled as a curved Euler-Bernoulli beam.}
	\label{fig:ring}
\end{figure}

\subsection{The mixed formulation}
\begin{subequations}
\label{eq:eigenvalue_problem}
Let $\domain$ denote the circular ring of radius $R$ displayed in Figure \ref{fig:ring}, parameterized by arc-length parameter $s \in [0,2\pi R)$. The field variables are the circumferential displacement component $u$, the transverse displacement component $w$, the membrane strain resultant $\varepsilon$, and the bending strain resultant $\chi$. Displacements are considered as elements from the space $V = H^2(\domain) \times H^2(\domain)$ and strains as elements in $S = L^2(\domain) \times L^2(\domain)$.

We consider the following eigenvalue problem: find $\vect{u} = (u,w) \in V$, $\vect{\varepsilon} = (\varepsilon,\chi) \in S$ and $\lambda \in \mathbb{R}^{+}_{0}$ such that
\begin{align}
	- \cc{\vect{\varepsilon}}{\vect{\delta \varepsilon}} + \cc{B(\vect{u})}{\vect{\delta \varepsilon}} &= 0 & \forall \vect{\delta \varepsilon} &= (\delta \varepsilon, \delta \chi) \in S \label{eq:eigenvalue_problem_kinematics} \\
	  \cc{\vect{\varepsilon}}{B(\vect{\delta u})} &= \lambda \bb{\vect{u}}{\vect{\delta u}} 	 & \forall \vect{\delta u} &= (\delta u, \delta w) \in V. \label{eq:eigenvalue_problem_equation}
\end{align}
where
\begin{align}
	\cc{\vect{\varepsilon}}{\vect{\delta \varepsilon}} &:= \int_{0}^{2\pi R} \left(EA \, \varepsilon \, \delta \varepsilon + EI \, \chi \, \delta \chi \right) \, ds &  \vect{\varepsilon}, \vect{\delta \varepsilon} \in S 	\label{eq:stiffness}	\\
	\bb{\vect{u}}{\vect{\delta u}} &:= \int_{0}^{2\pi R} \rho A \left( u \, \delta u + w \, \delta w\right) \, ds
&  \vect{u}, \vect{\delta u} \in V. \label{eq:mass}
\end{align}
Here $E, A, I$ and $\rho$ denote, respectively, the Young's modulus, the cross sectional area and moment of inertia, and the density of the ring. Furthermore, $B \; : \; V \mapsto S$ is defined via the kinematic relations
\begin{align}
	B (\vect{u}) := \begin{cases}
	\varepsilon 	(\vect{u}) 	&= \frac{\partial u}{\partial s} + \frac{1}{R} w \\
	\chi 				(\vect{u}) &= - \frac{\partial^2 w}{\partial s^2}  + \frac{1}{R} \frac{\partial u}{\partial s}
\end{cases}
\end{align}

Membrane locking occurs due to the disparate nature of work due to membrane stresses (first term in \eqref{eq:stiffness}) and work due to bending (second term in \eqref{eq:stiffness}). To study the phenomenon we are interested only in the relative magnitude of these terms. In the following we transform the problem to polar coordinates $(r, \theta)$, factor out the density and Young's modulus, and assume a rectangular cross section with constant thickness, that is,
\begin{align}
	&E/\rho = 1,&
	&t = \text{const},&
	&I = \frac{A t^2}{12},&
	&
	B (\vect{u}) := \begin{cases}
	\varepsilon 	(\vect{u}) 	&= \frac{1}{R} \frac{\partial u}{\partial \theta} + \frac{1}{R} w \\
	\chi 				(\vect{u}) &= -  \frac{1}{R^2} \frac{\partial^2 w}{\partial \theta^2}  + \frac{1}{R^2} \frac{\partial u}{\partial \theta}
\end{cases}&
\end{align}
The bilinear forms, $c \, : \, S \times S \mapsto \mathbb{R}$ and $b \, : \, V \times V \mapsto \mathbb{R}$, then simplify to
\begin{align}
	\cc{\vect{\varepsilon}}{\vect{\delta \varepsilon}} &:= \int_{0}^{2\pi} \left(\varepsilon \, \delta \varepsilon + \tfrac{t^2}{12} \, \chi \, \delta \chi \right) \, R \, d \theta &  \vect{\varepsilon}, \vect{\delta \varepsilon} \in S 	\label{eq:stiffness2} \\
	\bb{\vect{u}}{\vect{\delta u}} &:= \int_{0}^{2\pi} \left( u \, \delta u + w \, \delta w\right) \, R \, d \theta
&  \vect{u}, \vect{\delta u} \in V. \label{eq:mass2}
\end{align}
\end{subequations}

\subsection{The standard formulation}
Equivalently, we may consider the standard formulation obtained from the mixed formulation by substituting $\vect{\varepsilon} = B(\vect{u})$ into equation \eqref{eq:eigenvalue_problem_equation}. We seek $\vect{u} = (u,w) \in V$ and $\lambda \in \mathbb{R}^{+}_{0}$ such that
\begin{subequations}
\label{eq:standard_formulation}
\begin{align}
	\aa{\vect{u}}{\vect{\delta u}} = \lambda \bb{\vect{u}}{ \vect{\delta u}} \qquad \forall \vect{\delta u} = (\delta u, \delta w) \in V.
\end{align}
Here $a \, : \, V \times V \mapsto \mathbb{R}$ denotes the semi-definite symmetric bilinear form defined as
\begin{align}
	\aa{\vect{u}}{\vect{\delta u}} := \cc{B(\vect{u})}{B(\vect{\delta u})}, \qquad  \vect{u}, \vect{\delta u} \in V.
	\label{eq:bilinearform}
\end{align}
\end{subequations}
The standard and mixed formulations are equivalent, that is, they have the same eigenvalues and eigenfunctions.

\subsection{Analytical computation of exact eigenvalues and eigenfunctions}
\label{sec:analyticalEVEFs}
\begin{figure}[!h]
	\centering
	\captionsetup[subfigure]{labelformat=empty}
	\subfloat[$n = 1$]{\includegraphics[width=0.11\textwidth]{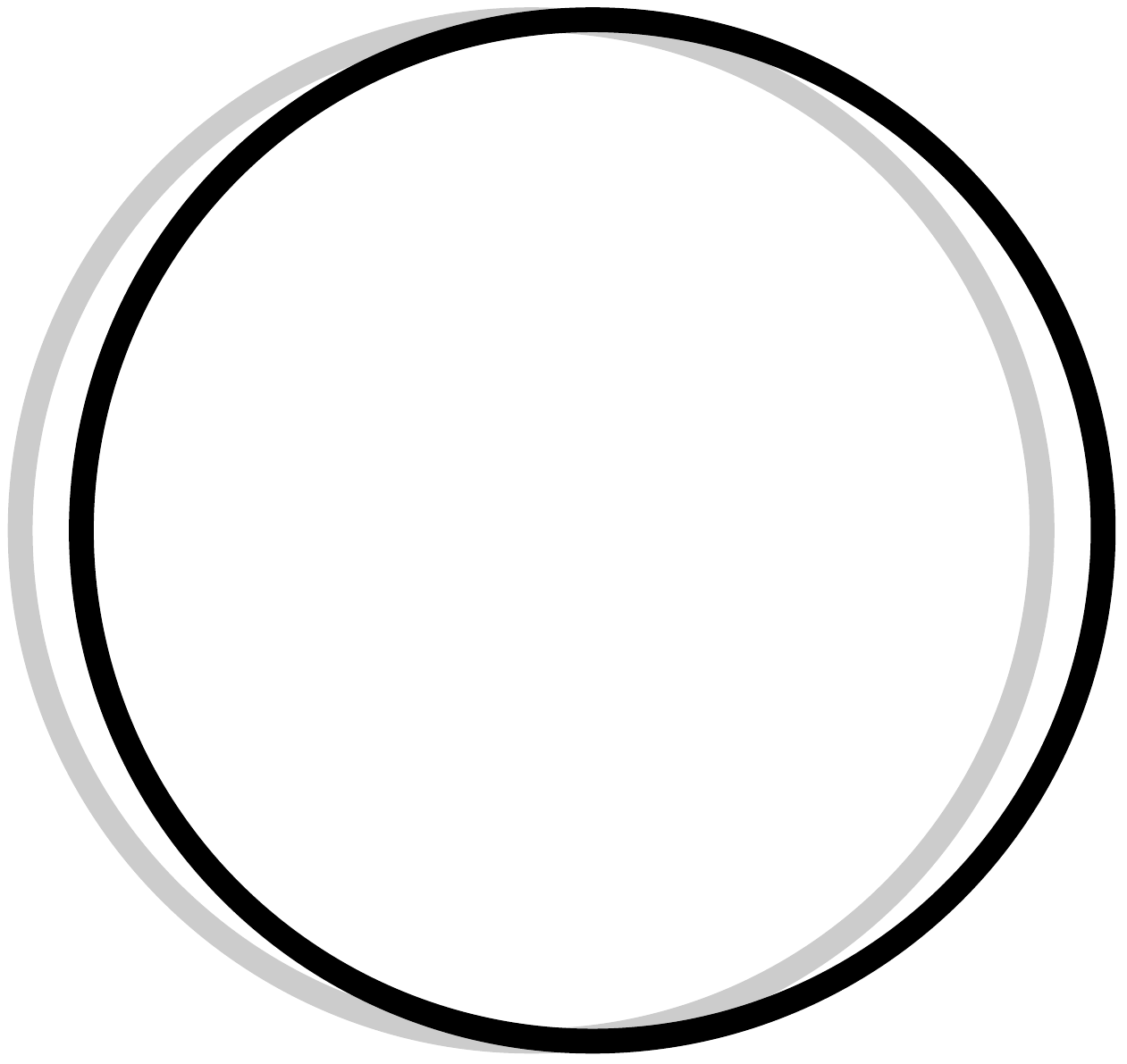}}	\hspace{0.5cm}
	\subfloat[$n = 2$]{\includegraphics[width=0.11\textwidth]{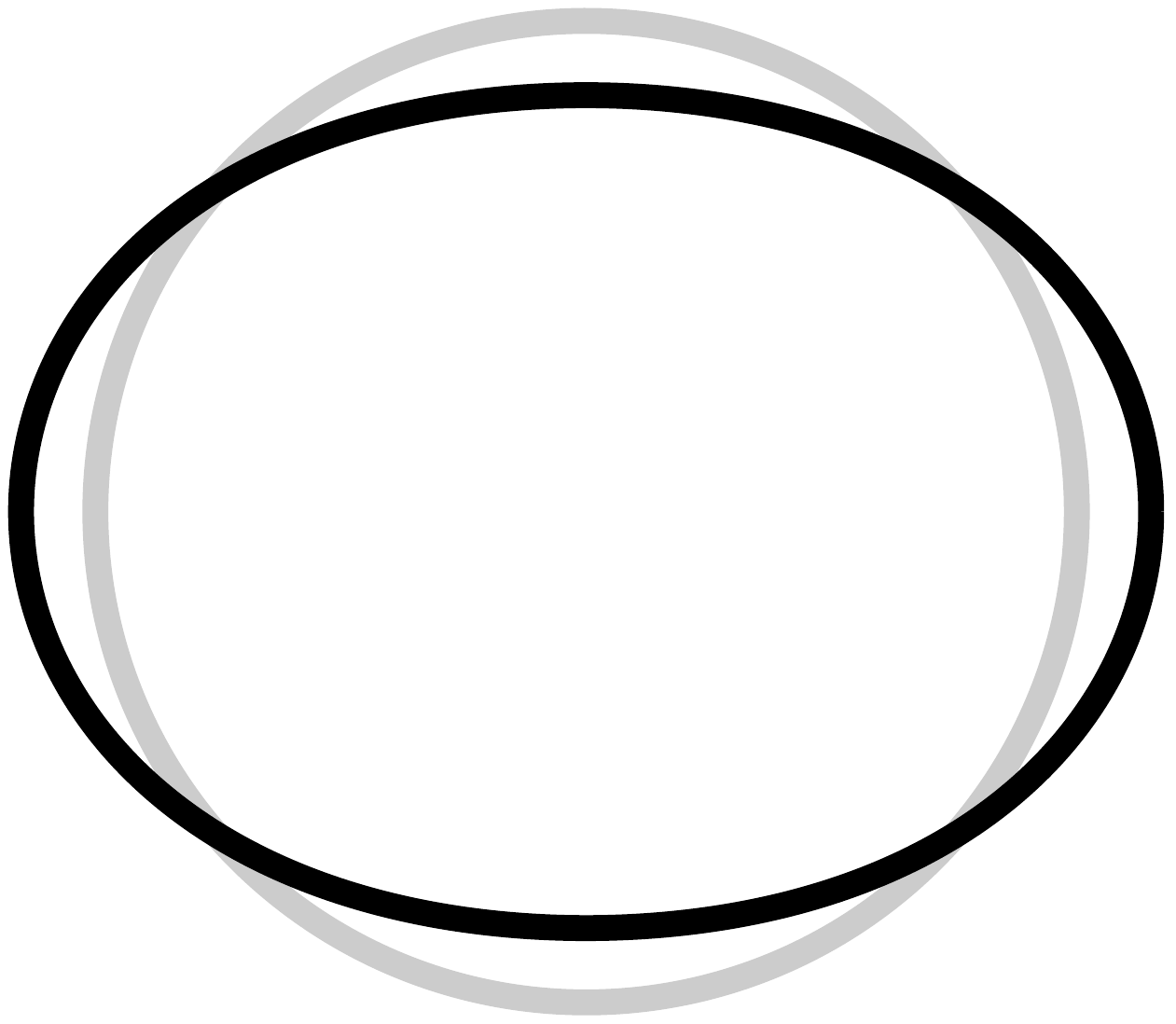}}	\hspace{0.5cm}
	\subfloat[$n = 3$]{\includegraphics[width=0.11\textwidth]{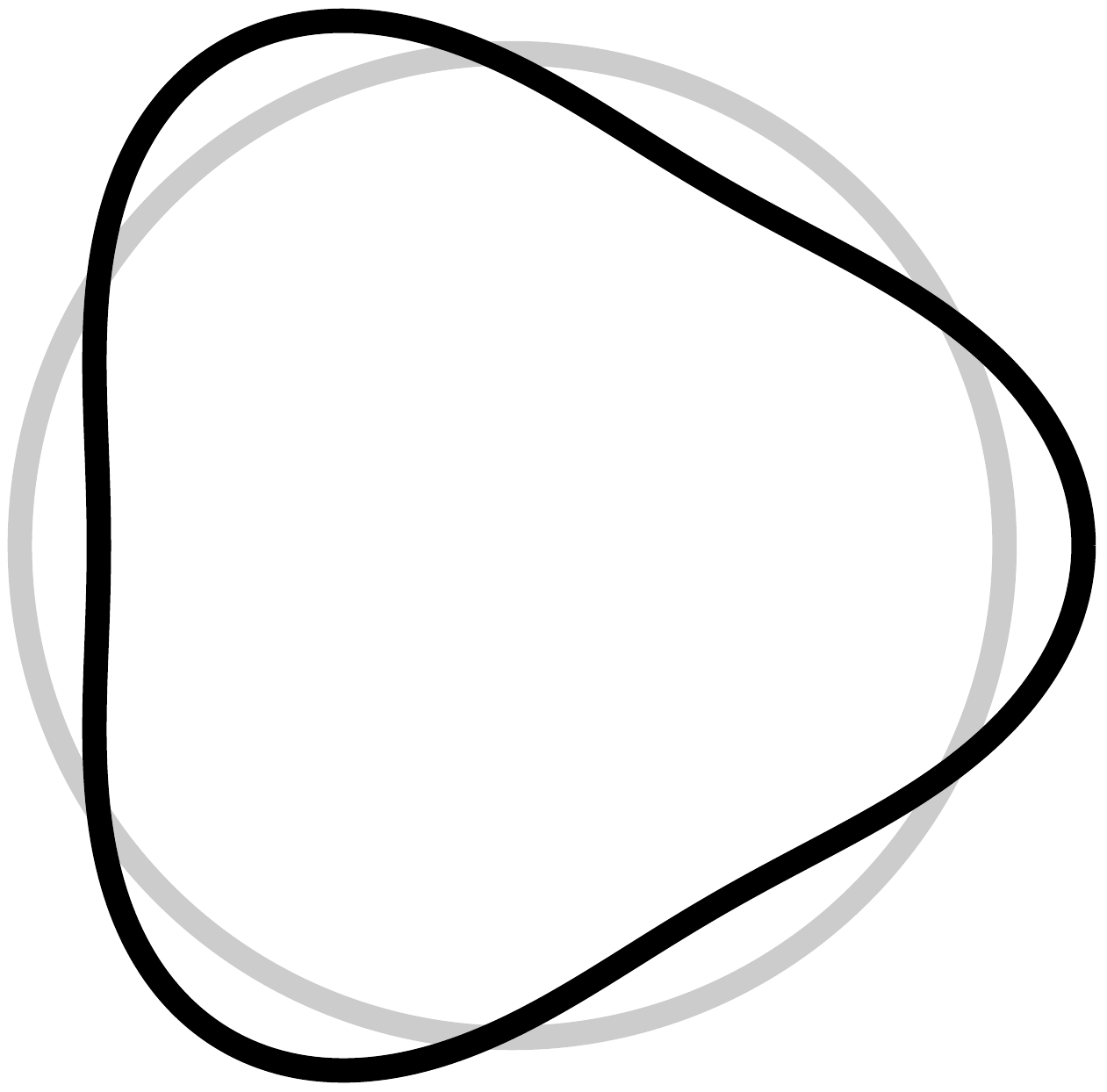}}	\hspace{0.5cm}
	\subfloat[$n = 4$]{\includegraphics[width=0.11\textwidth]{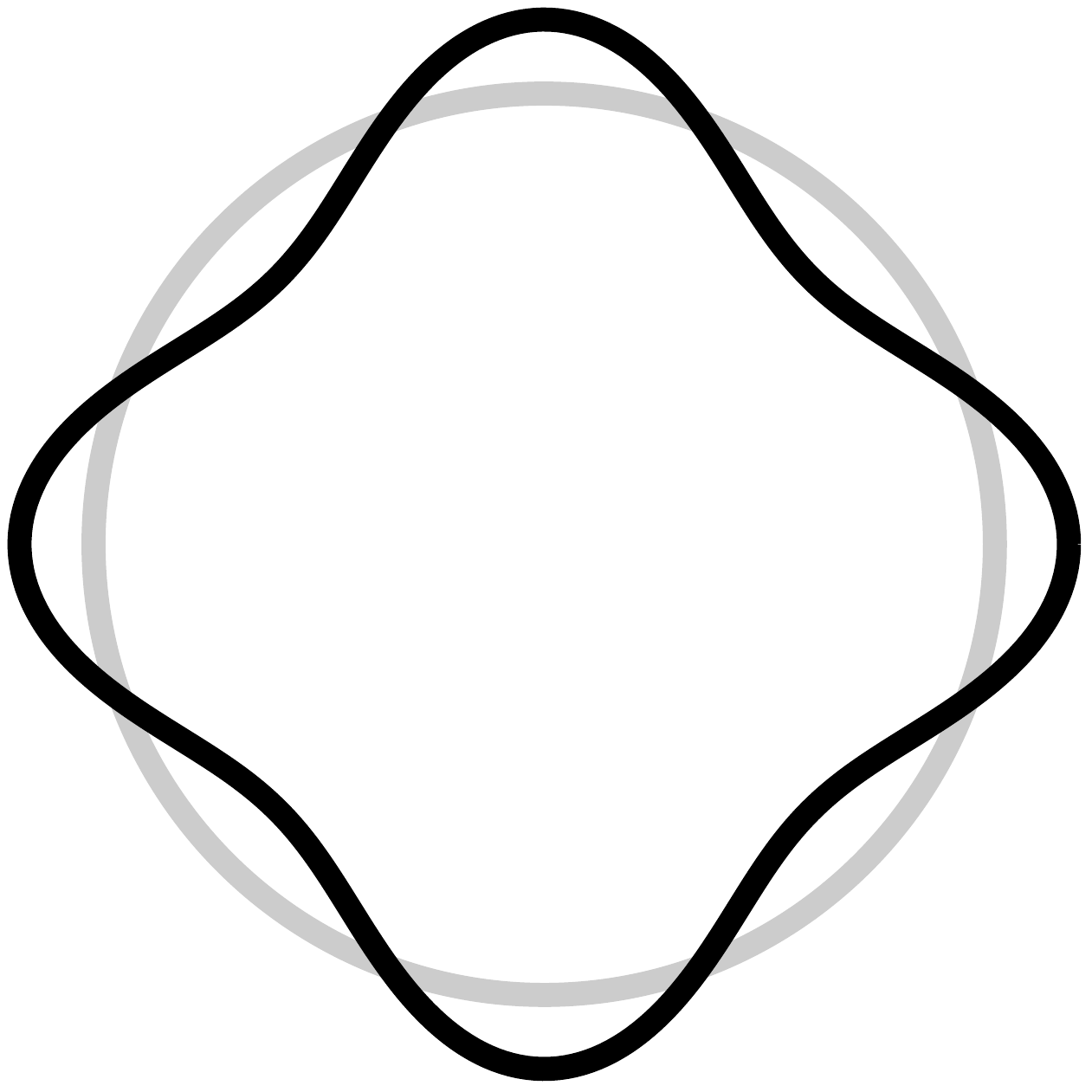}}  \hspace{0.5cm}
	\subfloat[$n = 5$]{\includegraphics[width=0.11\textwidth]{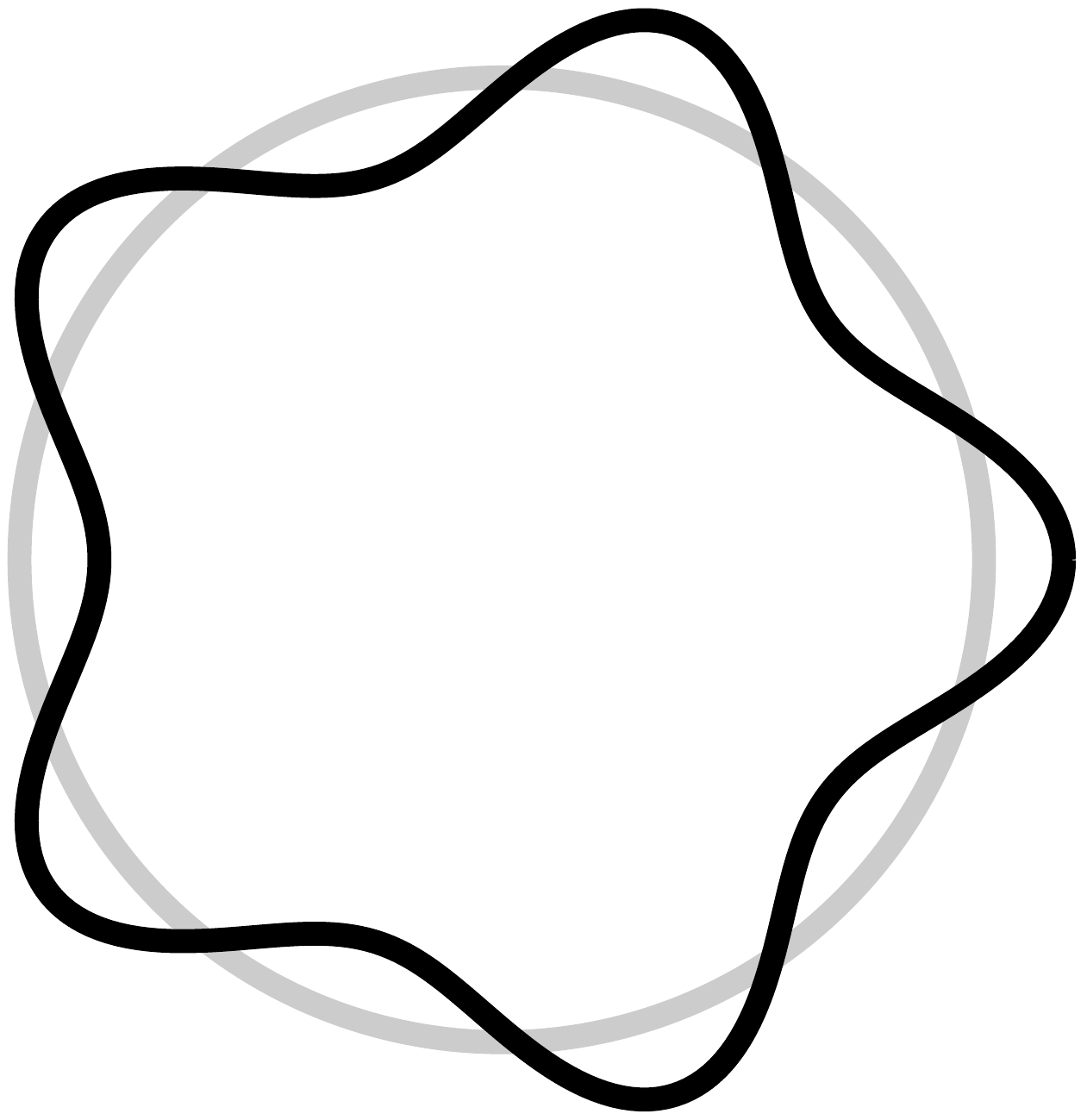}}	\hspace{0.5cm}
	\subfloat[$n = 6$]{\includegraphics[width=0.11\textwidth]{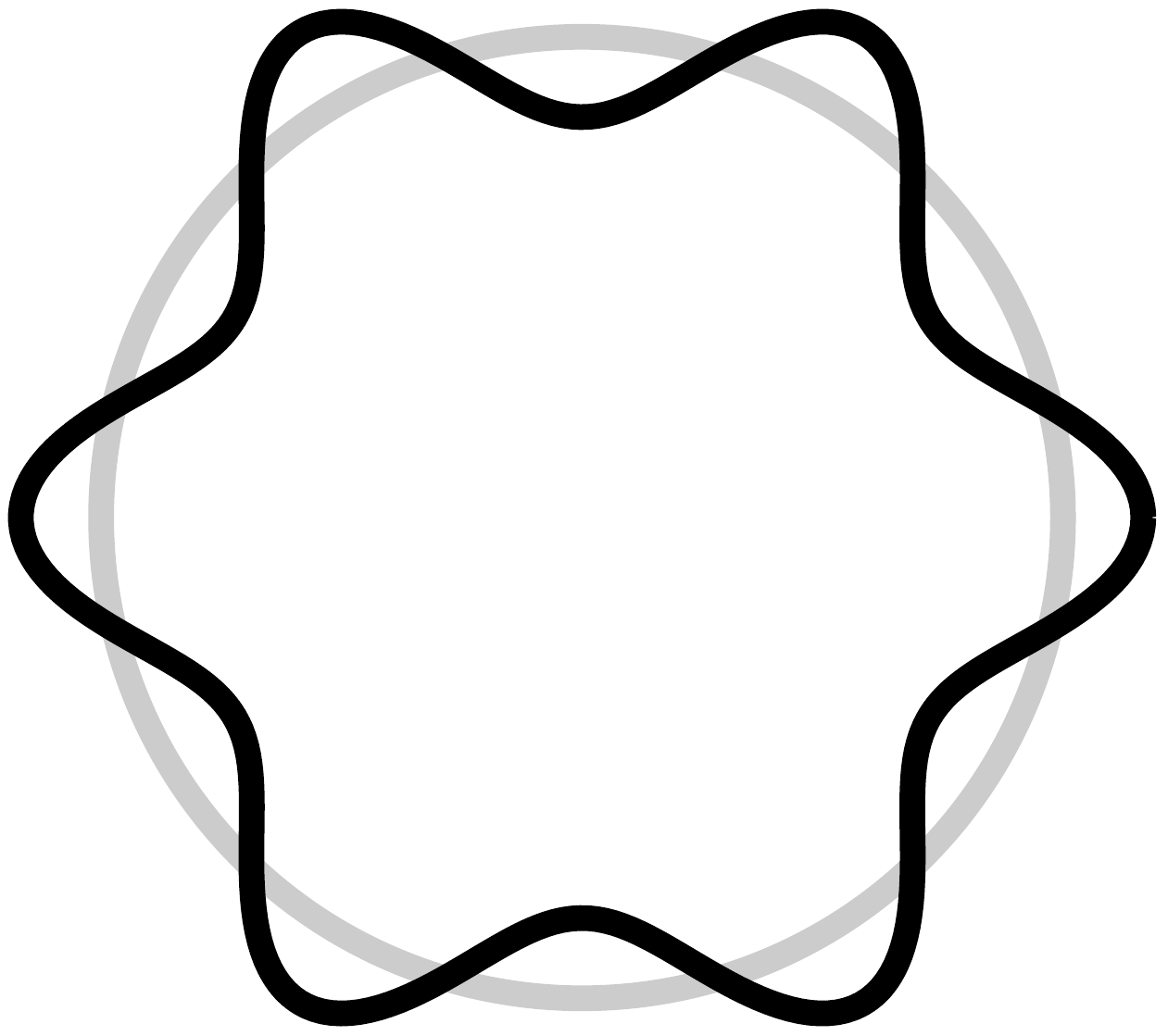}}
	\caption{The first eight analytical bending-dominated eigenfunctions.}
	\label{fig:ring_first12modes1}
	\centering
	\captionsetup[subfigure]{labelformat=empty}
	\subfloat[$n = 0$]{\includegraphics[width=0.11\textwidth]{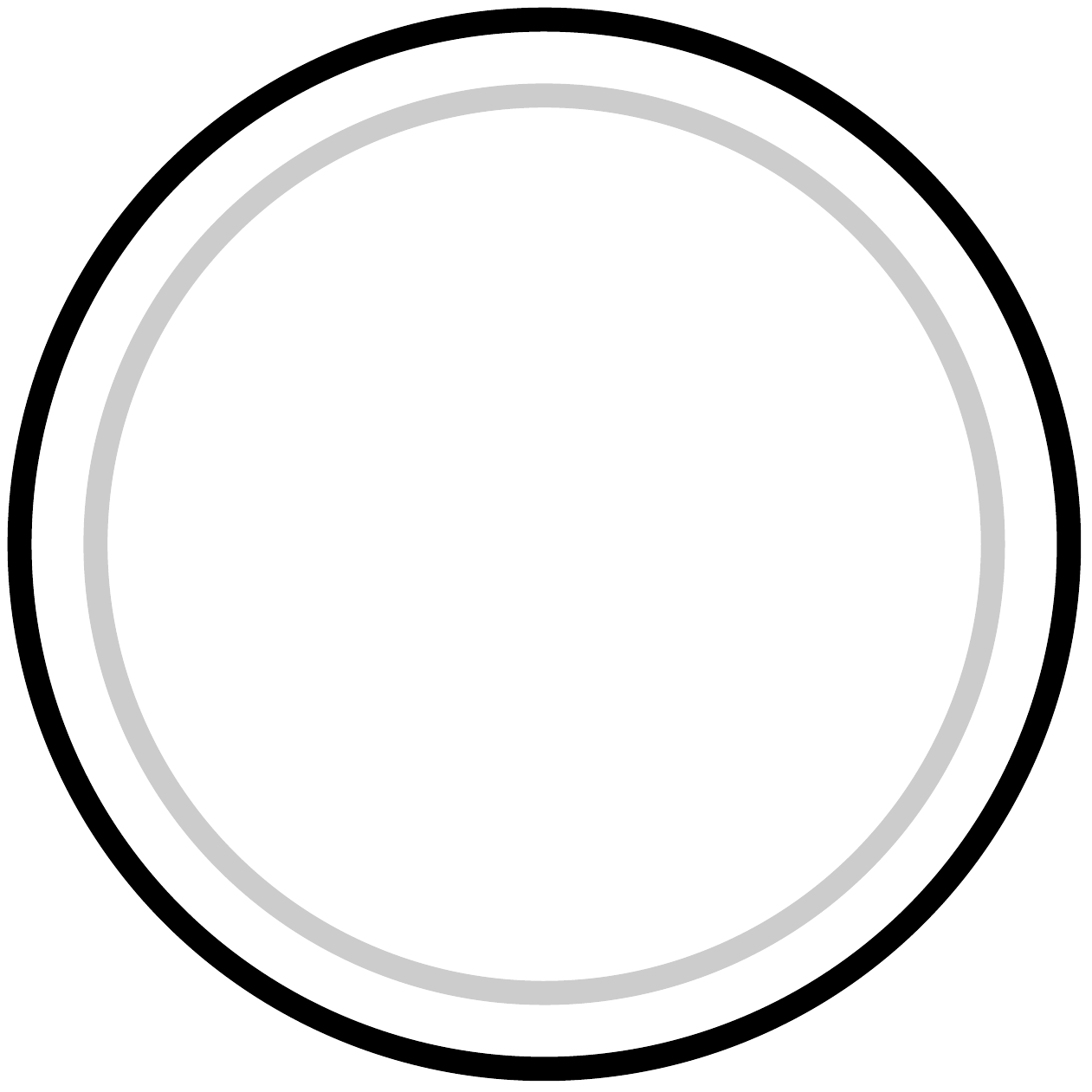}}	\hspace{0.5cm}
	\subfloat[$n = 1$]{\includegraphics[width=0.11\textwidth]{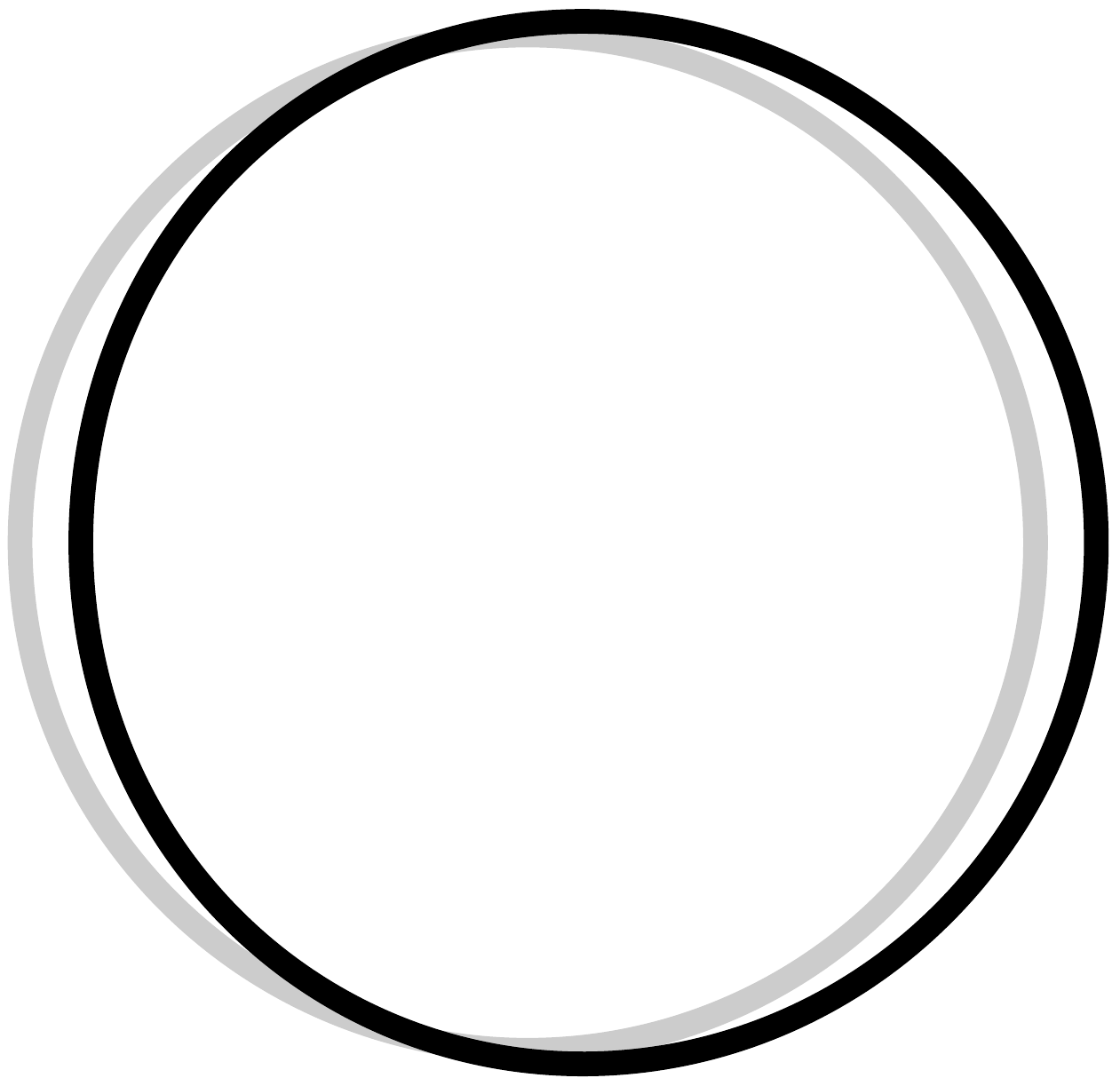}}	\hspace{0.5cm}
	\subfloat[$n = 2$]{\includegraphics[width=0.11\textwidth]{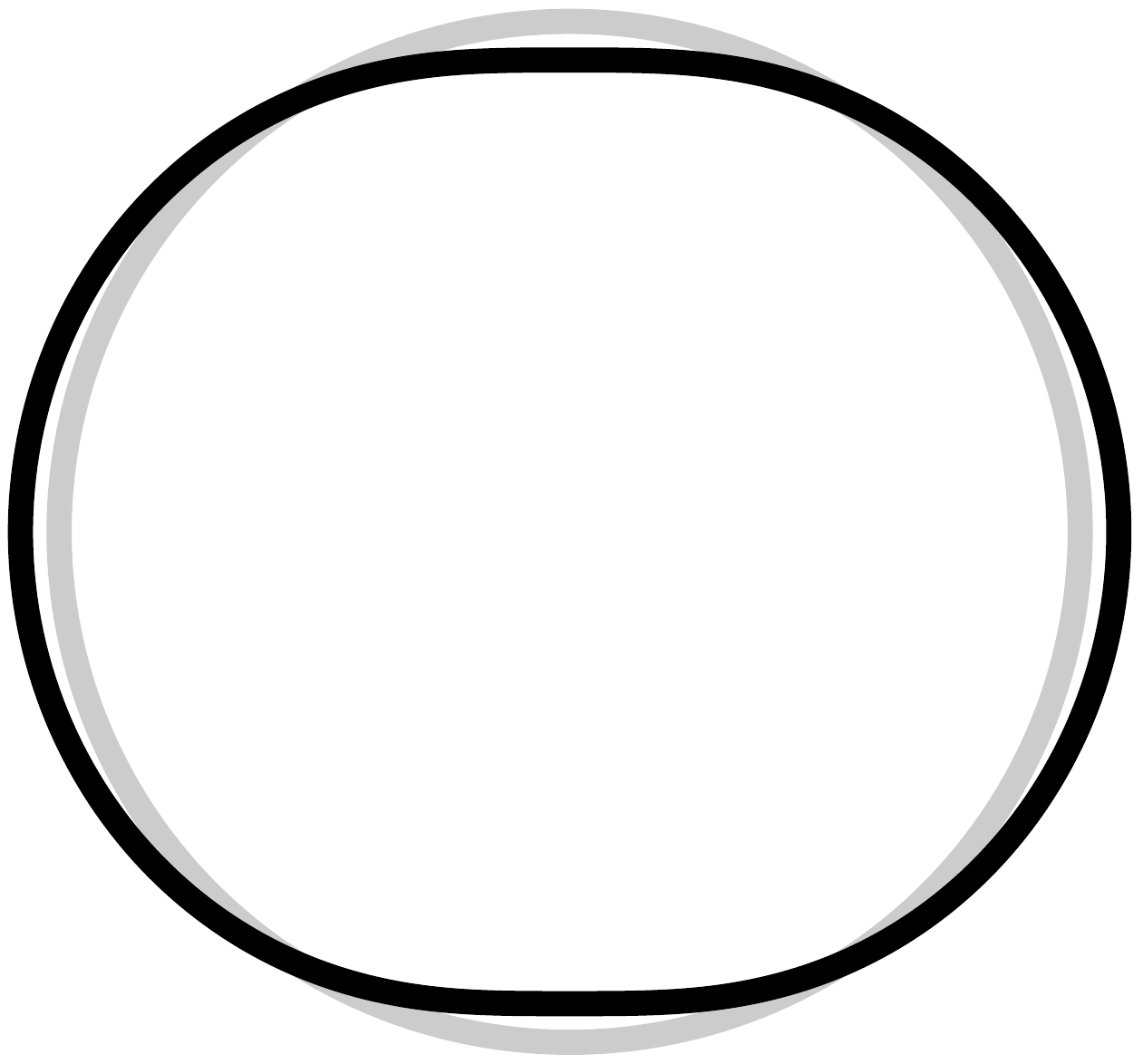}}	\hspace{0.5cm}
	\subfloat[$n = 3$]{\includegraphics[width=0.11\textwidth]{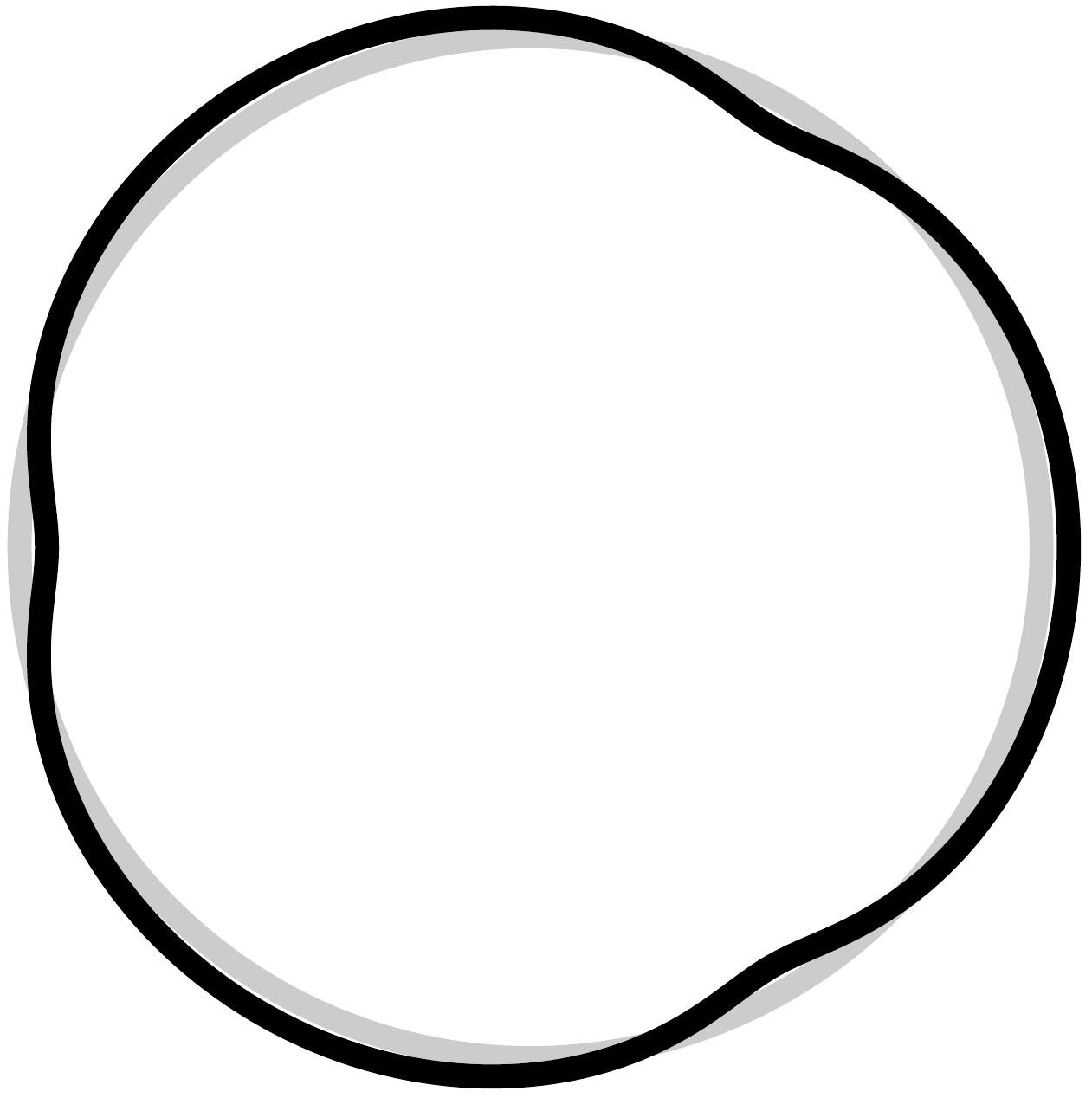}}	\hspace{0.5cm}
	\subfloat[$n = 4$]{\includegraphics[width=0.11\textwidth]{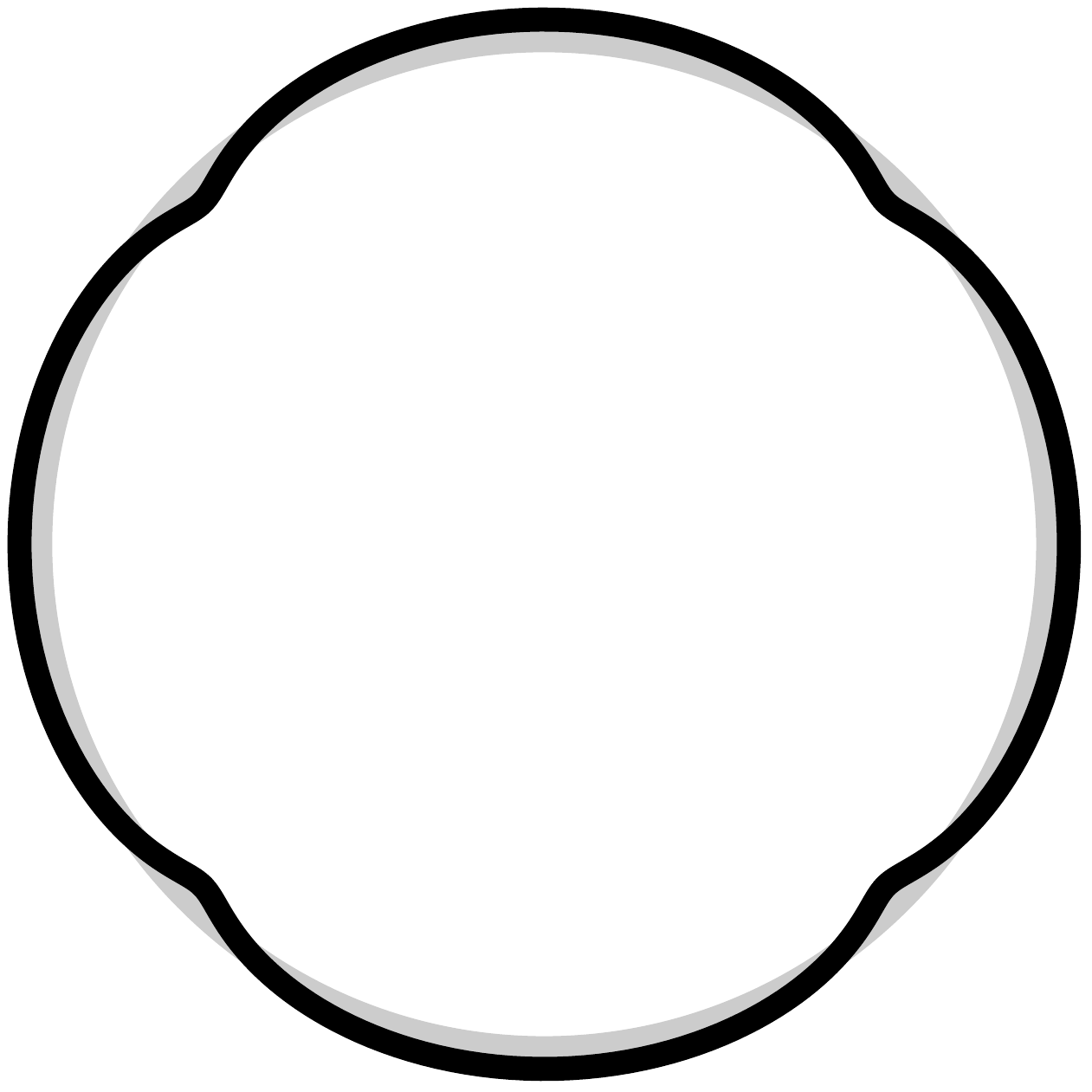}}	\hspace{0.5cm}
	\subfloat[$n = 5$]{\includegraphics[width=0.11\textwidth]{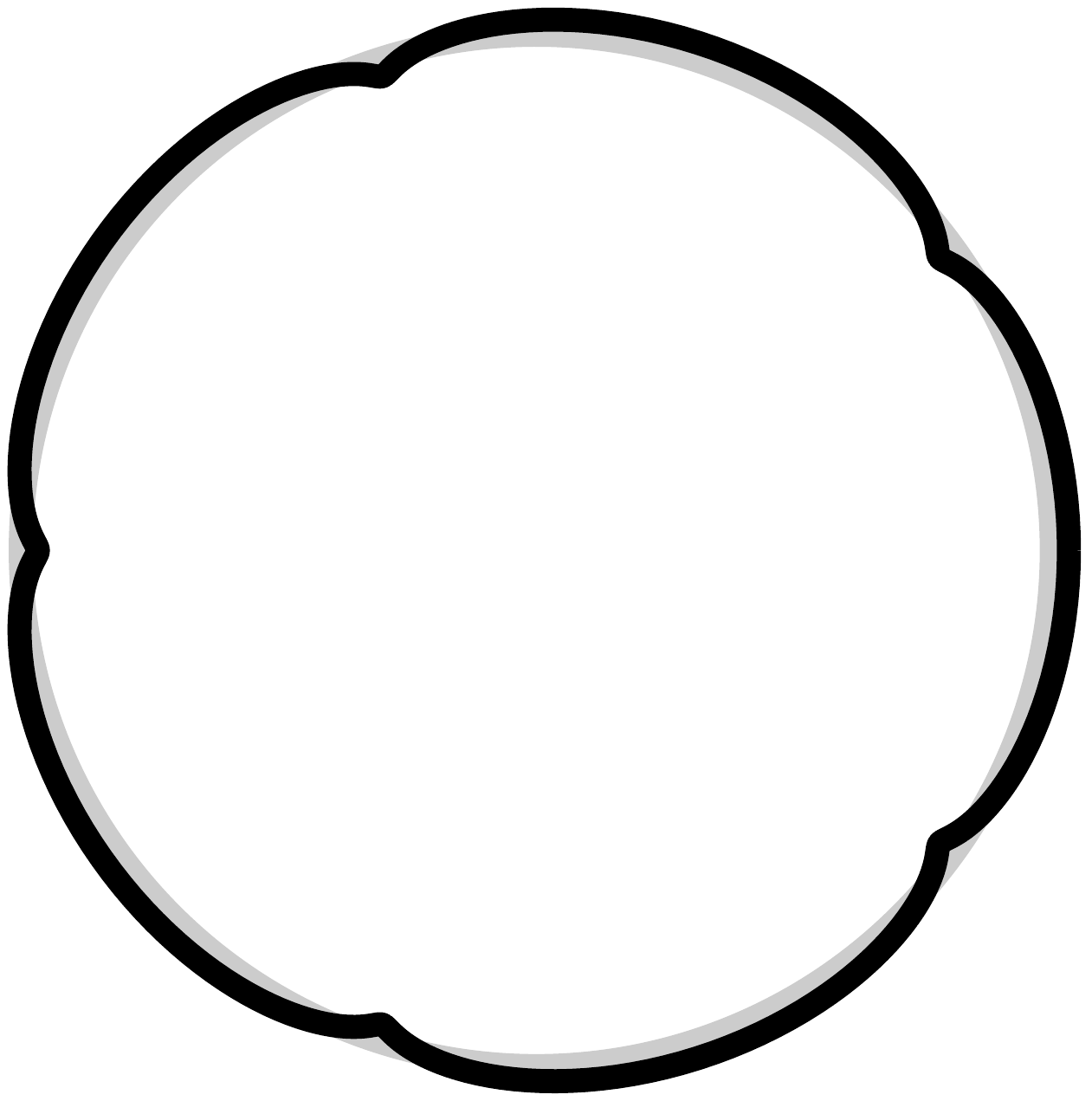}}
	\caption{The first eight analytical membrane-dominated eigenfunctions.}
	\label{fig:ring_first12modes2}
\end{figure}

Let $\phi \in [0, 2\pi)$ denote an arbitrary phase shift. An analytical solution to \eqref{eq:eigenvalue_problem} may be obtained by employing the following informed guess for $\vect{u}_n  = (u_n, \, w_n)$,
\begin{subequations}
\begin{align}
	u_n(\theta) 		&= U_n \sin{\left(n(\theta - \phi) \right)}, \\
	w_n(\theta) 		&= W_n \cos{\left(n(\theta - \phi) \right)}.
\end{align}
\end{subequations}
The eigenvalues and amplitudes of the eigenfunctions satisfy, for each $n$, a two by two matrix eigenvalue problem  (see \ref{app:fourier_series} or \cite[Chapter 5, Section 5.3, page 82]{soedel_vibrations_2004})
\begin{align}
\begin{pmatrix} 
	\tilde{A}_n - \lambda_n 	& \tilde{B}_n  \\ 
	\tilde{B}_n^*	& \tilde{D}_n - \lambda_n  
\end{pmatrix}
\begin{pmatrix} 
	U_n  \\ 
	W_n  
\end{pmatrix}
= 
\begin{pmatrix} 
	0  \\ 
	0  
\end{pmatrix}.
\label{eq:two_by_two_problem}
\end{align}
The values $\tilde{A}_n$, $\tilde{B}_n$ and $\tilde{D}_n$ are listed in \ref{app:fourier_series} in terms of the beam thickness $t$ and the radius $R$. For each $n$, there are two eigenfunctions and corresponding eigenvalues. The eigenvalues are computed as the roots of the characteristic equation
\begin{subequations}
\label{eq:lambaanalyticKL}
\begin{align}
	\lambda_{n+} 	&= K_n \left(1 + L_n \right)			\\
	\lambda_{n-} 		&= K_n \left(1 - L_n \right)
\end{align}
\end{subequations}
where $K_n$ and $L_n$ are constants that depend on $n$ , the normalized thickness $\bar{t} = t /R$ and radius $R$, and are given by
\begin{align}\label{eq:exactKML}	
	K_n = \frac{\left(n^2+1\right)\,\left(n^2\, \bar{t}^2+12\right)}{24 R^2}, \qquad
	M_n  = \frac{n^2\,\bar{t}^2\,{\left(n^2-1\right)}^2}{12 R^4}, \qquad
	L_n = \sqrt{1 - \frac{M_n}{K_n^2}}.
\end{align}
The eigenfunctions are defined up to a constant. However, the ratio of the amplitudes of the eigenfunctions, $U_n / W_n$, satisfy
\begin{subequations}
\begin{align}
	\frac{U_{n+}}{W_{n+}} = \frac{\tilde{B}_n}{\lambda_{n+} - \tilde{A}_n}			\\
	\frac{U_{n-}}{W_{n-}} = \frac{\tilde{B}_n}{\lambda_{n-} - \tilde{A}_n}
\end{align}
\end{subequations}
When $U_n / W_n < 1$ the transverse displacement is larger than the circumferential displacement, thus the bending response dominates. When $U_n / W_n > 1$, the circumferential displacements are larger than the transverse displacements and, consequently, the membrane response dominates. In what follows we will refer to these two types of modes as bending-dominated and membrane-dominated. Figure \ref{fig:ring_first12modes1} and \ref{fig:ring_first12modes2} depict several examples of the two different types of modes. 

\begin{figure}[H]
	\centering
	\subfloat[$\bar{t} = t/R = 0.1$]{\includegraphics[width=0.48\textwidth]{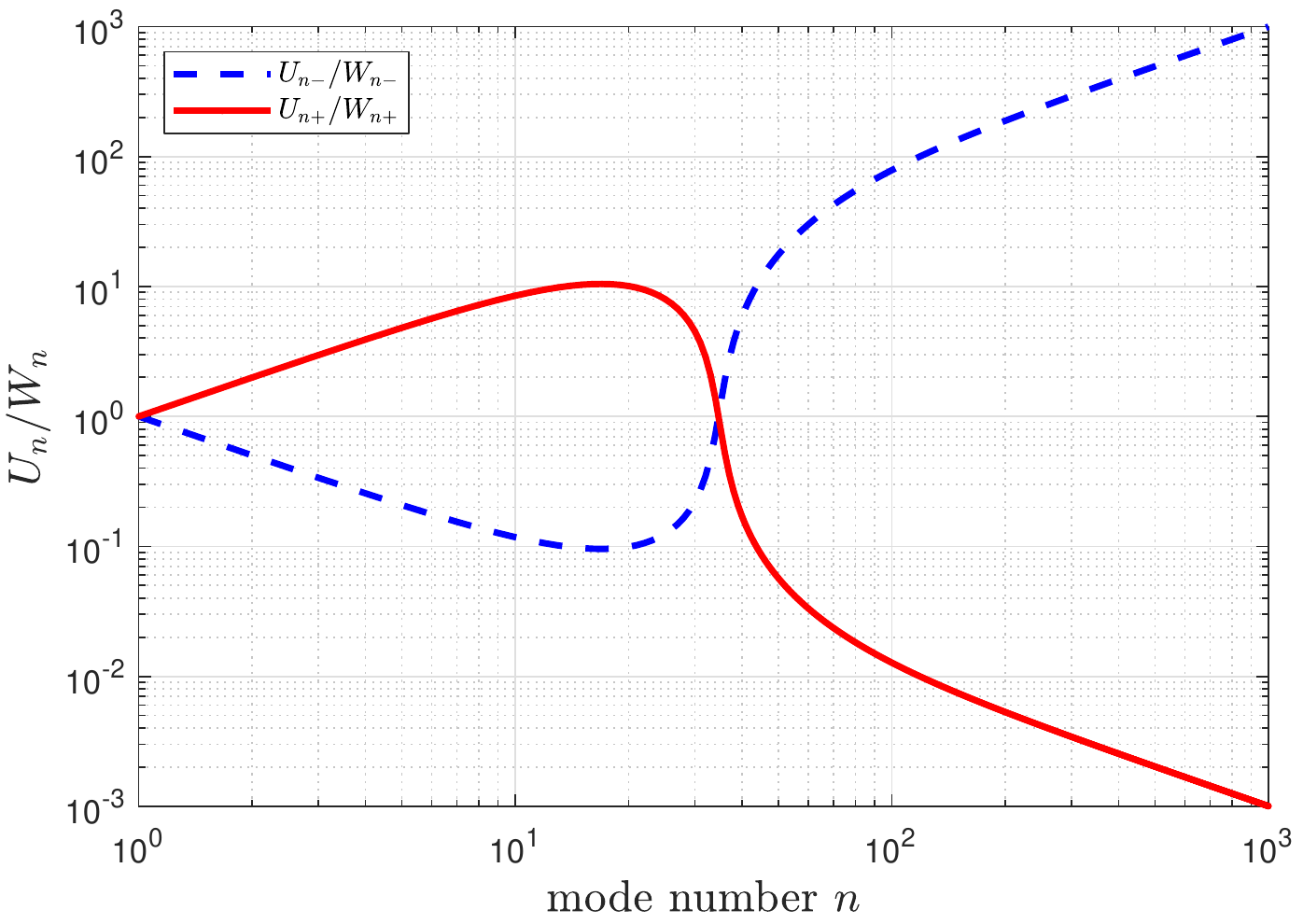}} 
	\subfloat[$\bar{t} = t/R = 0.01$]{\includegraphics[width=0.48\textwidth]{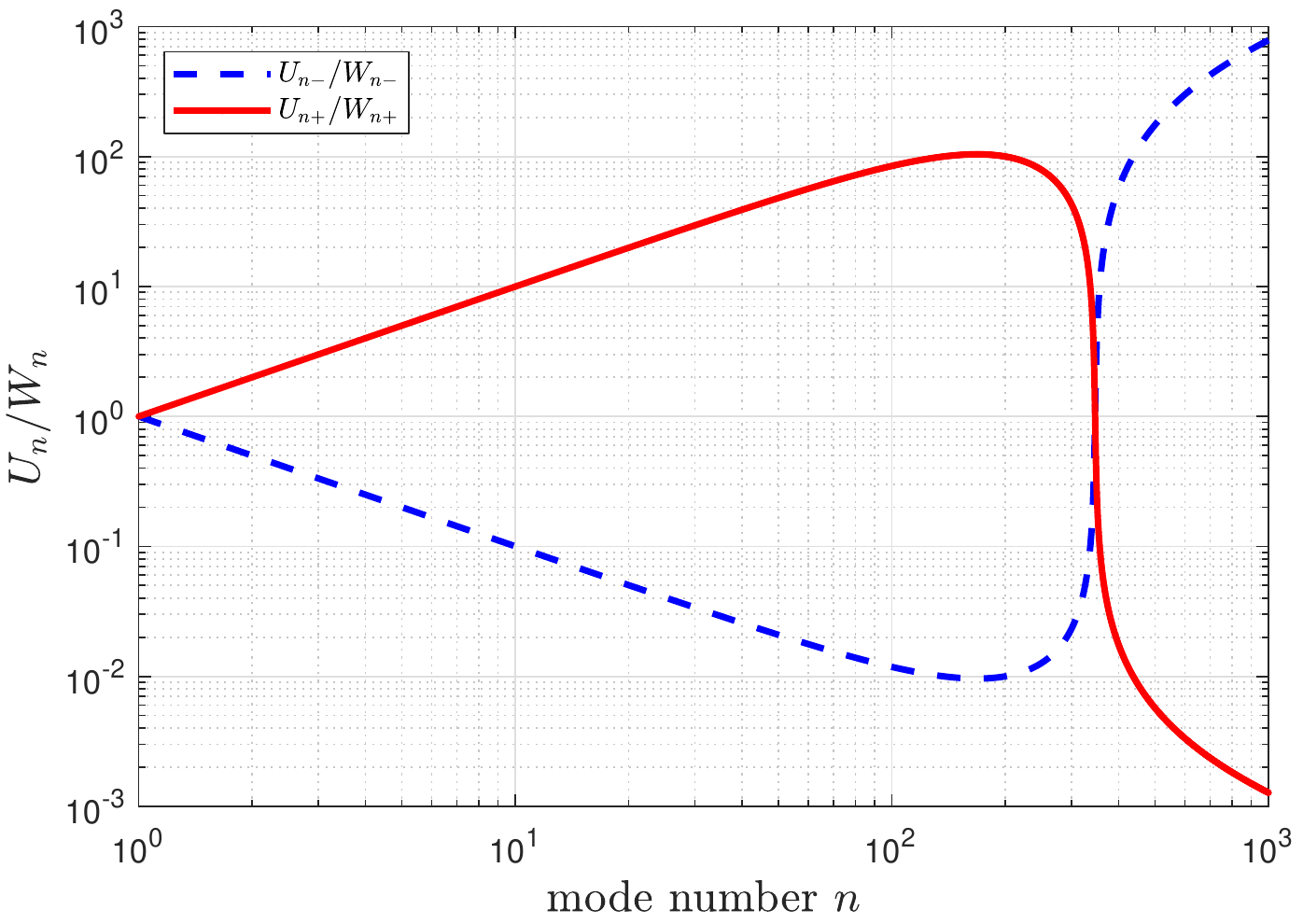}}
	\caption{For $n < \hat{n}$ the modes corresponding to $\lambda_{n+}$ are membrane-dominated since $U_{n+}/W_{n+} > 1$ and the modes corresponding to $\lambda_{n-}$ are bending-dominated since $U_{n-} / W_{n-}<1$. For $n > \hat{n}$ this behavior is reversed. This is shown for two different normalized thicknesses $\bar{t}$, since $\hat{n}$ depends on $\bar{t}$.}
	\label{fig:analytical_eigenmodes}
\end{figure}

Figure \ref{fig:analytical_eigenmodes} provides a more detailed view of the behavior of the ratio of the amplitudes, $U_n / W_n$, at two different choices of the normalized thickness. Modes $\vect{u}_{n-}$, with a ratio of amplitudes  $U_{n-} / W_{n-}$ and corresponding eigenvalue $\lambda_{n-}$, have dominating transverse displacements for $n<\hat{n}$, where $\hat{n}$ is a value that depends on the normalized thickness $\bar{t}$. In contrast, modes $\vect{u}_{n+}$, with a ratio of amplitudes  $U_{n+} / W_{n+}$ and corresponding eigenvalue $\lambda_{n+}$, have dominating circumferential displacements for $n<\hat{n}$. When $n>\hat{n}$ this behavior reverses.

\begin{figure}[!h]
	\centering
	\subfloat[$\bar{t} = t/R = 0.1$]{\includegraphics[width=0.48\textwidth]{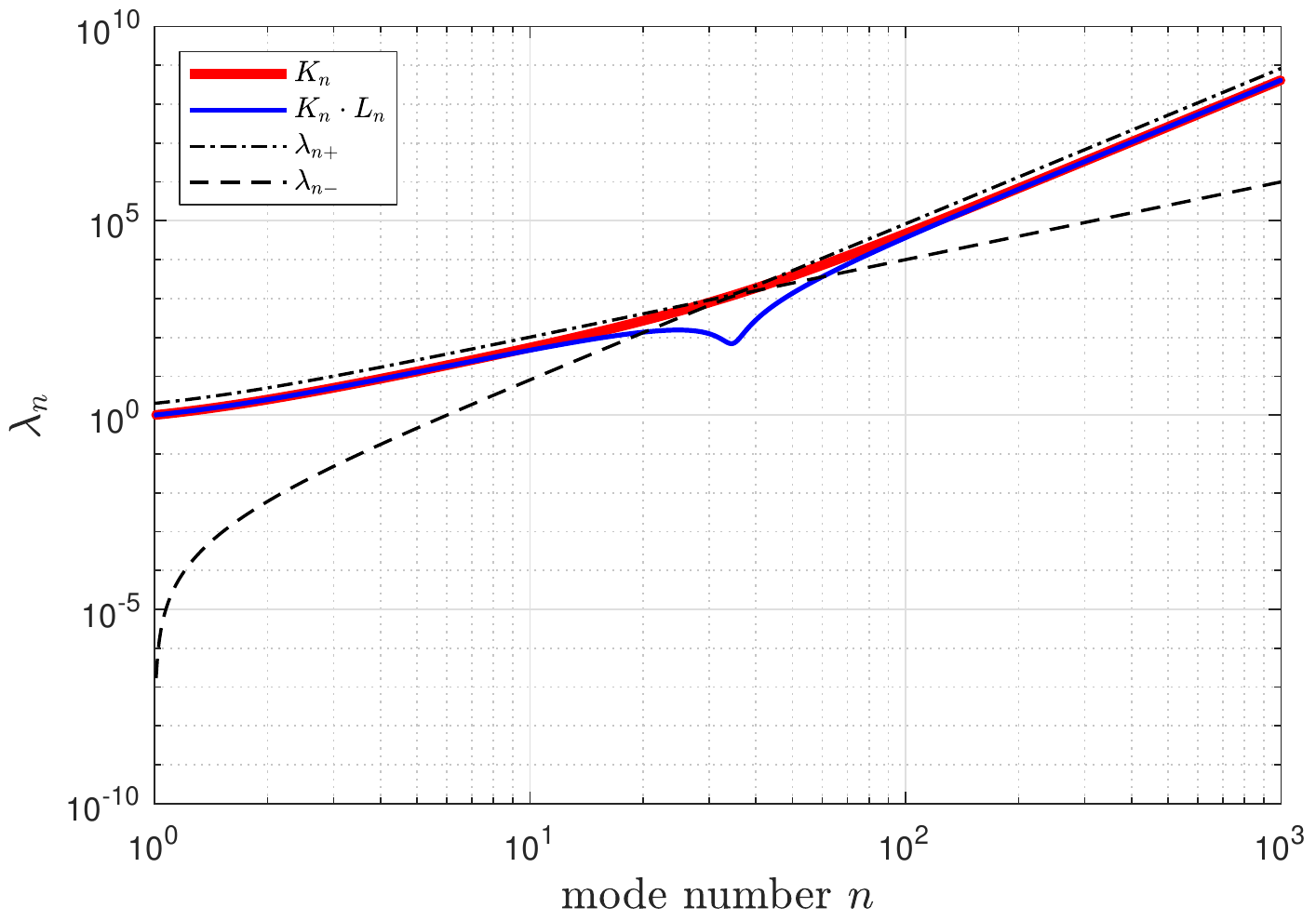}} 
	\subfloat[$\bar{t} = t/R = 0.01$]{\includegraphics[width=0.48\textwidth]{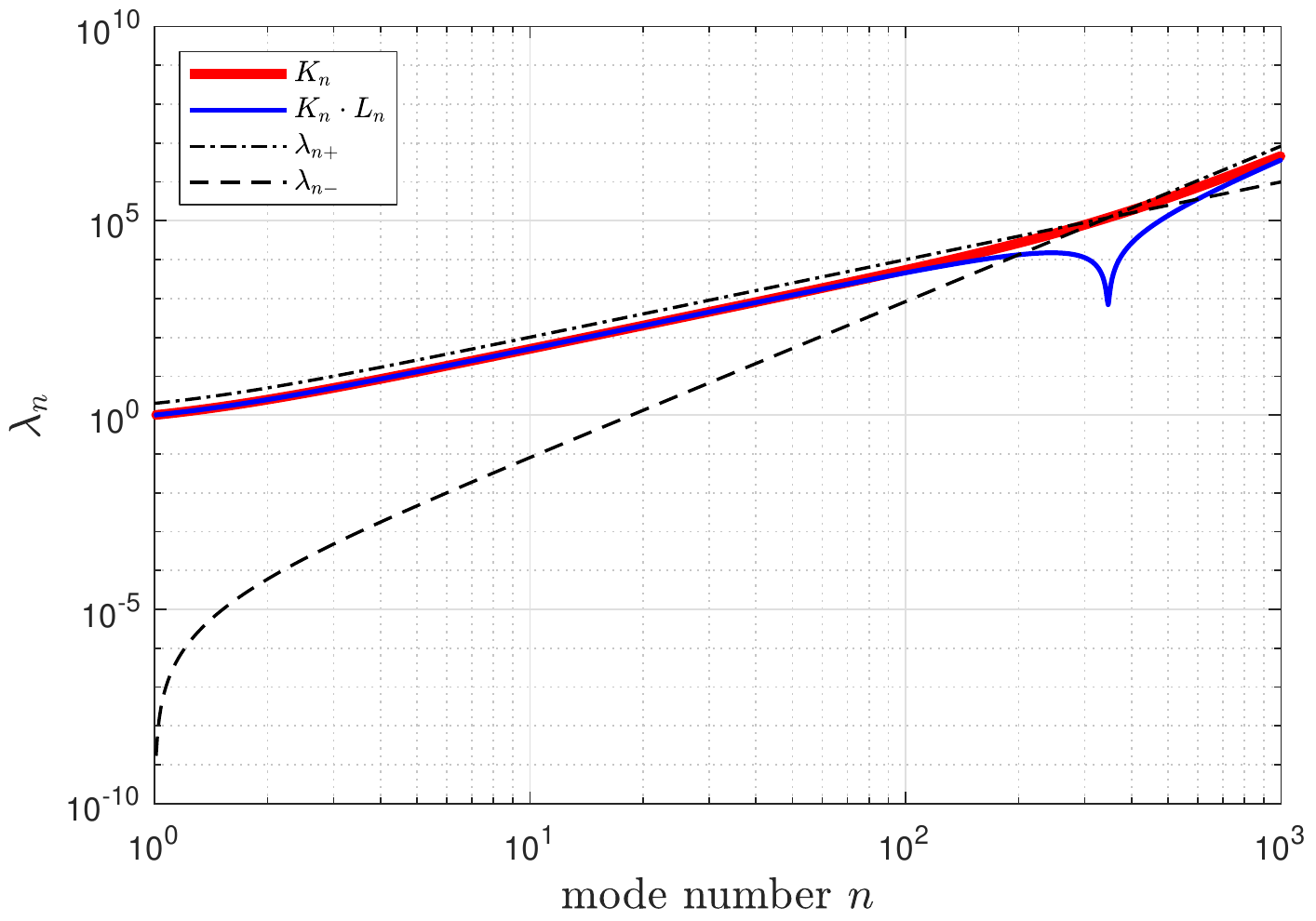}}
	\caption{The analytical spectrum plotted on a log-log scale. Typically $\lambda_{n+} >> \lambda_{n-}$ except near $n = \hat{n}$, which is the location of the local minimum of the function $K_n \cdot L_n$, denoted in blue.}
	\label{fig:analytical_eigenvalues}
\end{figure}

\begin{figure}[!h]
	\centering
	\includegraphics[width=0.75\textwidth]{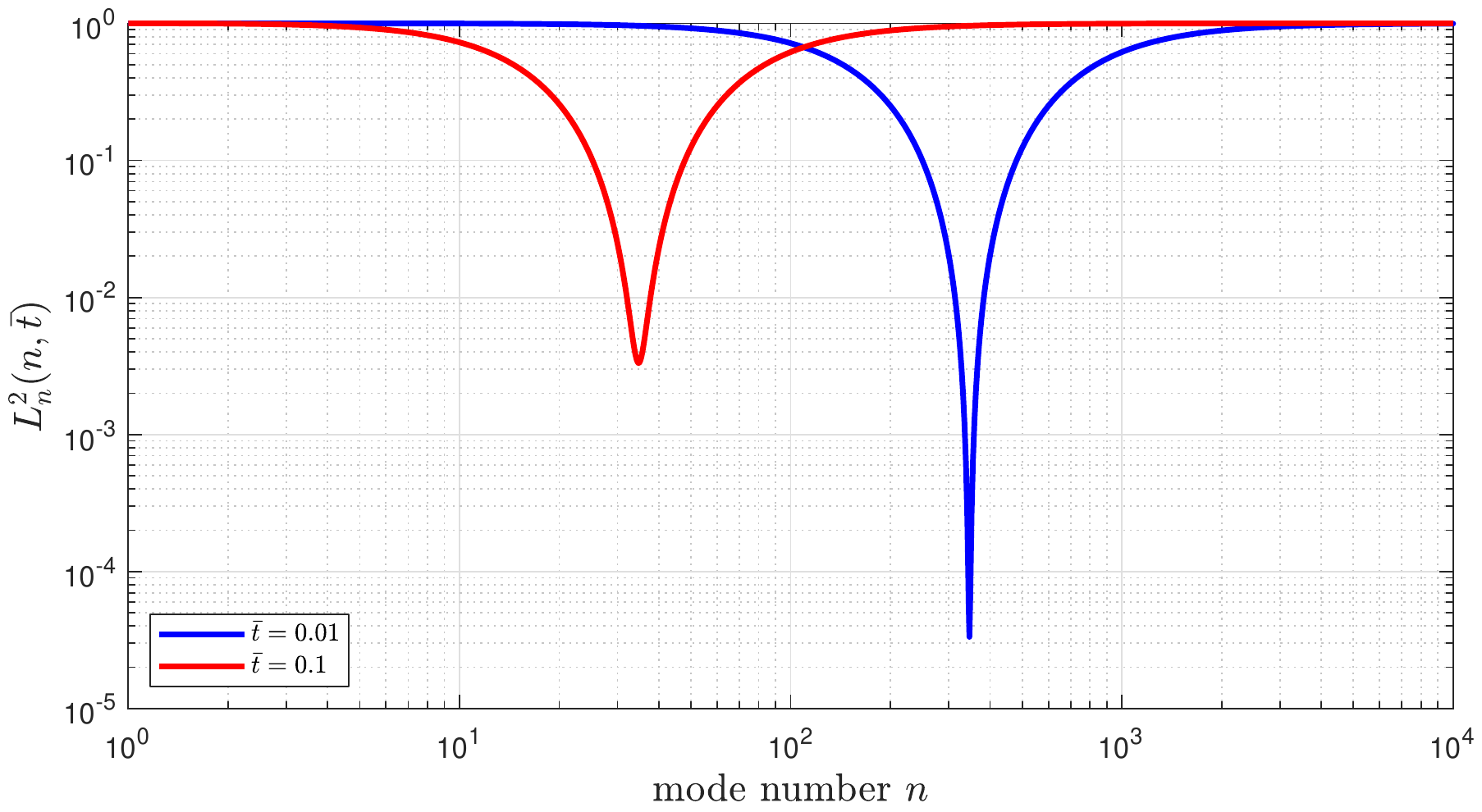}
	\caption{The function $L_n^2(n,\bar{t}) = 1 - M_n / K^2_n$ is a smooth positive function that is less than or equal to $1$. For $n\geq1$, the function attains a single minimum at mode number $\hat{n}$, which depends on the value of $\bar{t} = t / R$.}
	\label{fig:minimal_point}
\end{figure}

The corresponding branches of eigenvalues, $\lambda_{n+}$ and $\lambda_{n-}$, also behave differently, as may be observed in Figure \ref{fig:analytical_eigenvalues}, where they are plotted as a continuous (rather than discrete) function of $n$. For reasonable values of the material and cross sectional parameters it follows that $\lambda_{n+} >> \lambda_{n-}$, except near $n = \hat{n}$, which is the location of the local minimum of the function $K_n \cdot L_n$, shown in blue in Figure \ref{fig:analytical_eigenvalues}. At this point the argument inside the square root, $L^2_n = 1 - M_n / K_n^2$, attains a global minimum, which itself approaches zero when viewed as a function of $\bar{t}$ with $\bar{t} \rightarrow 0$, see Figure \ref{fig:minimal_point}.

\begin{figure}[!h]
	\centering
	\includegraphics[width=0.6\textwidth]{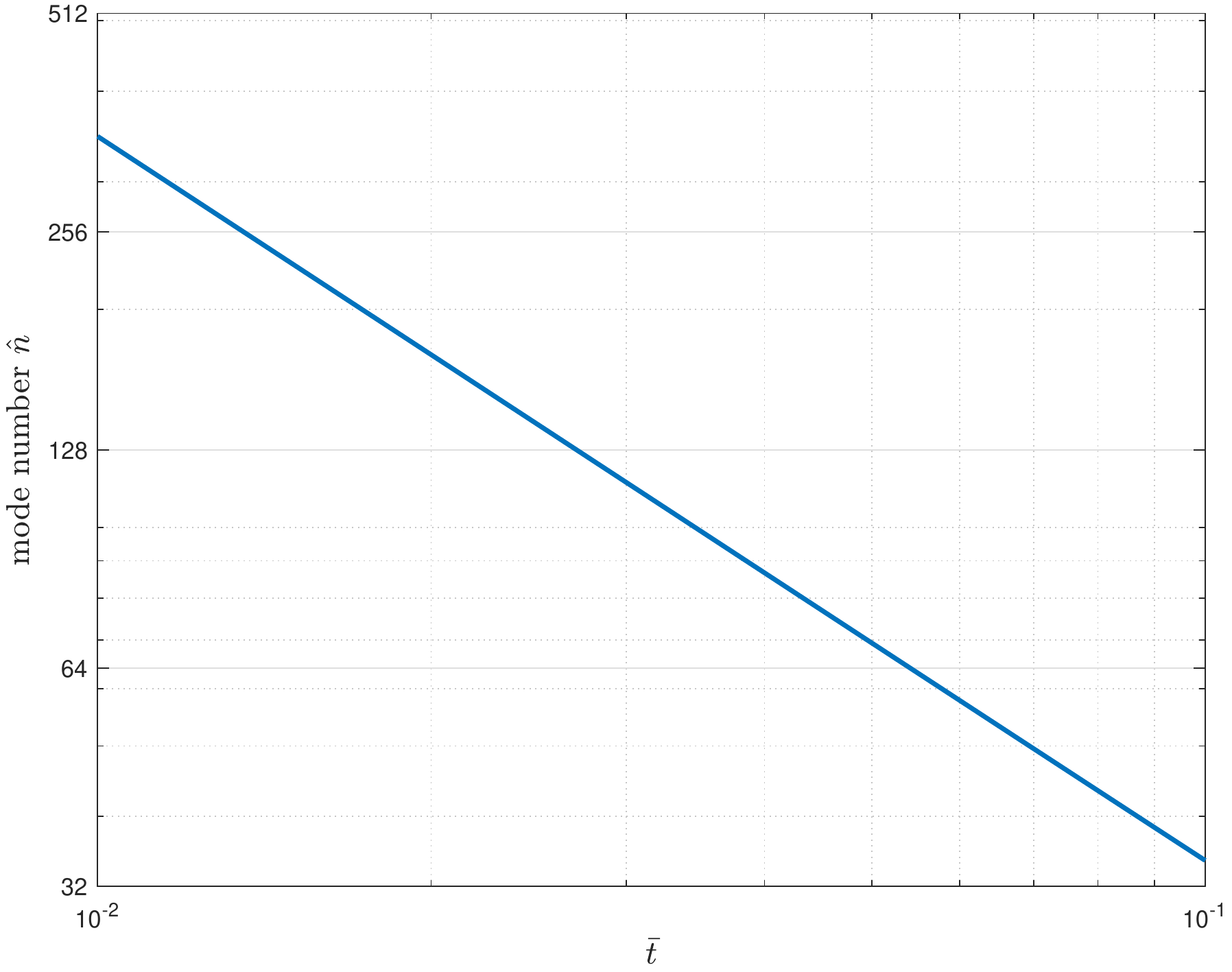}
	\caption{The position of $\hat{n}\geq1$, where $L^2_n(n,\bar{t})$ achieves a minimum, as a function of the normalized thickness $\bar{t} = t/R$.}
	\label{fig:minimal_point2}
\end{figure}

More precisely, $L_n^2$ depends on $n$ and on the normalized thickness $\bar{t}$, but not on the radius, 
\begin{align}	
	L^2_n = 1 - 
\frac{48n^2 \, \bar{t}^2 \, (n^2 - 1)^2}{(n^2 + 1)^2  (n^2 \, \bar{t}^2 + 12)^2}
\label{eq:Ln}
\end{align}
Thus, the mode number at which the minimum of $L_n^2$ is attained, denoted by $\hat{n}$, can be determined as a function of the normalized thickness $\bar{t} = t/R$. The minimum coincides with the maximum of the fraction in \eqref{eq:Ln}. Let $f$ denote its numerator and $g$ the denominator. At an optimum we require that $f^\prime(n) g(n) - f(n) g^\prime(n) = 0$, or, since $g(n)>0$ for any $n\geq0$, $f^\prime(n) - f(n) g^\prime(n) / g(n) = 0$, meaning that we must seek the zeros of the function
\begin{align*}
z(n) &= f^\prime(n) - f(n) g^\prime(n) / g(n)  \\ &= 
\frac{96\,n\,\bar{t}^2\,\left(n^2-1\right)\,\left(-\bar{t}^2 \, n^6 + (4\bar{t}^2+12) n^4 + (\bar{t}^2+48) n^2 -12 \right)}{(n^2+1)(n^2 \, \bar{t}^2 + 12)}\,.
\end{align*}
By inspection, the roots of $z(n)$ are $n=0$, $n = \pm 1$, and the roots of the sixth order polynomial in the numerator. Focusing on the latter we substitute $x = n^2$, and obtain a cubic polynomial, $h(x) = -\bar{t}^2\,x^3 + (4\,\bar{t}^2 + 12) x^2  + (\bar{t}^2+48) x - 12$, whose only root above unity, $n^2=x>1$, is given by Cardano's formula, and corresponds precisely to $\hat{n}$ (i.e., the minimum of $L_n$ for $n\geq1$). We do not show the explicit expression for $\hat{n}$, as it is very lengthy, but the main point is that it is an algebraic expression depending solely on the normalized thickness $\bar{t}$. Figure \ref{fig:minimal_point2} plots the value of $\hat{n}$ as a function of $\bar{t}$ in a log-scale, and we note the function is very well approximated (with negligible mean squared error in a broad range of $\bar{t}$) by a linear function in log-space given by
\begin{align}\label{eq:nhatformula}
	\log_{10}\big(\hat{n}(\bar{t})\big) \approx -0.9985 \log_{10}(\bar{t}) + 0.5421
\end{align}

\section{Discretization by periodic uniform splines \label{sec:3}}
In this section we introduce the discrete function spaces used to discretize the standard and mixed eigenvalue problem presented in the previous section. We consider as trial space a maximally smooth periodic space of splines defined on a uniform partition of the circle. In this special case there exists a discrete Fourier basis that inherits many of the properties of analytical Fourier modes (see Figure \ref{fig:overview}), allowing us analytically compute expressions for discrete eigenvalues and the ratio of amplitudes of the modes.

\subsection{Periodic uniform splines}
This section introduces Cardinal and uniform splines and discusses how they are used to develop a discrete trial and test space of splines.

\subsubsection*{Cardinal B-spline function}
\begin{figure}[H]
	\centering	   
	\begin{tikzpicture}
		\node[inner sep=0pt] (russell) at (0,0)
    {\includegraphics[trim = 0cm 0cm 0cm 0cm, clip,width=0.5\textwidth]{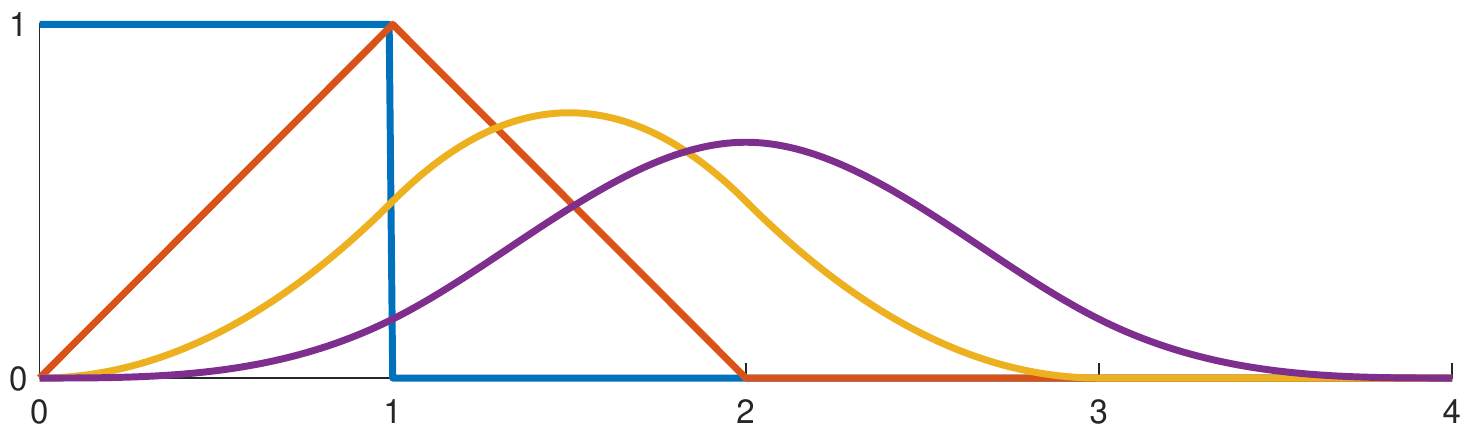}};
		\path (-4,1.75) node {$\phi_{0}$};
		\path (-2,1.4) node {$\phi_{1}$};
		\path (-0.85,1) node {$\phi_{2}$};
		\path (1,0.6) node {$\phi_{3}$};
	\end{tikzpicture}
	\caption{Convolution by the characteristic function increases smoothness and polynomial degree.}
	\label{fig:uniform_spline}
\end{figure}
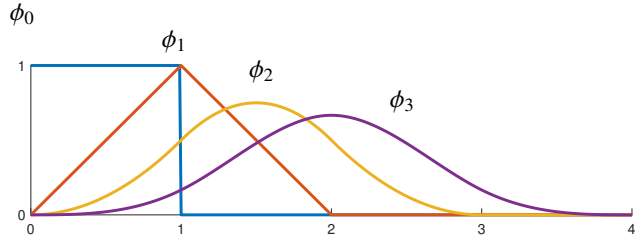

A convenient basis in which to represent polynomial splines is given by B-splines. An important special case is the so called Cardinal B-spline.
\begin{definition}[Cardinal B-spline] \label{def:uniform_bspline} Let $\vect{1}_{[0,1)}$ denote the characteristic function on $[0,1)$. The Cardinal B-spline of polynomial degree $\p$ is defined via the convolution
\begin{align}
	\phi_{\p} := \phi_{\p-1} \ast \vect{1}_{[0,1)}
\end{align}
\end{definition}
A graphical interpretation of Definition \ref{def:uniform_bspline} is depicted in Figure \ref{fig:uniform_spline}. Convolution increases the polynomial degree as well as the smoothness. Several useful properties follow directly. Here, $\phi_\p (x)$ is a degree $\p$ polynomial on every element $(k, k+1)$, is $C^{\p-1}$ smooth, positive inside the interval $(0,p+1)$ and zero elsewhere.

\subsubsection*{Uniform splines}
Consider a partitioning of the real line with mesh-size $h$
\begin{align}
	h \mathbb{Z} := \ldots -2h, \; -h, \; 0, \; h, \; 2h, \; \ldots
\end{align}
A space of uniform B-splines on $h \mathbb{Z}$ is obtained by translated and scaled copies of $\phi_{\p}$.
\begin{definition}[Uniform B-splines] For $h > 0$
\begin{align}
	B_{i,p}(x) = \phi_p(x/h - i), \qquad i \in \mathbb{Z}
\end{align}
is the $i$-th B-spline defined on the grid $h \mathbb{Z}$. Linear combinations $\sum_{i \in \mathbb{Z}} \alpha_i \, B_{i,p}(x)$ are called uniform splines of polynomial degree $\p$ and mesh-size $h$.
\end{definition}

\subsubsection*{Periodic uniform splines on the unit ring}
\begin{figure}[h]
	\centering
	\includegraphics[trim = 25cm 3cm 0cm 0cm, clip,width=0.3\textwidth]{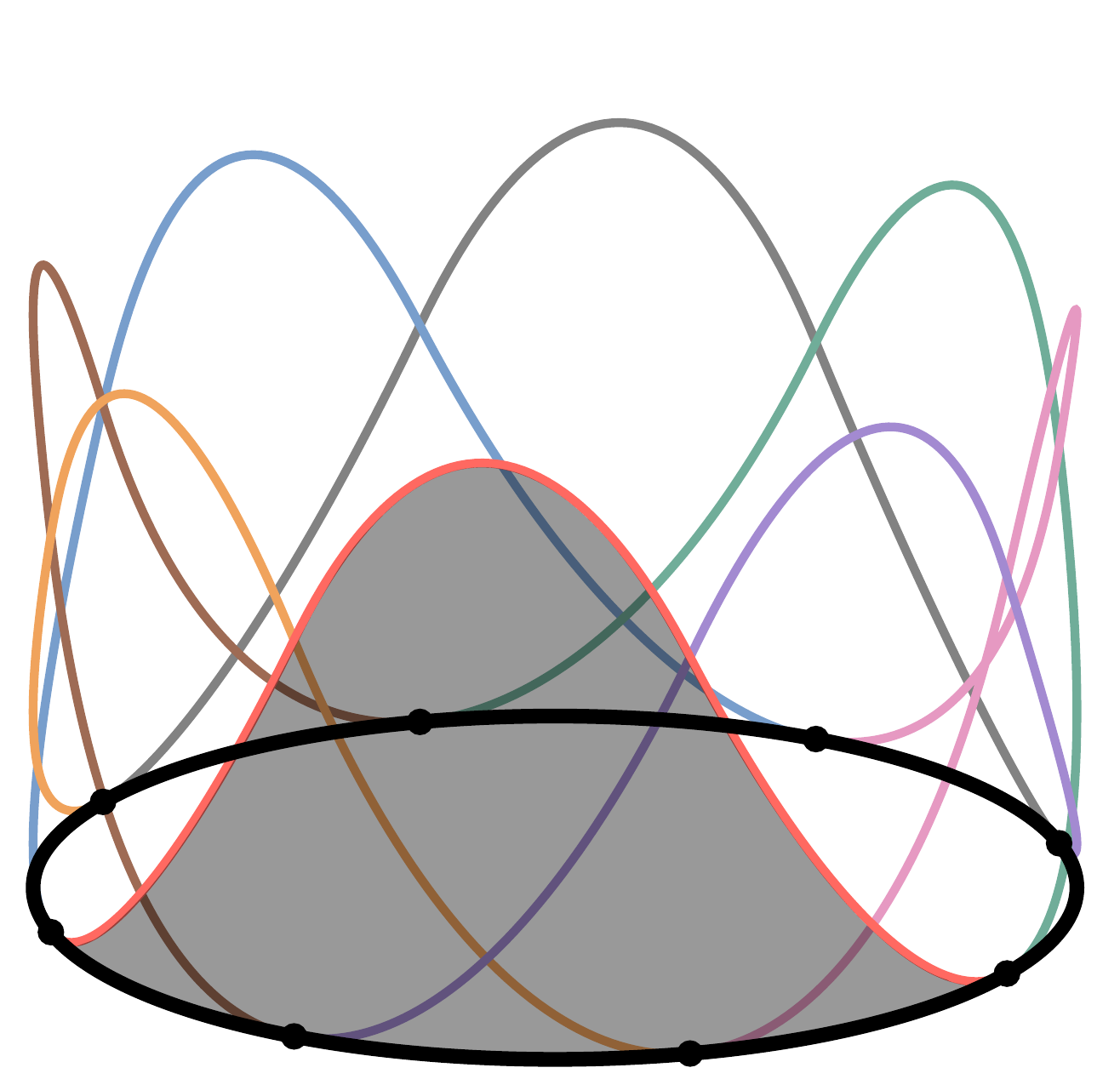}	   
	\caption{The periodic space $\sspace{2}{8}(0, 2 \pi)$, consisting of $8$ quadratic uniform B-splines on the circle}
	\label{fig:periodic_splines_circle}
\end{figure}

Periodic spaces of splines of dimension $N$ can be constructed by taking $N+\p$ sequential B-splines and applying suitable end-conditions to the last $\p$ B-splines. In Definition \ref{def:splinespace} this is done for a smooth uniform space of splines defined on the unit circle. Figure \ref{fig:periodic_splines_circle} shows a graphical interpretation.
\begin{definition}[Periodic uniform spline space \label{def:splinespace}] \label{def:} Let $h = 2\pi / N$ denote the mesh size and $\theta \in [0,2 \pi)$ denote the angular coordinate. The periodic $N$-dimensional space of splines on the unit-circle is defined as
\begin{align}
	\sspace{\p}{N}(0, 2 \pi) = \left\{ s \; : \; [0, 2 \pi) \mapsto \mathbb{R} \; : \; s(\theta) = \sum_{i =-\p}^{N-1} \alpha_i \, B_{i,p}(\theta), \; \text{with } \alpha_{-\p} = \alpha_{N-\p}, \ldots ,  \alpha_{-1} = \alpha_{N-1} \right\}
\end{align}
\end{definition}
Uniform, periodic B-splines inherit local support, positivity, and are globally $C^{p-1}$ smooth. The periodic and translation invariant properties of the basis lead to circulant system matrices when used in weak forms. Such matrices are easily diagonalized via the Discrete Fourier Transform (DFT), see \ref{sec:EVcirculant}. Since we are dealing with piecewise polynomials, differentiation and integration can be performed analytically using the standard spline formulae, see \cite{boor_practical_2001}. 

\subsection{The mixed formulation}
We consider discrete function spaces of the form
\begin{subequations}
\begin{align}
	&V^h = \mathbb{S}^{p}_N(\domain) \times \mathbb{S}^{p}_N(\domain), &&	
	S^h = \mathbb{S}^{p-1}_N(\domain) \times \mathbb{S}^{p-1}_N(\domain).&
\end{align}
\end{subequations}
For $p>1$ it holds that $V^h$ and $S^h$ are conforming subspaces of $V$ and $S$, respectively. The Galerkin mixed formulation reads: find $\vect{u}^h = (u^h,w^h) \in V^h$, $\vect{\varepsilon}^h = (\varepsilon^h,\chi^h) \in S^h$ and $\lambda^h \in \mathbb{R}^{+}_{0}$ such that
\begin{subequations}
\label{eq:eigenvalue_problem_discrete}
\begin{align}
	- \cc{\vect{\varepsilon}^h}{\vect{\delta \varepsilon}^h} + \cc{B(\vect{u}^h)}{\vect{\delta \varepsilon}^h} &= 0 & \forall \vect{\delta \varepsilon}^h &= (\delta \varepsilon^h, \delta \chi^h) \in S^h \label{eq:eigenvalue_problem_kinematics_discrete} \\
	  \cc{\vect{\varepsilon}^h}{B(\vect{\delta u}^h)} &= \lambda^h \, \bb{\vect{u}^h}{\vect{\delta u}^h} 	 & \forall \vect{\delta u}^h &= (\delta u^h, \delta w^h) \in V^h. \label{eq:eigenvalue_problem_equation_discrete}
\end{align}
\end{subequations}
The mixed Galerkin formulation \eqref{eq:eigenvalue_problem_discrete} results in a $4N \times 4N$ system of equations. As shown in \ref{app:splines}, this system can be reduced to a $2N \times 2N$ system in standard form using static condensation.

\subsection{The standard formulation}
As usual, the Galerkin method restricts the formulation to the finite dimensional setting as follows: find $\vect{u}^h = (u^h,w^h) \in V^h$ and $\lambda^h \in \mathbb{R}^{+}_{0}$ such that
\begin{align}
	\aa{\vect{u}^h}{\vect{\delta u}^h} = \lambda^h \, \bb{\vect{u}^h}{\vect{\delta u}^h} \qquad \forall \vect{\delta u}^h = (\delta u^h, \delta w^h) \in V^h.
	\label{eq:standard_formulation_discrete}
\end{align}
The standard formulation leads directly to a $2N \times 2N$ matrix eigenvalue problem, see \ref{app:splines}. Although the standard and mixed formulation are equivalent at the continuous level, they are not in the discrete setting. The mixed Galerkin method offers additional flexibility in the discretization of $S^h$. 

\subsection{Analytical computation of discrete eigenvalues and eigenfunctions}
The matrix system corresponding to the standard and mixed Galerkin discretization are discussed in \ref{app:standard} and \ref{app:mixed}, respectively. The Fourier modes decouple the $2N \times 2N$ system of equations to $N$ systems of $2 \times 2$ equations, as discussed in \ref{app:similarity_transforms} and  \ref{sec:circulent-block-matrix}.  The eigenvalues and amplitudes of the eigenfunctions satisfy, for each $n$, a two by two matrix eigenvalue problem
\begin{align}
\begin{pmatrix} 
	\tilde{A}^h_n - \lambda^h_n 	& \tilde{B}^h_n  \\ 
	(\tilde{B}_n^h)^\ast	& \tilde{D}^h_n - \lambda^h_n  
\end{pmatrix}
\begin{pmatrix} 
	U^h_n  \\ 
	W^h_n  
\end{pmatrix}
= 
\begin{pmatrix} 
	0  \\ 
	0  
\end{pmatrix}.
\label{eq:two_by_two_problem_discrete}
\end{align}
The values $\tilde{A}^h_n$, $\tilde{B}^h_n$ and $\tilde{D}^h_n$ have been analytically computed for both the standard and mixed formulation (for $p=2,3,4$)  using the approach outlined in \ref{app:splines}. They are constants that depend on $n$ and $N$, the thickness $t$, and radius $R$. For each $n$ there are two discrete eigenvalues
\begin{subequations}
\label{eq:lambdadiscreteKL}
\begin{align}
	\lambda^h_{n+} 	&= K^h_n \left(1 + L^h_n \right)			\\
	\lambda^h_{n-} 	&= K^h_n \left(1 - L^h_n \right)
\end{align}
\end{subequations}
where $K^h_n$ and $L^h_n$ are derived in \ref{sec:circulent-block-matrix} in terms of $\tilde{A}^h_n, \tilde{B}^h_n$ and $\tilde{D}^h_n$. Again, the ratio of the amplitudes of the eigenfunctions, $U^h_n / W^h_n$, satisfies
\begin{subequations}
\begin{align}
	\frac{U^h_{n+}}{W^h_{n+}} = \frac{\tilde{B}^h_n}{\lambda^h_{n+} - \tilde{A}^h_n}			\\
	\frac{U^h_{n-}}{W^h_{n-}} = \frac{\tilde{B}^h_n}{\lambda^h_{n-} - \tilde{A}^h_n}
\end{align}
\end{subequations}

\begin{figure}[!h]
	\centering
	\subfloat[$\lambda^{h}_{-}$ - standard formulation]{\includegraphics[width=0.47\textwidth]{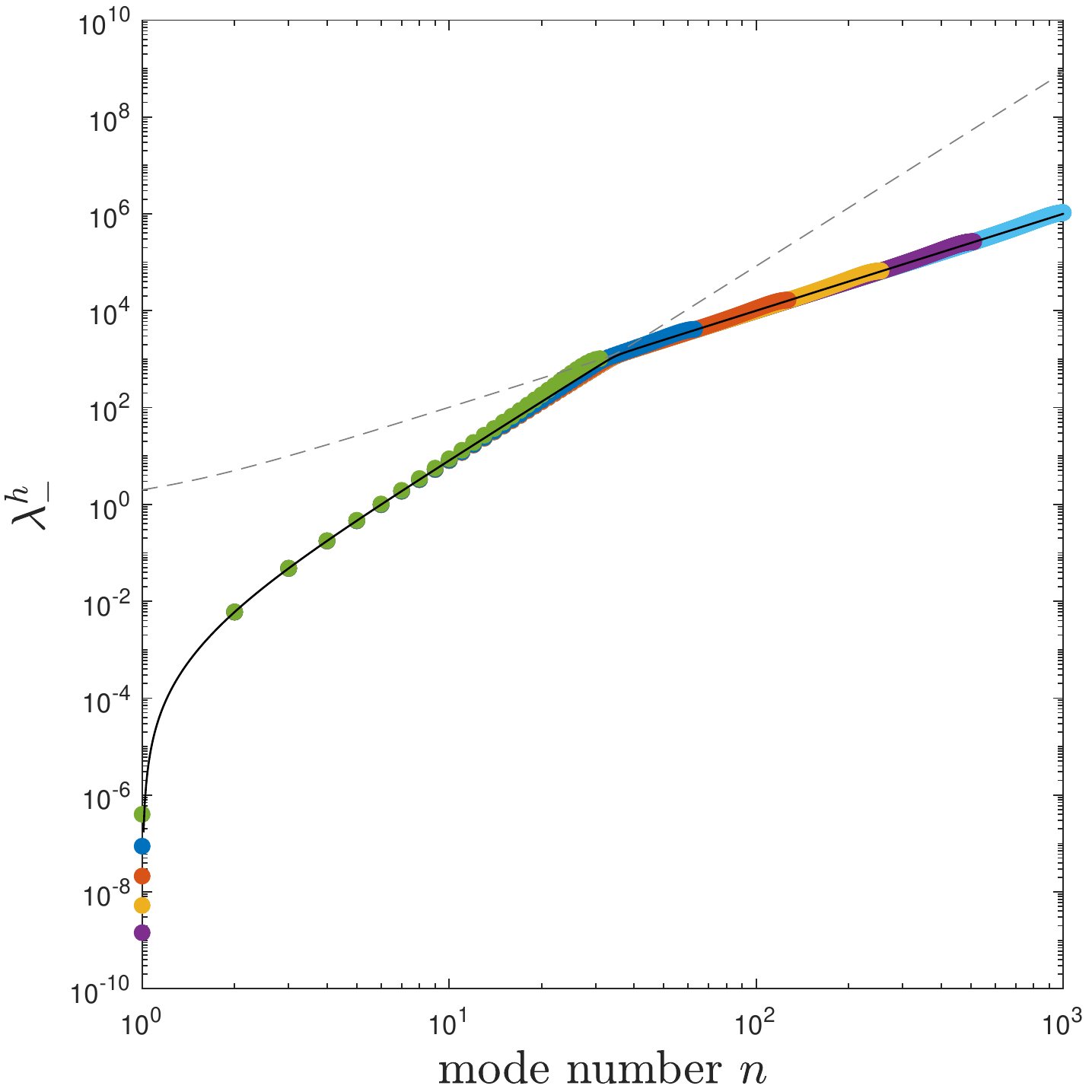}}  \hspace{0.5cm}
	\subfloat[$\lambda^{h}_{-}$ - mixed formulation]{\includegraphics[width=0.47\textwidth]{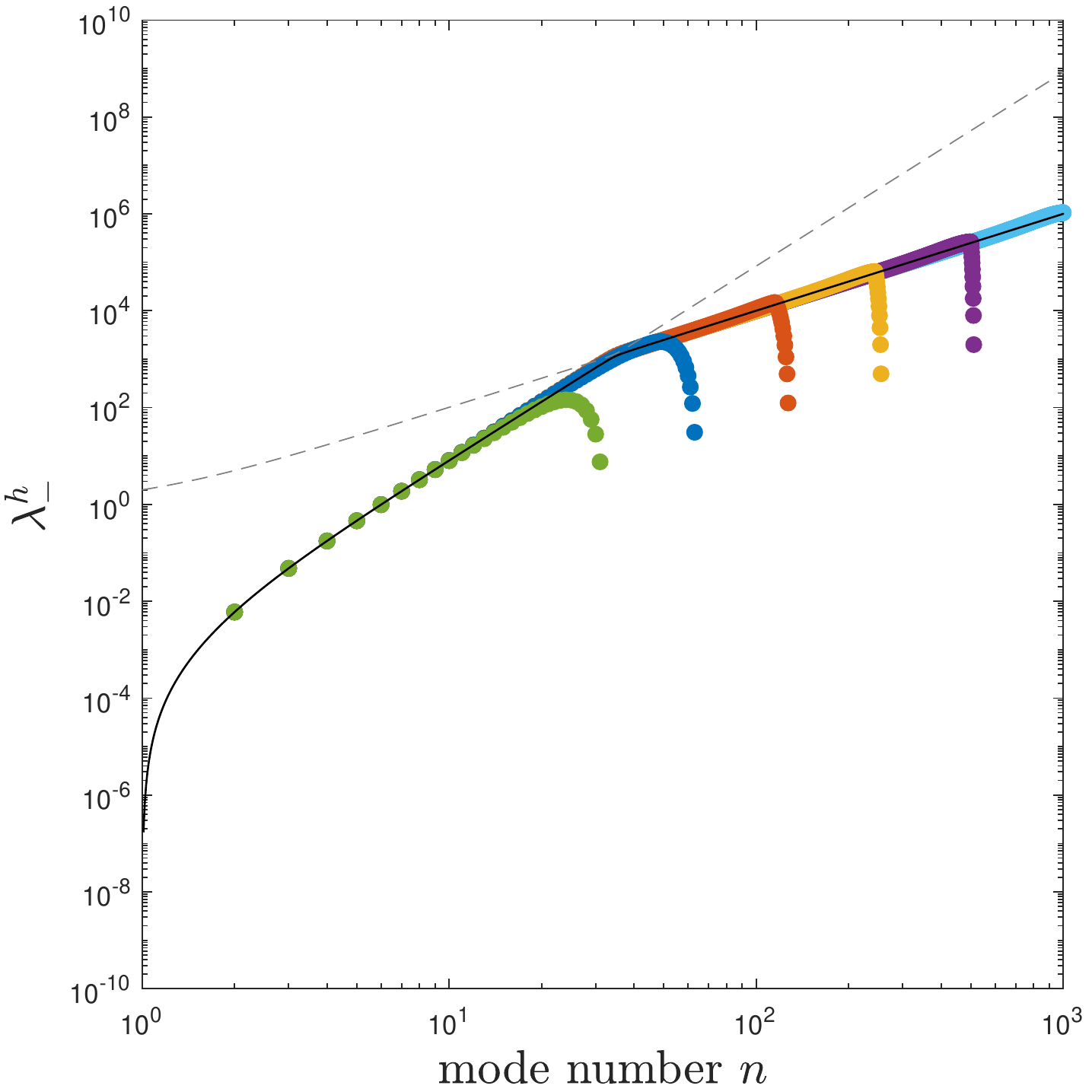}} \\
	\subfloat[$\rho^h_{n-}$ - standard formulation]{\includegraphics[width=0.47\textwidth]{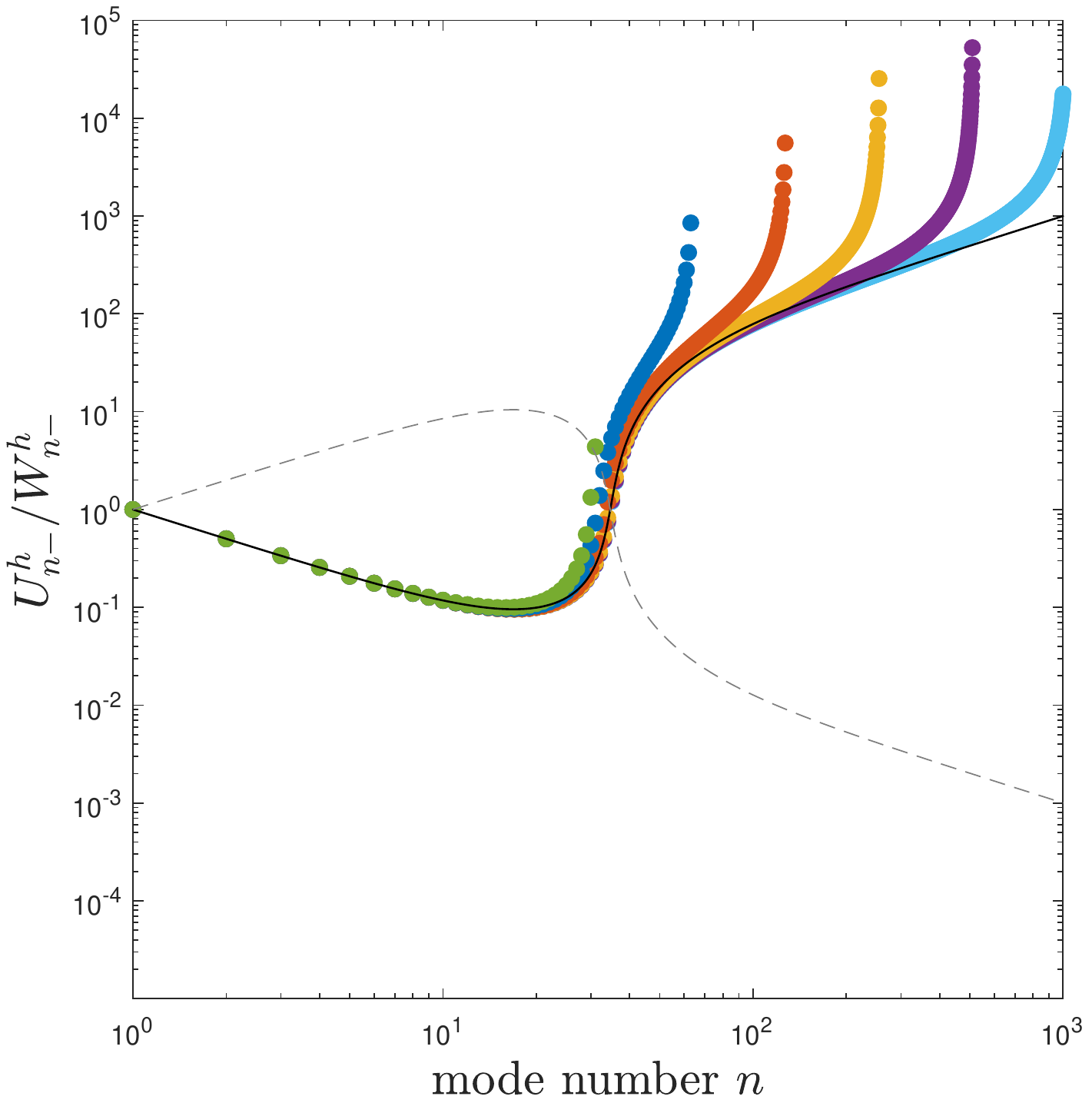}}  \hspace{0.5cm}
	\subfloat[$\rho^h_{n-}$ - mixed formulation]{\includegraphics[width=0.47\textwidth]{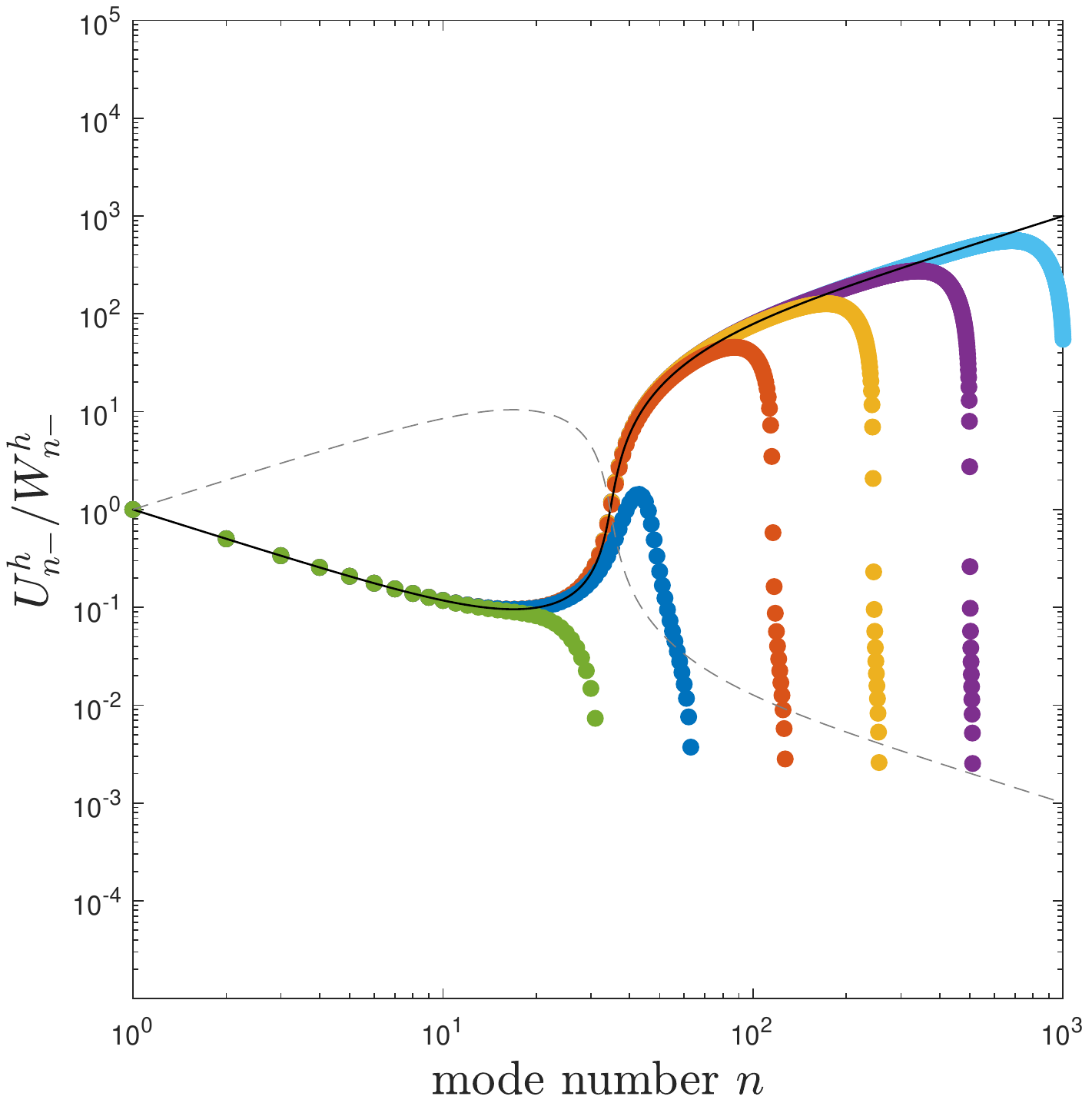}} \\
		\begin{tikzpicture}
		\filldraw[green1,line width=1pt] (2,0) circle (2pt);
		\filldraw[green1,line width=1pt] (2,0) node[right]{\scriptsize $N = 16$};	
		\filldraw[blue1,line width=1pt] (4,0) circle (2pt);
		\filldraw[blue1,line width=1pt] (4,0) node[right]{\scriptsize$N = 32$};	
		\filldraw[red1,line width=1pt] (6,0) circle (2pt);
		\filldraw[red1,line width=1pt] (6,0) node[right]{\scriptsize$N = 64$};	
		\filldraw[yellow1,line width=1pt] (8,0) circle (2pt);
		\filldraw[yellow1,line width=1pt] (8,0) node[right]{\scriptsize$N = 128$};
		\filldraw[purple1,line width=1pt] (10,0) circle (2pt);
		\filldraw[purple1,line width=1pt] (10,0) node[right]{\scriptsize$N = 256$};		
		\filldraw[lightblue1,line width=1pt] (12,0) circle (2pt);
		\filldraw[lightblue1,line width=1pt] (12,0) node[right]{\scriptsize$N = 512$};		
		\filldraw[black, line width=1pt] (13.9,0) -- (14.1,0) node[right]{\scriptsize exact};
	\end{tikzpicture}	
	\caption{Approximation of eigenvalues $\lambda_{-}$ and relative amplitudes $\rho_{n-} = U_{n-} / W_{n-}$ using $p=2$ at a normalized thickness of $\bar{t} = 0.1$.}
	\label{fig:discrete_eigenvalues_and_amplitudes}
\end{figure}
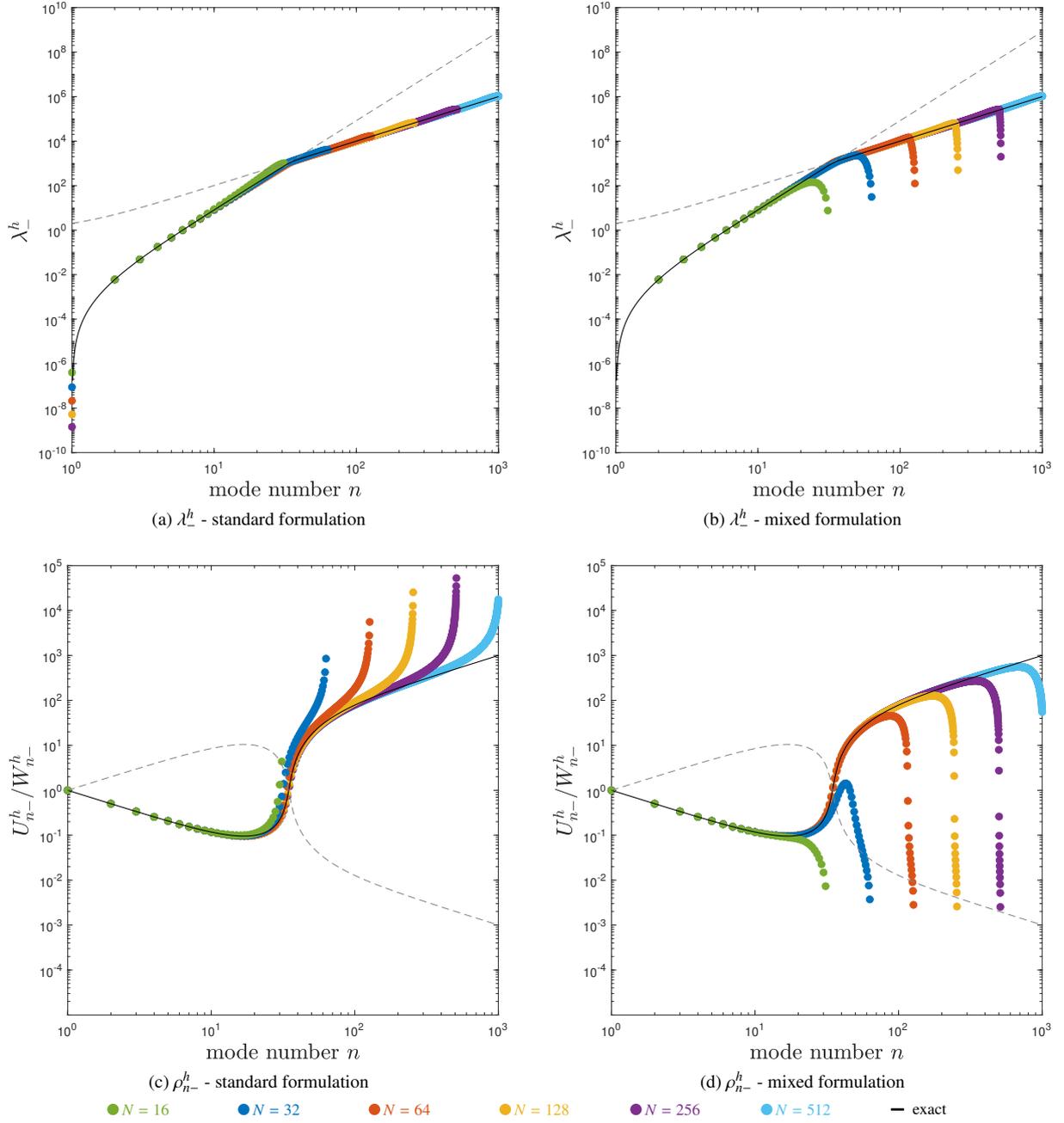

The discrete eigenvalues and the discrete ratio of the amplitudes may be directly compared with the analytical solutions developed in the previous section. Figure \ref{fig:discrete_eigenvalues_and_amplitudes} illustrates the global behavior of the approximation of the eigenvalues $\lambda_{n-}$ and relative amplitudes $U_{n-} /W_{n-}$ for a range of meshes at a normalized thickness of $\bar{t}=0.1$. We can make the following general observations
\begin{itemize}
	\item The eigenvalues, particularly $\lambda_{-}$, cover many orders of magnitude. It should be challenging to capture this with a numerical method, particularly, maintaining accuracy in the very small values. As we show later, the standard formulation struggles to accurately capture these small eigenvalues, which is a manifestation of membrane locking. Coarse discretizations suffer more than fine ones.
	\item At a discretization with $N > \hat{n}$ degrees of freedom, it is challenging to capture the two different types of physical behavior exhibited by the exact solution. We show later that the asymptotic accuracy in the amplitude ratio reduces, in relative terms, compared with results on coarse meshes. Hence, this phenomenon affects fine discretizations more than coarser ones.
\end{itemize}

In Section \ref{sec:5} we take a deeper look into the normalized eigenvalue and relative amplitude errors. In Section \ref{sec:6} we perform a mathematical analysis that sheds light on some of these observations.

\section{A criterion for assessing membrane locking \label{sec:4}}
In this section we present a criterion to assess the presence of membrane locking in thin beam and shell formulations. It is motivated by the lack of invariance of normalized eigenvalue spectra as a function of the normalized mode number. First, we show that the isolated membrane response or isolated bending response, modeled via a second order or fourth order eigenvalue problem, respectively, has a discrete normalized eigenvalue error that is invariant, that is, it is independent of the number of degrees of freedom of the discrete space. Then, we show an example where the combined membrane and bending response of a curved Euler-Bernoulli beam lead to a discrete normalized eigenvalue spectrum that lacks this invariance. This motivates the choice of the presented criterion to assess membrance locking rigorously in the end of this section.

\subsection{Invariance of the isolated membrane response}
Let $\domain$ denote the unit circle and consider the space of periodic functions $V = H^1(\domain)$. We seek $u \in V$ and $\lambda \in \mathbb{R}$ such that
\begin{align*}
	 (u_{,\theta},\, v_{,\theta}) = \lambda (u,v) \quad \forall v\in V
\end{align*}
Assuming analytical modes of the form $u_n(\theta) = A_n\cos(n \theta)$, $\theta \in [0,2\pi)$, we may determine the analytical eigenvalues as
\begin{align*}
	\lambda_n = n^2, \quad n = 0, 1, \ldots
\end{align*}

Let $V^h \subset V$ denote a space of periodic, quadratic splines of dimension $2N$, defined on a uniform partition of the unit circle with mesh size $h=\pi / N$. It can be shown that the matrix eigenvalue problem is represented by the $2N$ equations
\begin{align*}
	\frac{1}{6h} \left(-u_{A-2} - 2u_{A-1} + 6 u_{A} - 2 u_{A+1} - u_{A+2} \right) -
	\frac{\lambda^h h}{120} \left(u_{A-2} + 26 u_{A-1} + 66 u_{A} + 26 u_{A+1} + u_{A+2} \right)  = 0
\end{align*}
for $A = 1, 2, \ldots , 2N$. The discrete eigenvalues can then be computed using elementary properties of circulant matrices. The result is  
\begin{align*}
	\lambda^h_{n} = \frac{20}{h^2} \cdot \frac{6 - 4 \cos{(n h)} - 2\cos{(2 n h)}}{66 + 52 \cos{(n h)} + 2\cos{(2 n h)} }, \qquad n = 0, 1, \ldots ,N-1.
\end{align*}
The error in the non-zero eigenvalues may then be computed analytically as
\begin{align*}
	\frac{\lambda^h_{n} - \lambda_n}{\lambda_n}  = \frac{20}{(nh)^2} \cdot \frac{6 - 4 \cos{(n h)} - 2\cos{(2 n h)}}{66 + 52 \cos{(n h)} + 2\cos{(2 n h)} } - 1.
\end{align*}
In classical eigenvalue error analysis \cite[Chapter 6]{strang_analysis_2008}, one proceeds then by fixing the mode number $n$ and deriving the order of approximation in terms of the mesh size. Here it can be shown that the discrete eigenvalues converge at a rate of $h^4$ for fixed $n$. We are, however, more interested in the behavior of the eigenvalue error at fixed normalized mode number, $\xi = n / N \in [0,1]$. Using $n h = \pi \xi$, the eigenvalue error may be expressed in terms of the normalized mode number as
\begin{align}\label{eq:EVsecondorder}
	\frac{\lambda^h(\xi) - \lambda(\xi)}{\lambda(\xi)} = \frac{20}{(\pi \xi)^2} \cdot \frac{6 - 4 \cos{(\pi  \xi)} - 2\cos{(2 \pi  \xi)}}{66 + 52 \cos{(\pi   \xi)} + 2\cos{(2\pi \xi)} } - 1.
\end{align}
\emph{Note there is no dependence on mesh size!} Hence, the normalized discrete frequency spectrum is invariant, that is independent of the number of degrees of freedom $2N$. Figure \ref{fig:invariance}a illustrates this observation.

\subsection{Invariance of the isolated bending response}
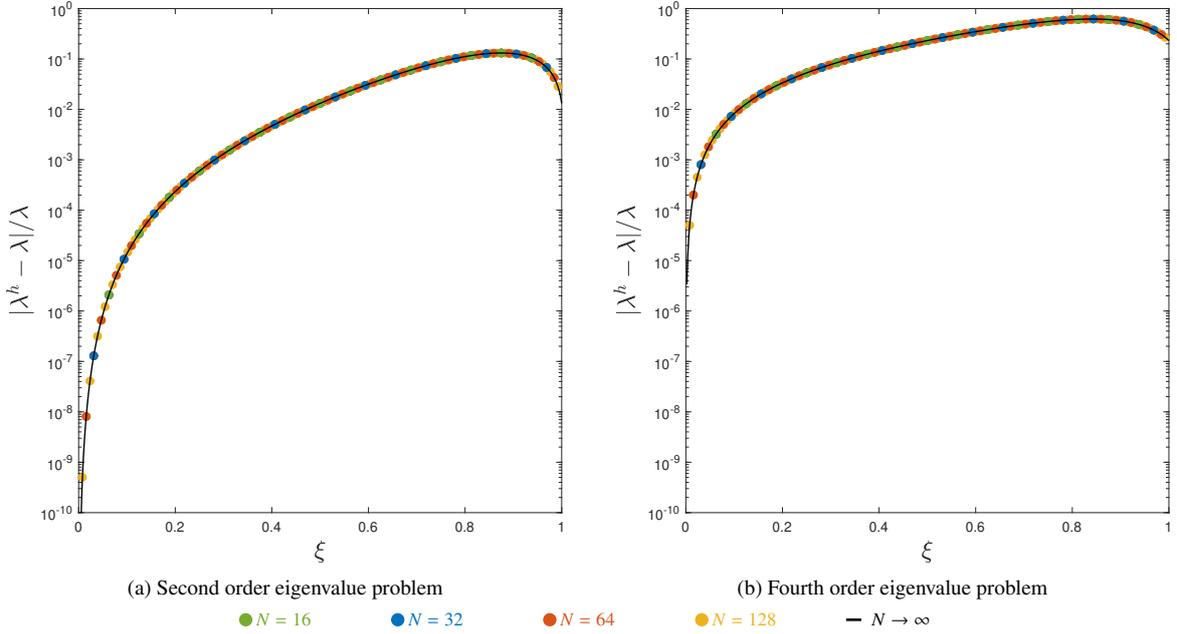
\begin{figure}
	\centering
	\subfloat[Second order eigenvalue problem]{\includegraphics[trim = 0cm 0cm 0cm 0cm, clip,width=0.45\textwidth]{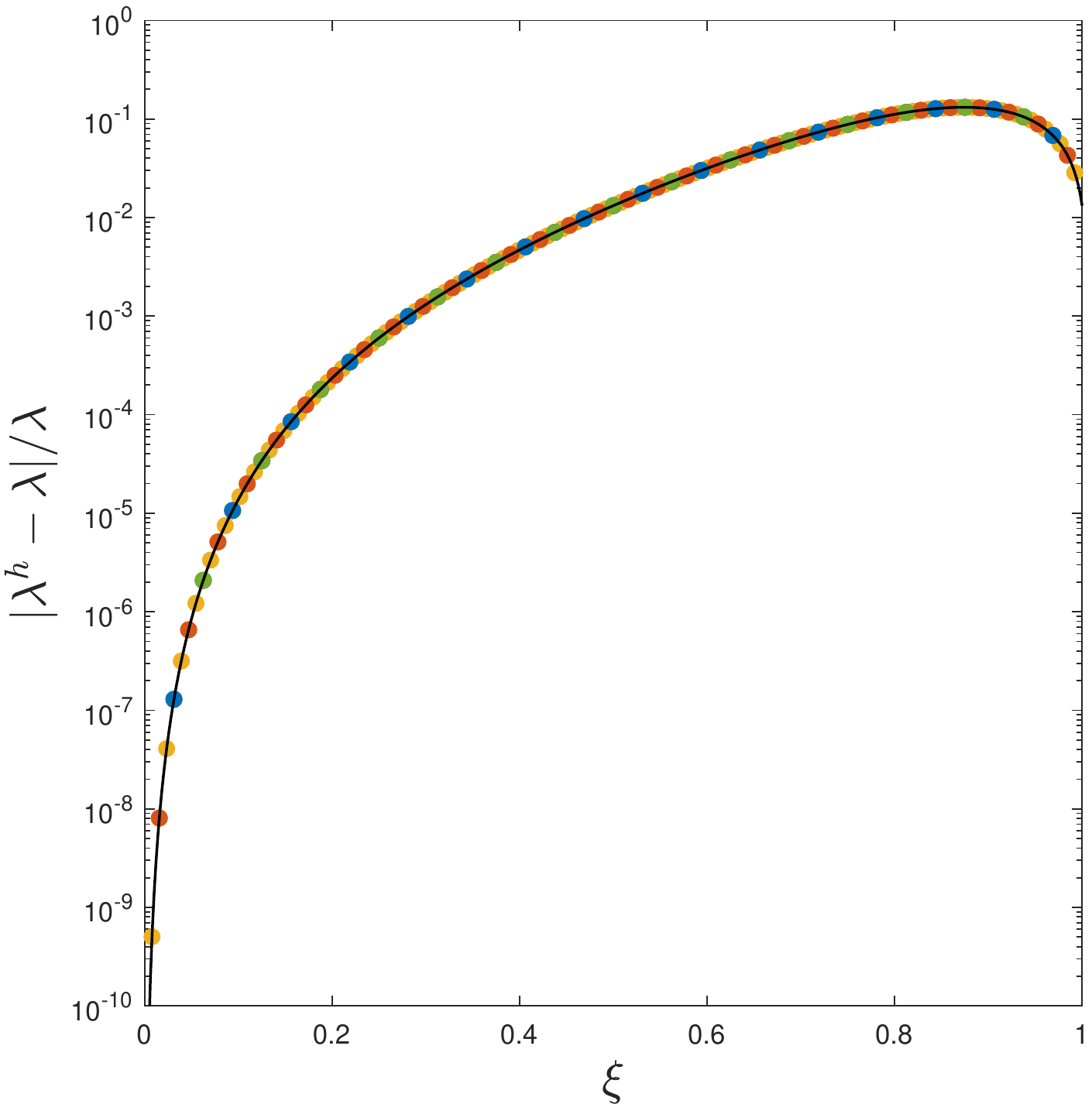}} \hspace{0.5cm}
	\subfloat[Fourth order eigenvalue problem]{\includegraphics[trim = 0cm 0cm 0cm 0cm, clip,width=0.45\textwidth]{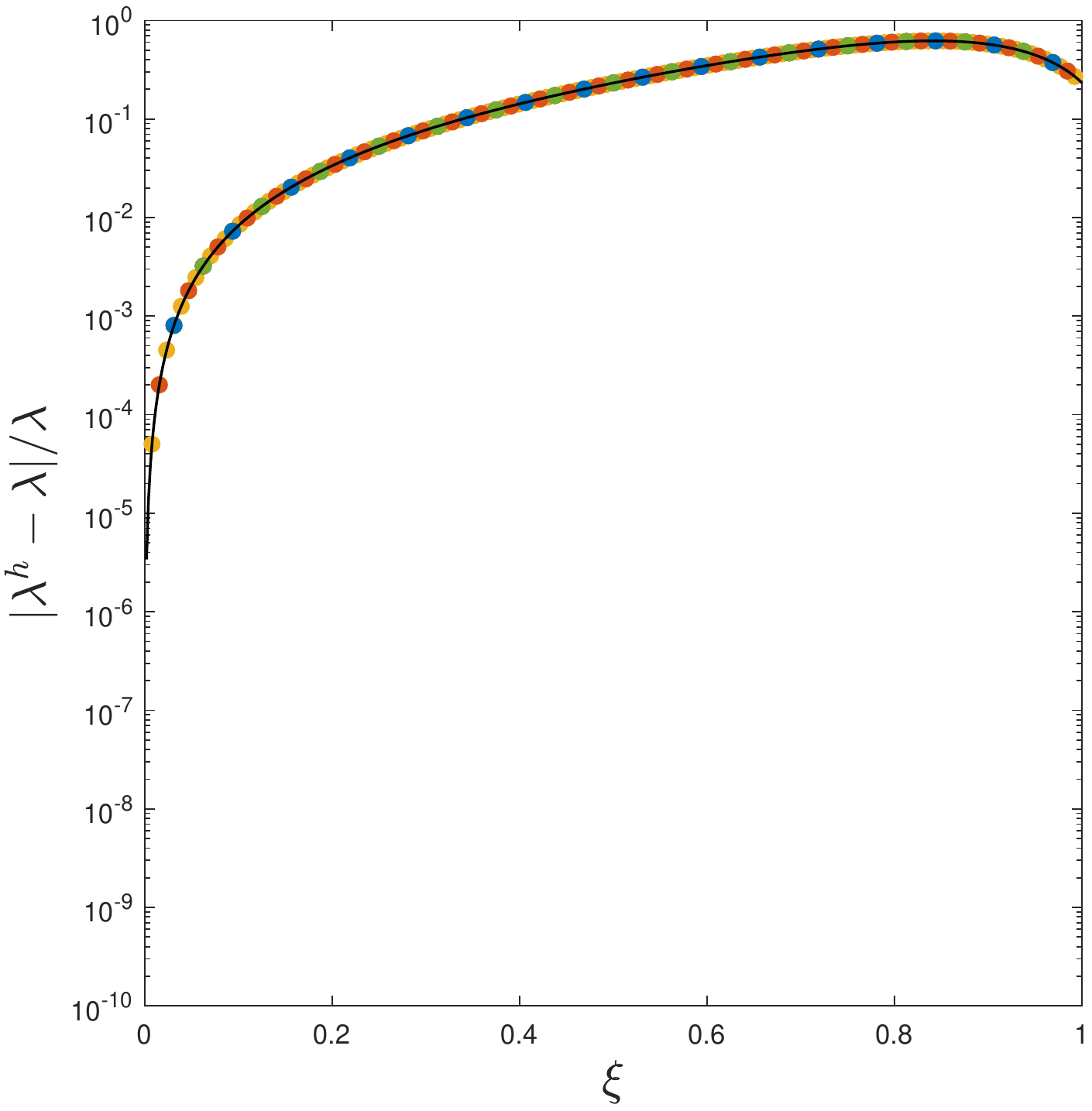}} \\
	\begin{tikzpicture}
		\filldraw[green1,line width=1pt] (2,0) circle (2pt);
		\filldraw[green1,line width=1pt] (2,0) node[right]{\scriptsize $N = 16$};	
		\filldraw[blue1,line width=1pt] (4,0) circle (2pt);
		\filldraw[blue1,line width=1pt] (4,0) node[right]{\scriptsize$N = 32$};	
		\filldraw[red1,line width=1pt] (6,0) circle (2pt);
		\filldraw[red1,line width=1pt] (6,0) node[right]{\scriptsize$N = 64$};	
		\filldraw[yellow1,line width=1pt] (8,0) circle (2pt);
		\filldraw[yellow1,line width=1pt] (8,0) node[right]{\scriptsize$N = 128$};
		\filldraw[black, line width=1pt] (9.9,0) -- (10.1,0) node[right]{\scriptsize $N \rightarrow \infty$};
	\end{tikzpicture}
	\caption{Invariance of spectra as a function of the normalized mode number.}
	\label{fig:invariance}
\end{figure}

To show that this is not merely a coincidence, we apply the same analysis to a fourth order eigenvalue problem. Consider periodic functions in $V = H^2(\domain)$. We seek $u \in V$ and $\lambda \in \mathbb{R}$ such that
\begin{align*}
	 (u_{,\theta \theta},\, v_{,\theta \theta}) = \lambda (u,v) \quad \forall v\in V
\end{align*}
Assuming analytical modes of the form $u_n(\theta) = A_n\cos(n \theta)$, $\theta \in [0,2\pi)$, we may determine the analytical eigenvalues as
\begin{align*}
	\lambda_n = n^4, \quad n = 0, 1, \ldots
\end{align*}
Let $V^h \subset V$ denote the same space of periodic, quadratic splines of dimension $2N$, defined on a uniform partition of the unit circle with mesh size $h=\pi / N$. It can be shown that the matrix eigenvalue problem is represented by the $2N$ equations
\begin{align*}
	\frac{1}{h^3} \left(u_{A-2} - 4u_{A-1} + 6u_{A} - 4u_{A+1} + u_{A+2} \right) -
	\frac{\lambda^h h}{120} \left(u_{A-2} + 26 u_{A-1} + 66 u_{A} + 26 u_{A+1} + u_{A+2} \right)  = 0
\end{align*}
for $A = 1, 2, \ldots , 2N$. The discrete eigenvalues are computed as
\begin{align*}
	\lambda^h_{n} = \frac{120}{h^4} \cdot \frac{6 - 8 \cos{(n h)} + 2\cos{(2 n h)}}{66 + 52 \cos{(n h)} + 2\cos{(2 n h)} }, \qquad n = 0, 1, \ldots ,N-1.
\end{align*}
Analoguous to the second order eigenvalue problem, using $n h = \pi \xi$, the eigenvalue error may be expressed as
\begin{align}\label{eq:EVfourthorder}
	\frac{\lambda^h(\xi) - \lambda(\xi)}{\lambda(\xi)} = \frac{120}{(\pi \xi)^4} \cdot \frac{6 - 8 \cos{(\pi  \xi)} + 2\cos{(2 \pi  \xi)}}{66 + 52 \cos{(\pi   \xi)} + 2\cos{(2\pi \xi)} } - 1.
\end{align}
Importantly, once again, the expression is invariant with respect to mesh refinement! This is illustrated in Figure \ref{fig:invariance}b.

\subsection{Definition of a rigorous criterion based on spectrum invariance}\label{sec:lockingcriterion}
\begin{figure}
	\centering
	\includegraphics[trim = 0cm 0cm 0cm 0cm, clip,width=0.75\textwidth]{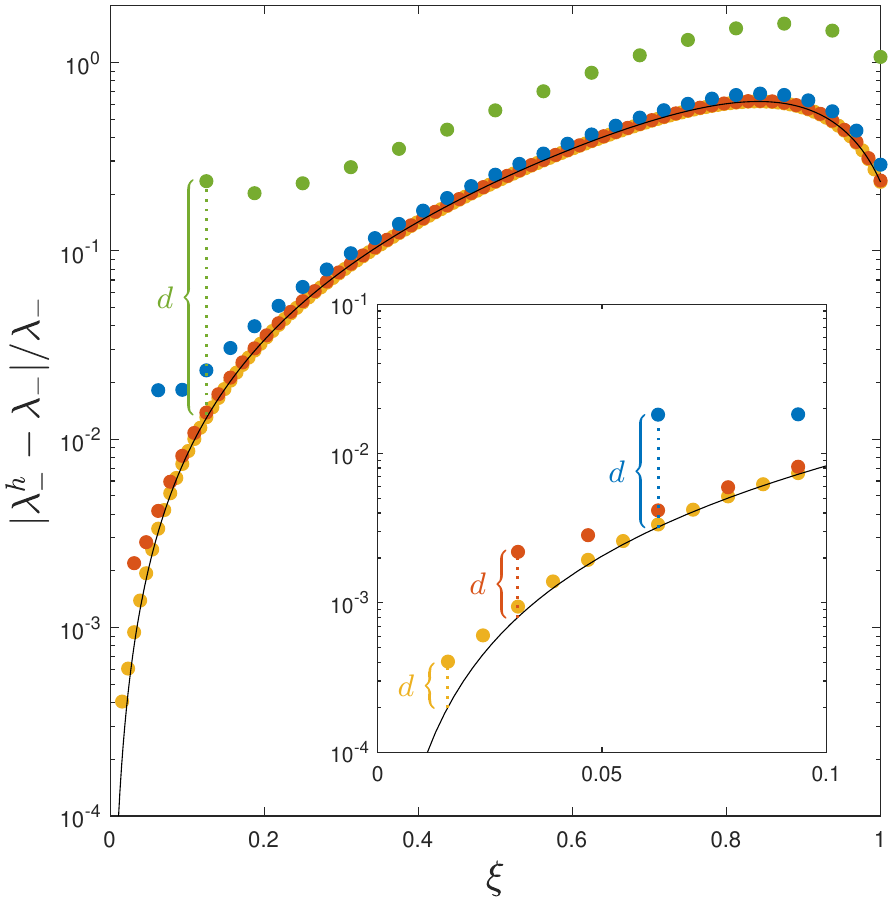} \\
	\begin{tikzpicture}
		\filldraw[green1,line width=1pt] (2,0) circle (2pt);
		\filldraw[green1,line width=1pt] (2,0) node[right]{\scriptsize $N = 16$};	
		\filldraw[blue1,line width=1pt] (4,0) circle (2pt);
		\filldraw[blue1,line width=1pt] (4,0) node[right]{\scriptsize$N = 32$};	
		\filldraw[red1,line width=1pt] (6,0) circle (2pt);
		\filldraw[red1,line width=1pt] (6,0) node[right]{\scriptsize$N = 64$};	
		\filldraw[yellow1,line width=1pt] (8,0) circle (2pt);
		\filldraw[yellow1,line width=1pt] (8,0) node[right]{\scriptsize$N = 128$};
		\filldraw[black, line width=1pt] (9.9,0) -- (10.1,0) node[right]{\scriptsize $N \rightarrow \infty$};
	\end{tikzpicture}
	\caption{Membrane locking is manifested by a lack of invariance of the discrete normalized spectrum. The discrepancy between the spectrum error and the asymptotic spectrum error, denoted by $d$, measures the extent of locking.}
	\label{fig:locking}
\end{figure}
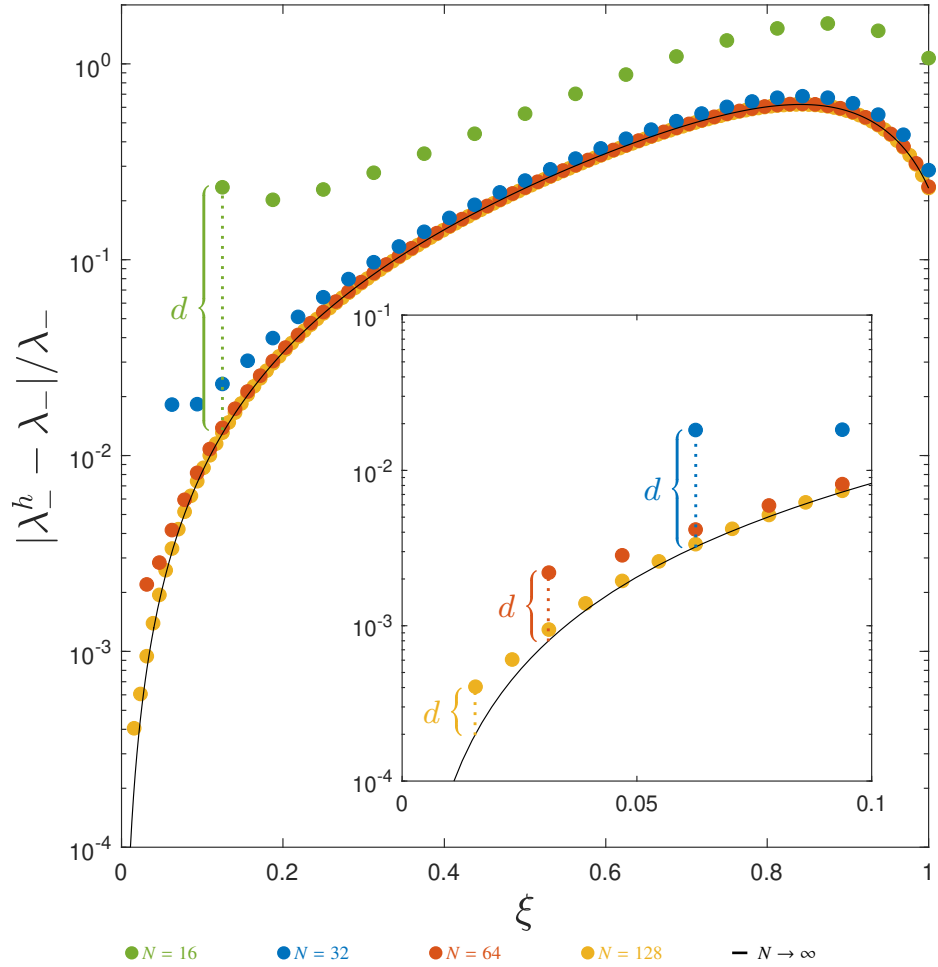

The isolated membrane and bending responses are invariant when viewed as functions of the normalized mode number $\xi$. This may no longer be the case when they are coupled in a curved Euler-Bernoulli beam. Figure \ref{fig:locking} depicts a preview of the normalized eigenvalue error obtained with quadratic splines. The eigenvalue error is clearly not invariant with respect to the number of degrees of freedom of the spline space. Particularly, the lower eigenvalues suffer from a visible discrepancy compared to the rest, which can manifest partly by a larger preasymptotic regime before converging. We identify this as membrane locking. The previous discussion motivates the following characterization of membrane locking

\begin{quote}
\emph{Membrane locking is measured, at fixed normalized mode number, as the discrepancy between the spectrum error and the spectrum error in the limit of asymptotic refinement.}
\end{quote}

To make this idea more precise, use $n = \xi N$, and let
\begin{align}\label{eq:relerrordef}
	e_N(\xi) := \frac{|\lambda^h_{\xi N} - \lambda_{\xi N}|}{\lambda_{\xi N}}.
\end{align}
In mathematical notation, the above statement may then be translated as follows.
 
\begin{figure}
	\centering
	\includegraphics[trim = 0cm 0cm 0cm 0cm, clip,width=0.75\textwidth]{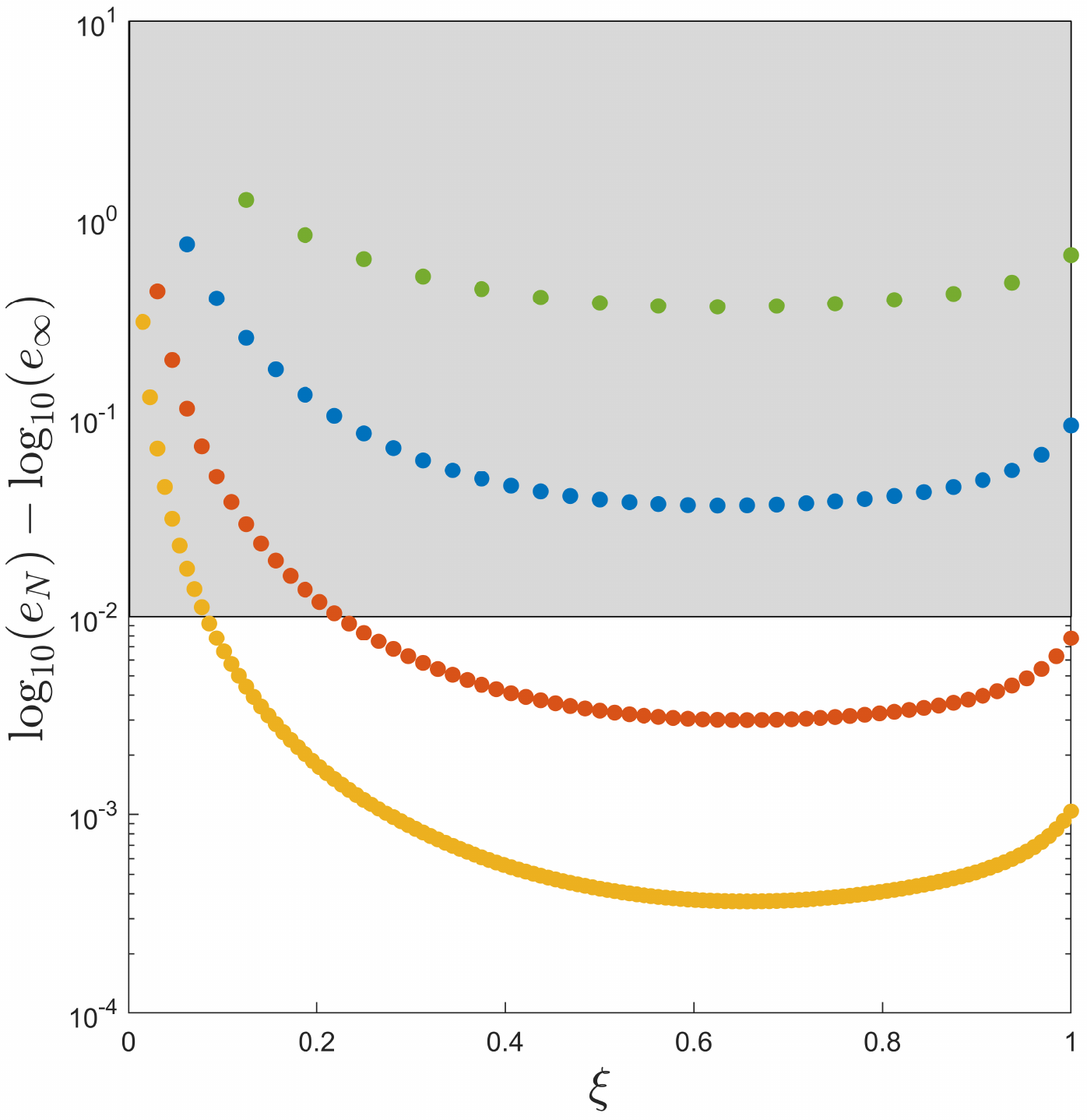} \\
	\begin{tikzpicture}
		\filldraw[green1,line width=1pt] (2,0) circle (2pt);
		\filldraw[green1,line width=1pt] (2,0) node[right]{\scriptsize $N = 16$};	
		\filldraw[blue1,line width=1pt] (4,0) circle (2pt);
		\filldraw[blue1,line width=1pt] (4,0) node[right]{\scriptsize$N = 32$};	
		\filldraw[red1,line width=1pt] (6,0) circle (2pt);
		\filldraw[red1,line width=1pt] (6,0) node[right]{\scriptsize$N = 64$};	
		\filldraw[yellow1,line width=1pt] (8,0) circle (2pt);
		\filldraw[yellow1,line width=1pt] (8,0) node[right]{\scriptsize$N = 128$};
		\filldraw[black, line width=1pt] (9.9,0) -- (10.1,0) node[right]{\scriptsize $N \rightarrow \infty$};
	\end{tikzpicture}
	\caption{Criterion to assess locking according to definition \ref{def:locking} with $\varepsilon = 1e-2$. Values in the gray area denote mode numbers for which there is measurable membrane locking.}
	\label{fig:locking_criterion}
\end{figure}

\begin{definition}[Spectral characterization of membrane locking \label{def:locking}] Let $0 < \varepsilon << 1$ denote a user-prescribed tolerance. A discretization is characterized as locking-free if and only if
\begin{align}
	\log_{10}{\left(e_N(\xi)\right)} - \lim_{\mathbb{Z}_N \ni M \rightarrow \infty} \log_{10}{\left(e_M(\xi)\right)} < \varepsilon
	\qquad \forall \xi = n / N, \qquad n = 1, 2, \ldots , N.
\end{align}
\end{definition}
The condition that $M \in \mathbb{Z}_N$, the set of integers modulo $N$, is needed to ensure there exists an integer $m$ such that $m/M = n/N = \xi$. The tolerance $\varepsilon$ is a distance in log-space, say $\varepsilon = 0.01$, beyond which we cannot observe (in log-space) a noticeable difference between the spectrum error and its asymptotic result. Figure \ref{fig:locking_criterion} depicts the distance between the spectrum error and the asymptotic spectrum error, $d = \log_{10}{\left(e_N(\xi)\right)} - \log_{10}{\left(e_\infty(\xi)\right)}$, corresponding to the results in Figure \ref{fig:locking} obtained using the standard formulation with $p=2$. According to our characterization of locking, locking is exhibited by all eigenmodes associated to the negative eigenvalue branch for discretizations with $2N=32$ and $2N=64$ degrees of freedom. In discretizations with $2N=128$ and $2N=256$, still about $20\%$ and $10\%$ of the modes lock, respectively.

\section{Spectral analysis results\label{sec:5}}
In this section, we assess the extent of membrane locking in the normalized eigenvalue spectra of both the standard and mixed formulation, employing the locking criterion outlined in Section \ref{sec:4}. While our primary focus lies on quadratic discretizations, we also present findings that validate membrane locking in cubic and quartic discretizations. Furthermore, we explore intriguing behaviors that manifest in fine discretizations and mode amplitudes.

\subsection{Presentation of the data}
We present the results pertaining to the two types of eigenvalues, namely $\lambda_{-}$ and $\lambda_{+}$, for both the standard and mixed formulation. In certain instances, we focus on results associated with the eigenvalues $\lambda_{-}$ of the standard formulation, as they exhibit pronounced locking effects. The presented figures share the following attributes:
\begin{itemize}
	\item Normalized eigenvalue errors and mode amplitudes are shown for meshes with $2N=32, \, 64, \, 128$, and $256$ degrees of freedom. The breathing mode ($n=0$) and the rigid body mode ($n=1$), see Figure \ref{fig:ring_first12modes2}, are exact or well approximated on any mesh and are not plotted. The remaining $N-2$ eigenvalues have multiplicity two. Only the unique ones are plotted. 
\item Eigenvalue error and mode amplitudes are normalized with respect to the normalized mode number $\xi = n / N$.
\item Four asymptotic curves are drawn in each graph that allow a succinct comparison: two for each method and two corresponding to each type of eigenvalue. The asymptotic curve under consideration is plotted as a solid curve, while the other three are dotted curves.
\item All results are obtained for a normalized thickness $\bar{t} = 0.015$. 
\end{itemize}

\subsection{Normalized error in eigenvalues}
\label{sec:eigenvalueerrors}
\begin{figure}[!h]
	\centering
	\subfloat[$\lambda^h_{-}$ - Standard formulation]{\includegraphics[trim = 0cm 0cm 0cm 0cm, clip,width=0.48\textwidth]{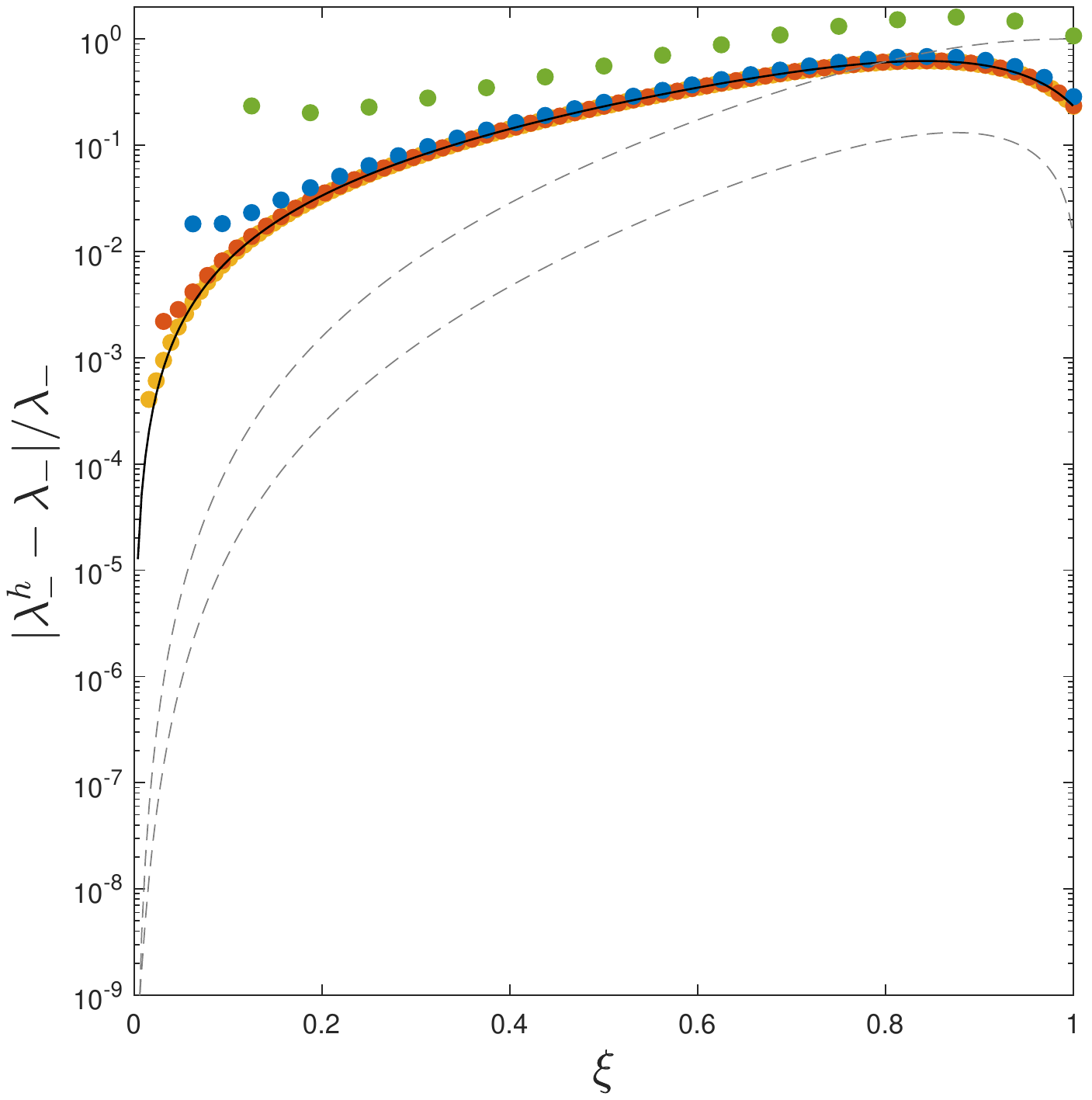}} \hspace{0.5cm}
	\subfloat[$\lambda^h_{-}$ - Mixed formulation]{\includegraphics[trim = 0cm 0cm 0cm 0cm, clip,width=0.48\textwidth]{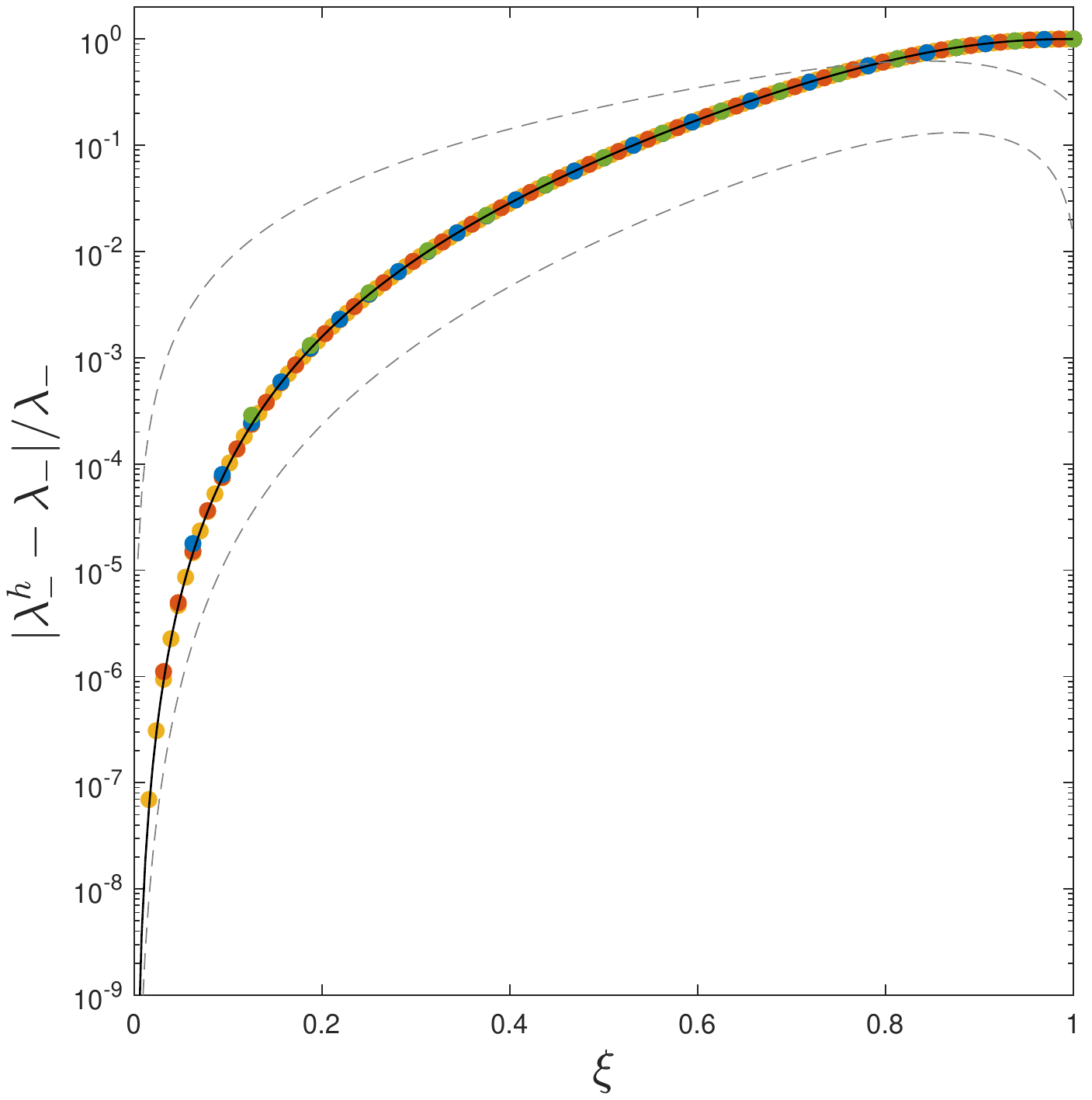}} \\
	\subfloat[$\lambda^h_{+}$ - Standard formulation]{\includegraphics[trim = 0cm 0cm 0cm 0cm, clip,width=0.48\textwidth]{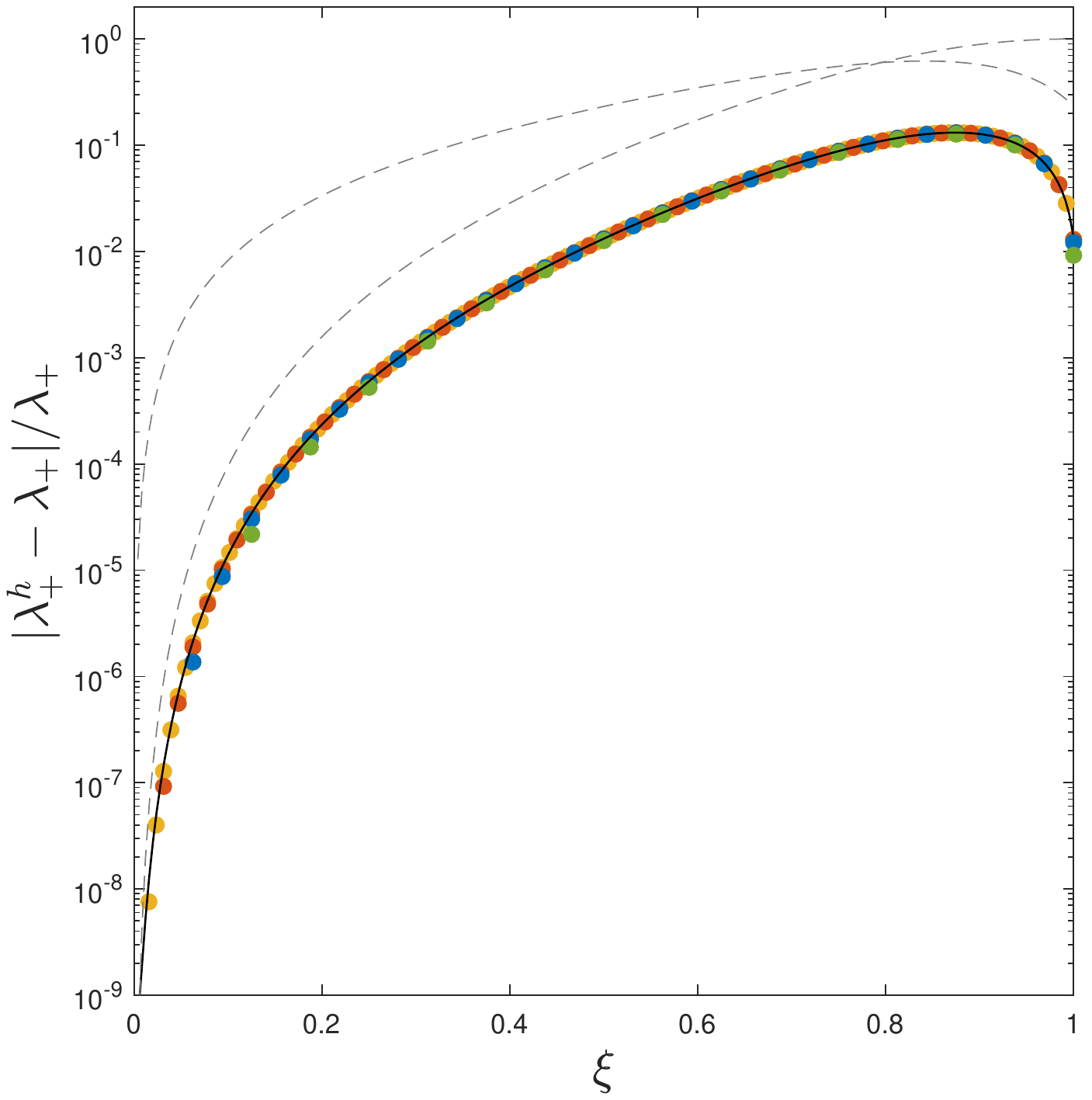}} \hspace{0.5cm}
	\subfloat[$\lambda^h_{+}$ - Mixed formulation]{\includegraphics[trim = 0cm 0cm 0cm 0cm, clip,width=0.48\textwidth]{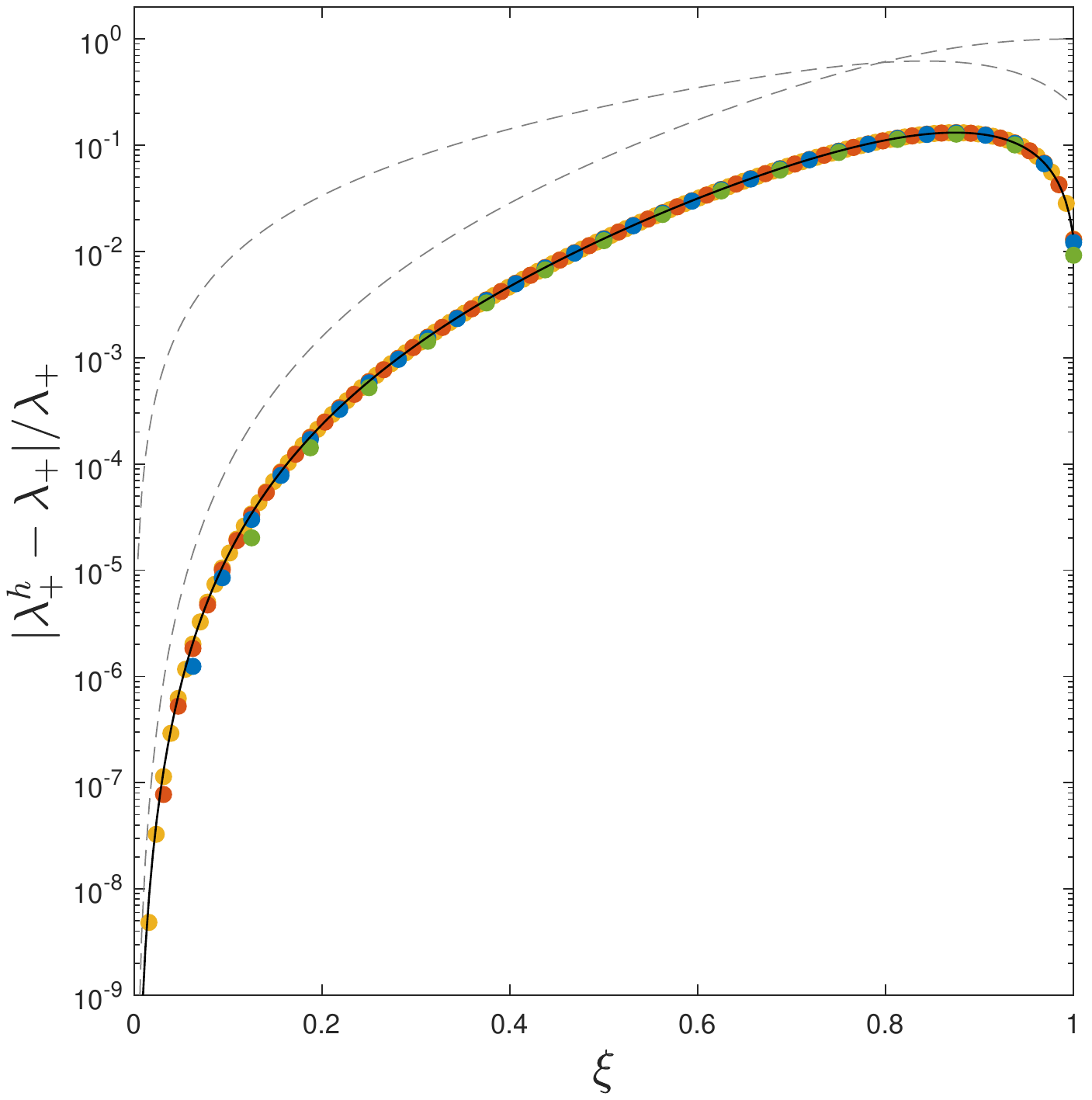}}	\\
	\begin{tikzpicture}
		\filldraw[green1,line width=1pt] (2,0) circle (2pt);
		\filldraw[green1,line width=1pt] (2,0) node[right]{\scriptsize $N = 16$};	
		\filldraw[blue1,line width=1pt] (4,0) circle (2pt);
		\filldraw[blue1,line width=1pt] (4,0) node[right]{\scriptsize$N = 32$};	
		\filldraw[red1,line width=1pt] (6,0) circle (2pt);
		\filldraw[red1,line width=1pt] (6,0) node[right]{\scriptsize$N = 64$};	
		\filldraw[yellow1,line width=1pt] (8,0) circle (2pt);
		\filldraw[yellow1,line width=1pt] (8,0) node[right]{\scriptsize$N = 128$};
		\filldraw[black, line width=1pt] (9.9,0) -- (10.1,0) node[right]{\scriptsize $N \rightarrow \infty$};
	\end{tikzpicture}
	\caption{Normalized eigenvalue error obtained with (a) standard and (b) mixed quadratic discretizations.}
	\label{fig:normalized_spectra_error}
\end{figure}
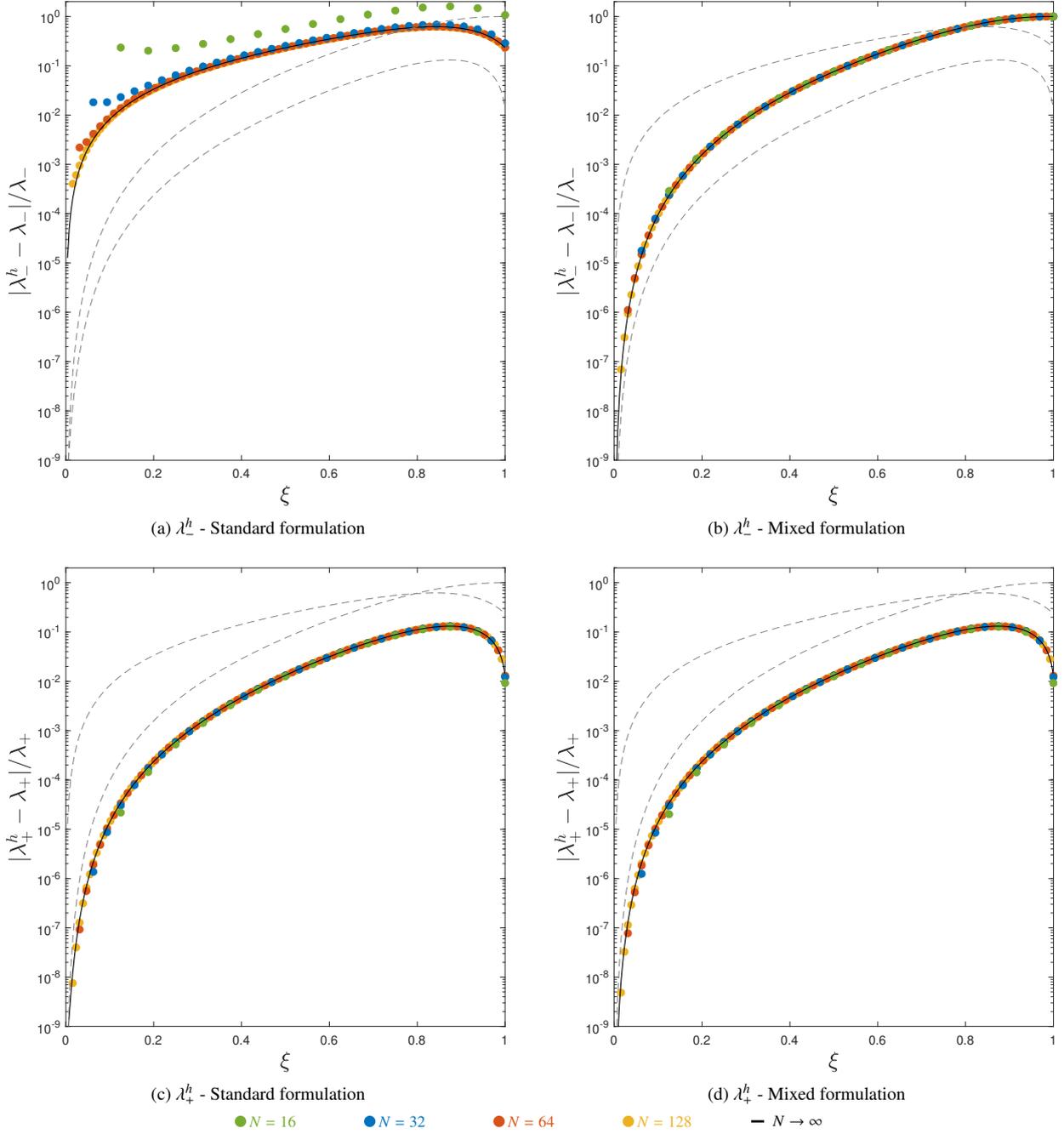

Figure \ref{fig:normalized_spectra_error} showcases the normalized eigenvalue error on a logarithmic scale. As discussed in the previous section, the normalized eigenvalue error, $\frac{\lambda^h_{-} - \lambda_-}{\lambda_{-}}$, of the standard formulation reveals significant membrane locking, as the errors deviate considerably from the asymptotic limit curves. Conversely, in the mixed method, these eigenvalues are accurately approximated, exhibiting minimal signs of locking. Additionally, it is worth noting that the asymptotic accuracy of the mixed method surpasses that of the standard formulation, evident from a comparison of the limit curves between the two methods. For the eigenvalues $\lambda_{+}$, both methods yield equally satisfactory approximations, with no indications of locking as the results align with or slightly fall below the limit curve.

\begin{figure}[!h]
	\centering
	\subfloat[Standard formulation]{\includegraphics[trim = 0cm 0cm 0cm 0cm, clip,width=0.48\textwidth]{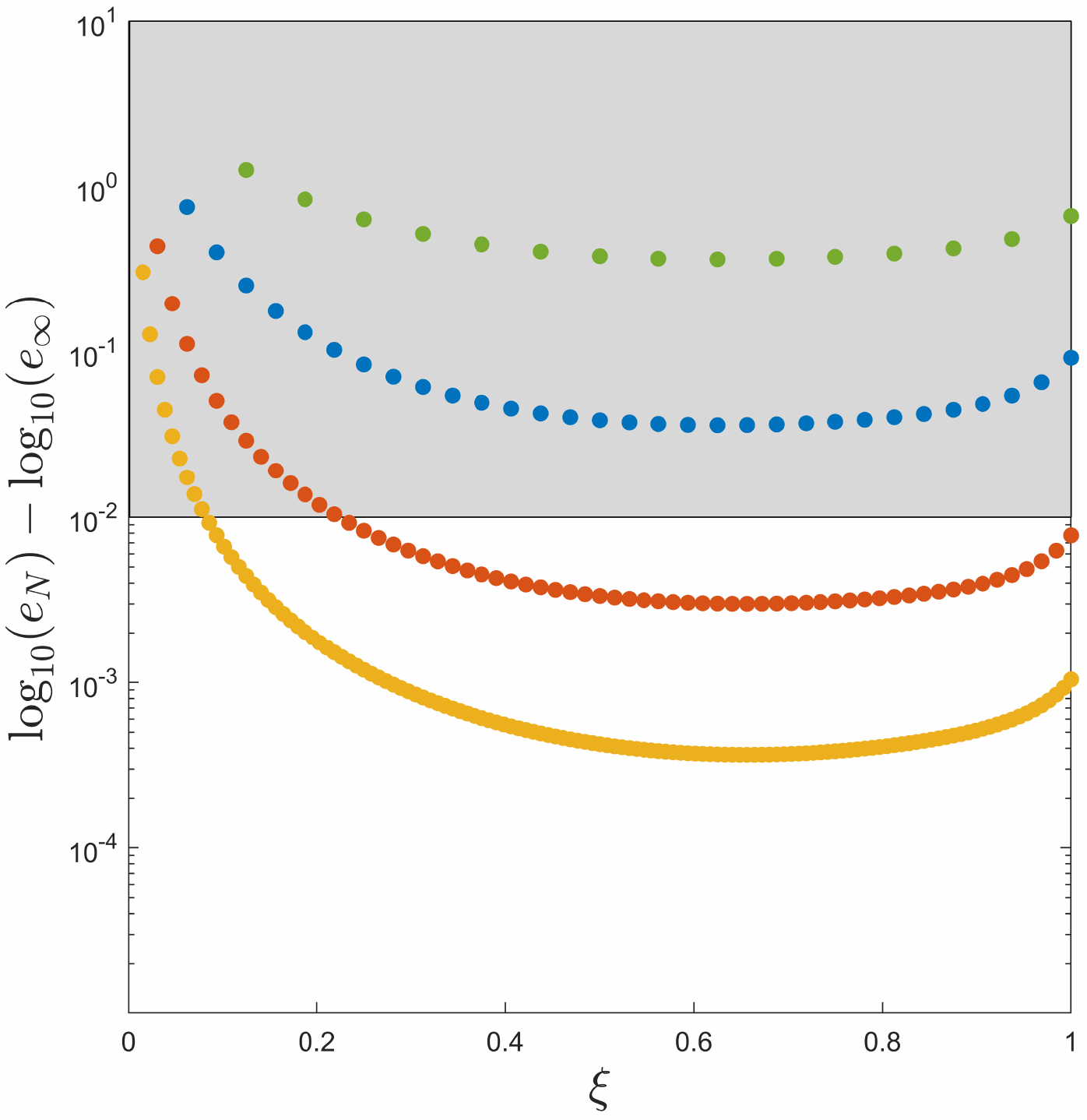}} \hspace{0.5cm}
	\subfloat[Mixed formulation]{\includegraphics[trim = 0cm 0cm 0cm 0cm, clip,width=0.48\textwidth]{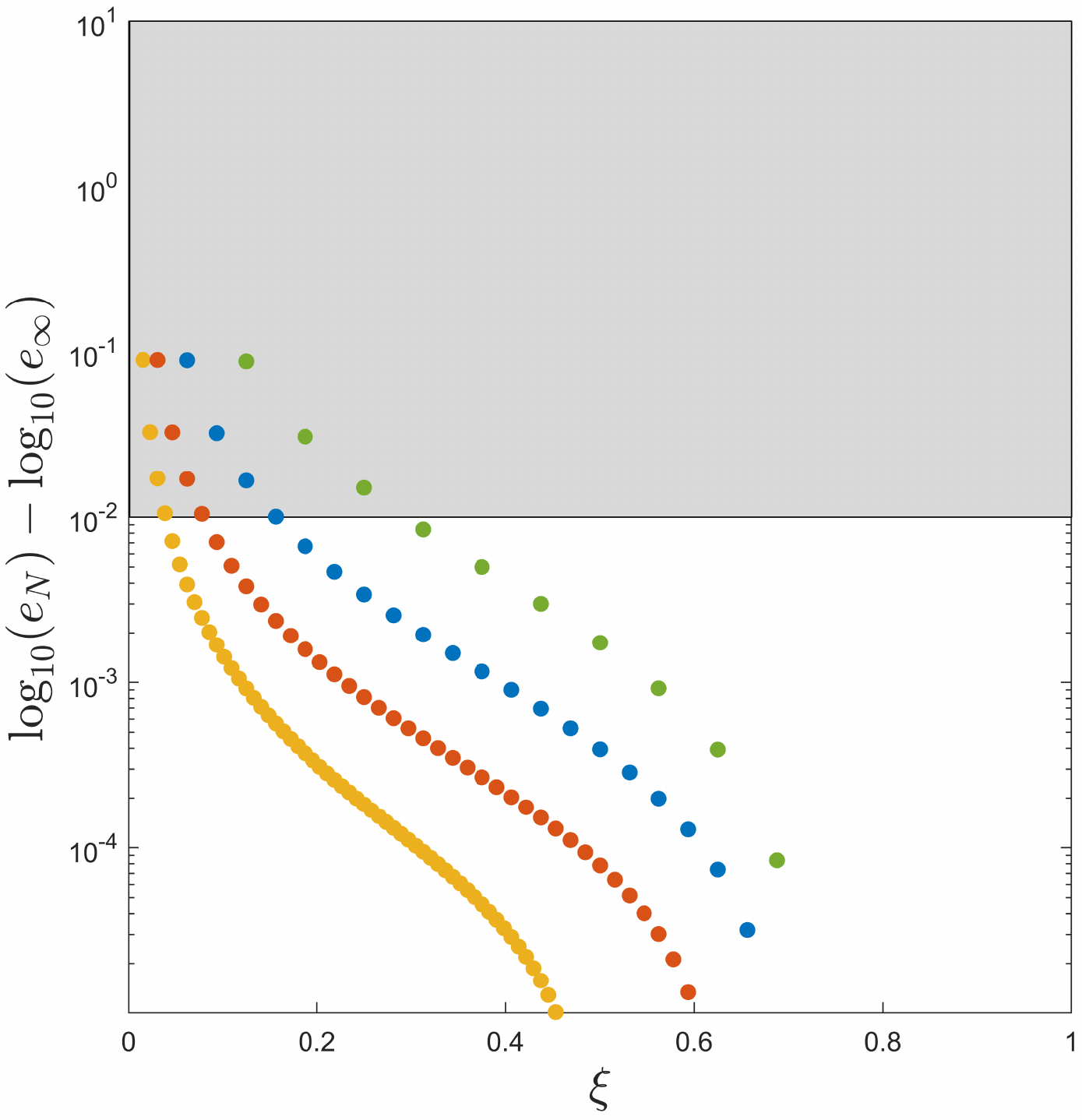}} \\
	\begin{tikzpicture}
		\filldraw[green1,line width=1pt] (2,0) circle (2pt);
		\filldraw[green1,line width=1pt] (2,0) node[right]{\scriptsize $N = 16$};	
		\filldraw[blue1,line width=1pt] (4,0) circle (2pt);
		\filldraw[blue1,line width=1pt] (4,0) node[right]{\scriptsize$N = 32$};	
		\filldraw[red1,line width=1pt] (6,0) circle (2pt);
		\filldraw[red1,line width=1pt] (6,0) node[right]{\scriptsize$N = 64$};	
		\filldraw[yellow1,line width=1pt] (8,0) circle (2pt);
		\filldraw[yellow1,line width=1pt] (8,0) node[right]{\scriptsize$N = 128$};
		\filldraw[black, line width=1pt] (9.9,0) -- (10.1,0) node[right]{\scriptsize $N \rightarrow \infty$};
	\end{tikzpicture}
	\caption{Identification of eigenvalues $\lambda^h_{-}$ that display membrane locking ($\varepsilon = 1 \cdot 10^{-2}$) in quadratic discretizations.}
	\label{fig:locking_p2}
\end{figure}

By calculating the logarithmic distance between the eigenvalue errors and the asymptotic curves, we can rigorously assess the extent of membrane locking using the criterion outlined in the previous section. Figure \ref{fig:locking_p2} illustrates the results. As previously discussed, the standard formulation exhibits significant locking for coarse to medium refined discretizations. In contrast, the mixed method demonstrates only minor locking in the first three eigenvalues, regardless of the level of refinement.

\begin{figure}[!h]
	\centering
	\subfloat[Standard formulation $p=3$]{\includegraphics[trim = 0cm 0cm 0cm 0cm, clip,width=0.45\textwidth]{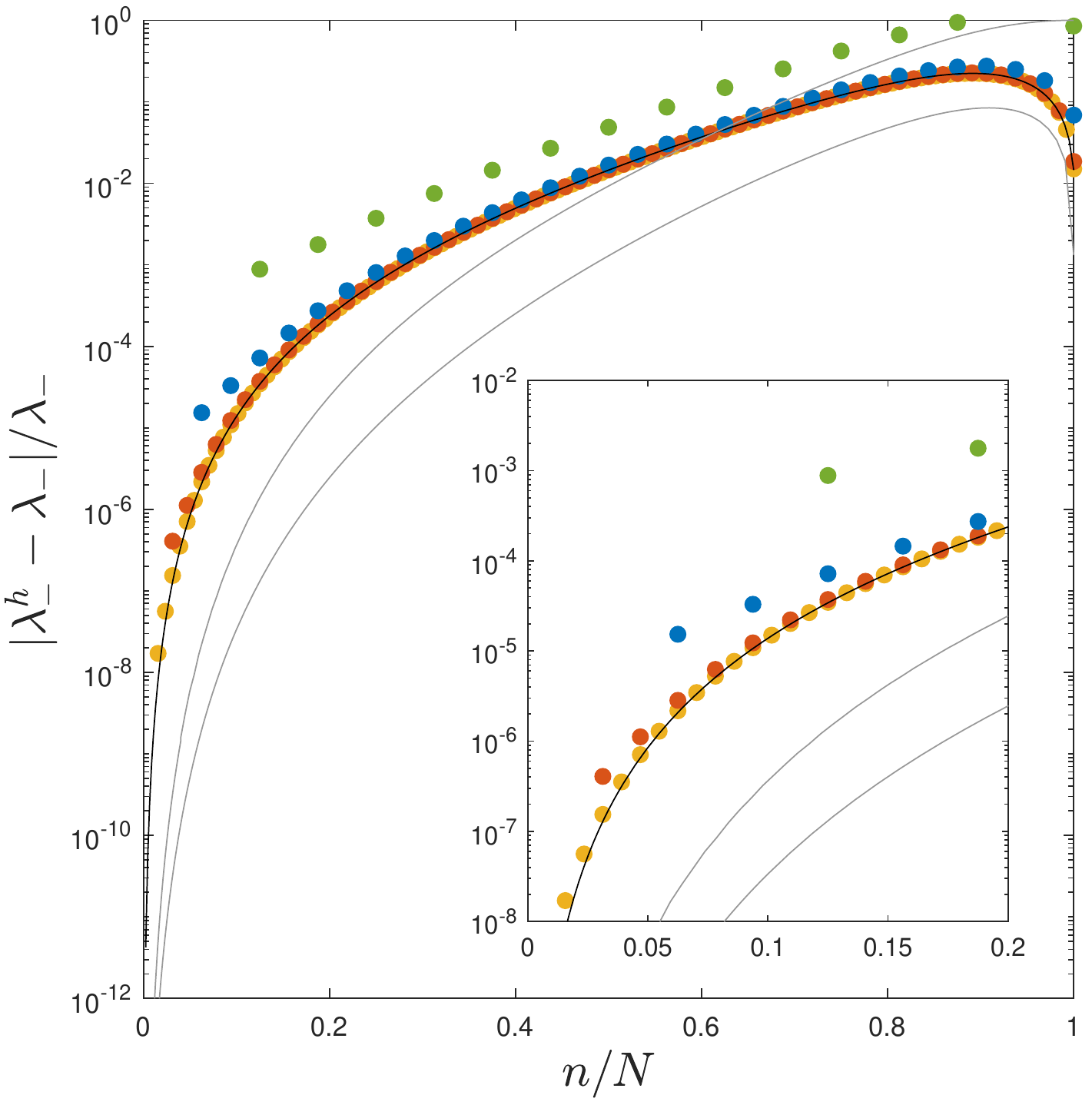}} \hspace{0.5cm}
	\subfloat[Standard formulation $p=4$]{\includegraphics[trim = 0cm 0cm 0cm 0cm, clip,width=0.45\textwidth]{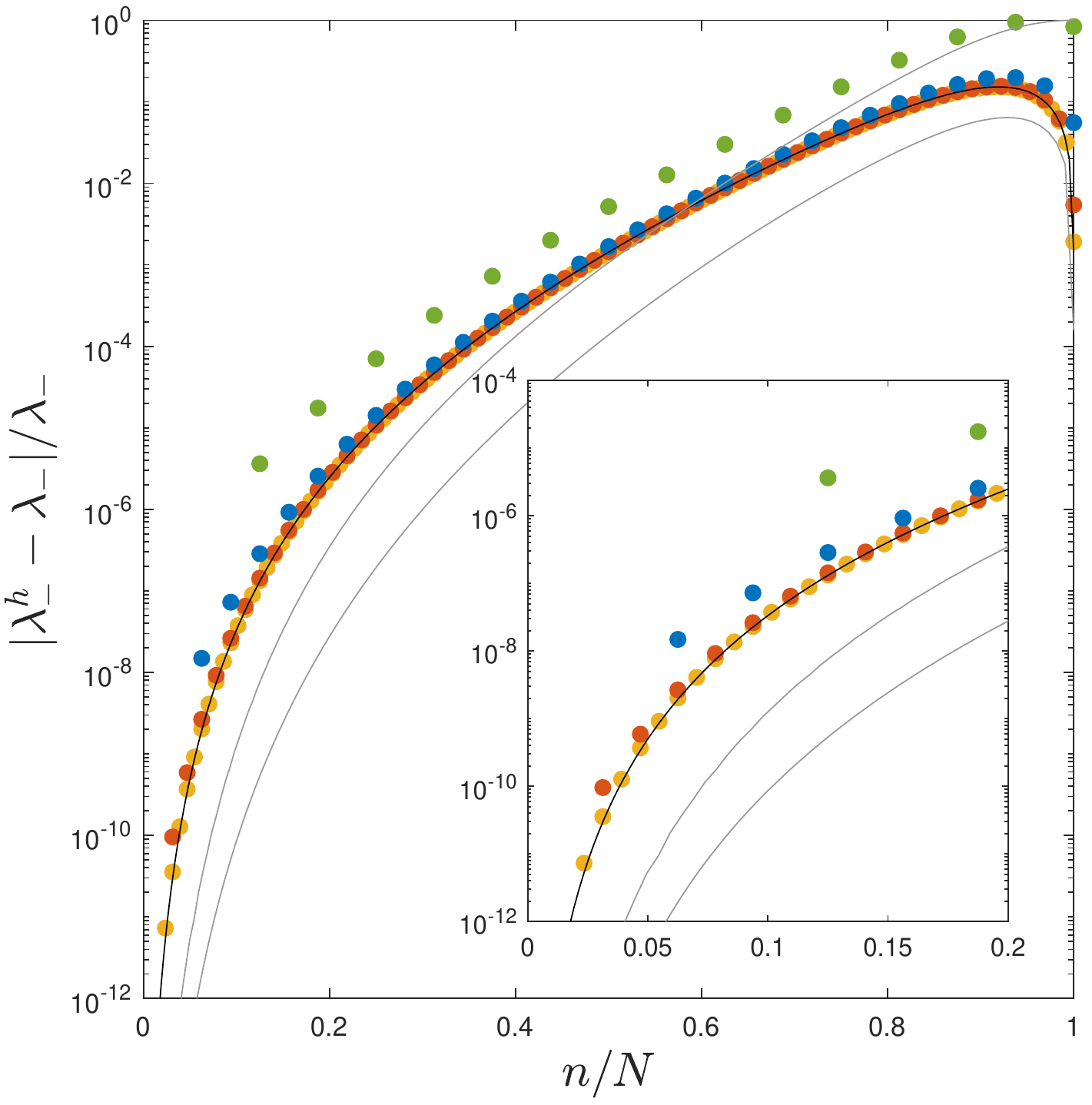}} \\
	\begin{tikzpicture}
		\filldraw[green1,line width=1pt] (2,0) circle (2pt);
		\filldraw[green1,line width=1pt] (2,0) node[right]{\scriptsize $N = 16$};	
		\filldraw[blue1,line width=1pt] (4,0) circle (2pt);
		\filldraw[blue1,line width=1pt] (4,0) node[right]{\scriptsize$N = 32$};	
		\filldraw[red1,line width=1pt] (6,0) circle (2pt);
		\filldraw[red1,line width=1pt] (6,0) node[right]{\scriptsize$N = 64$};	
		\filldraw[yellow1,line width=1pt] (8,0) circle (2pt);
		\filldraw[yellow1,line width=1pt] (8,0) node[right]{\scriptsize$N = 128$};
		\filldraw[black, line width=1pt] (9.9,0) -- (10.1,0) node[right]{\scriptsize $N \rightarrow \infty$};
	\end{tikzpicture}
	\caption{Normalized error in the eigenvalues $\lambda^h_{-}$ for the standard discretization obtained with cubics and quartics.}
	\label{fig:spectra_negative_p34_t0015}
\end{figure}

\begin{figure}[!h]
	\centering
	\subfloat[Standard formulation $p=3$]{\includegraphics[trim = 0cm 0cm 0cm 0cm, clip,width=0.45\textwidth]{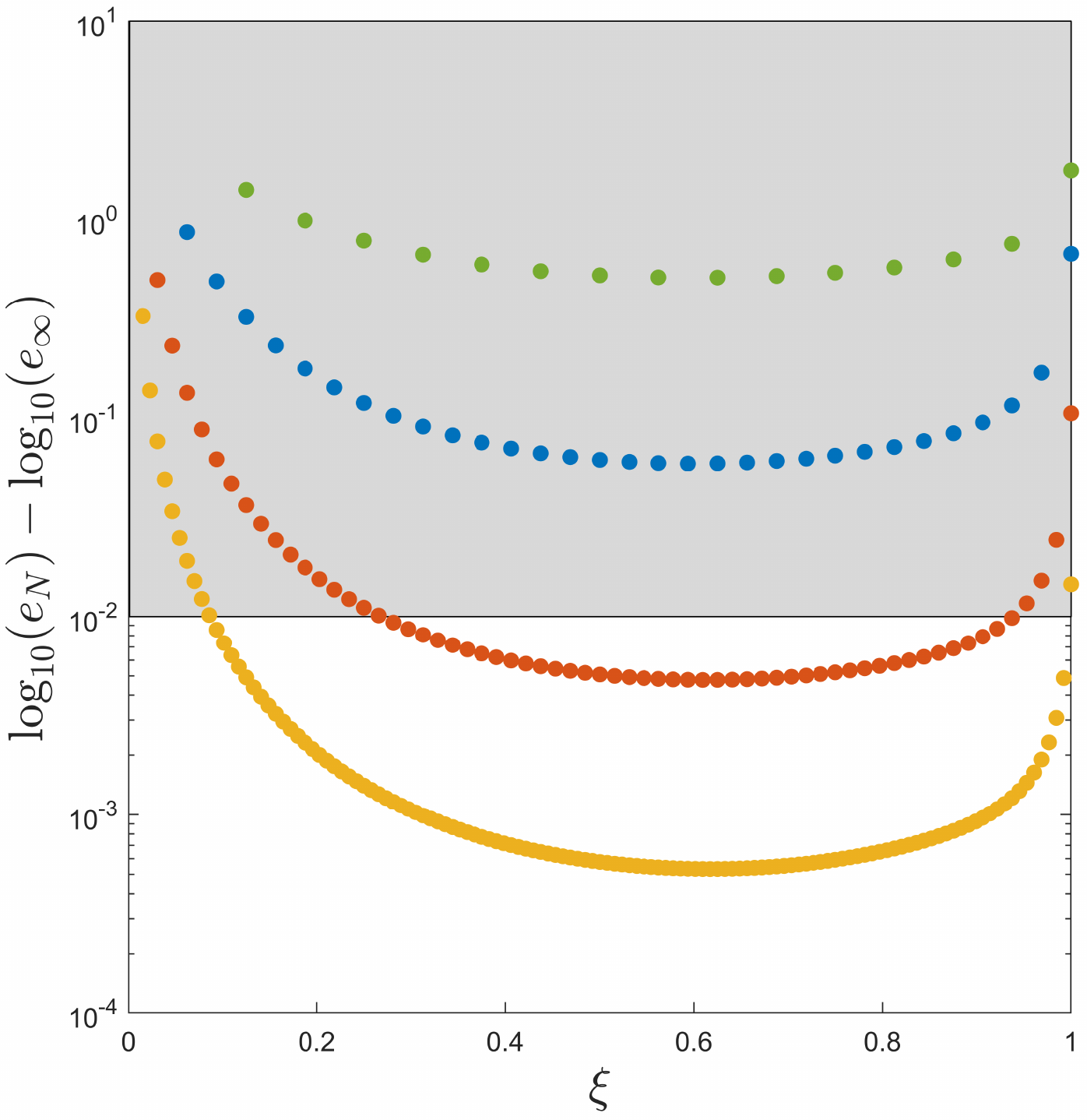}} \hspace{0.5cm}
	\subfloat[Standard formulation $p=4$]{\includegraphics[trim = 0cm 0cm 0cm 0cm, clip,width=0.45\textwidth]{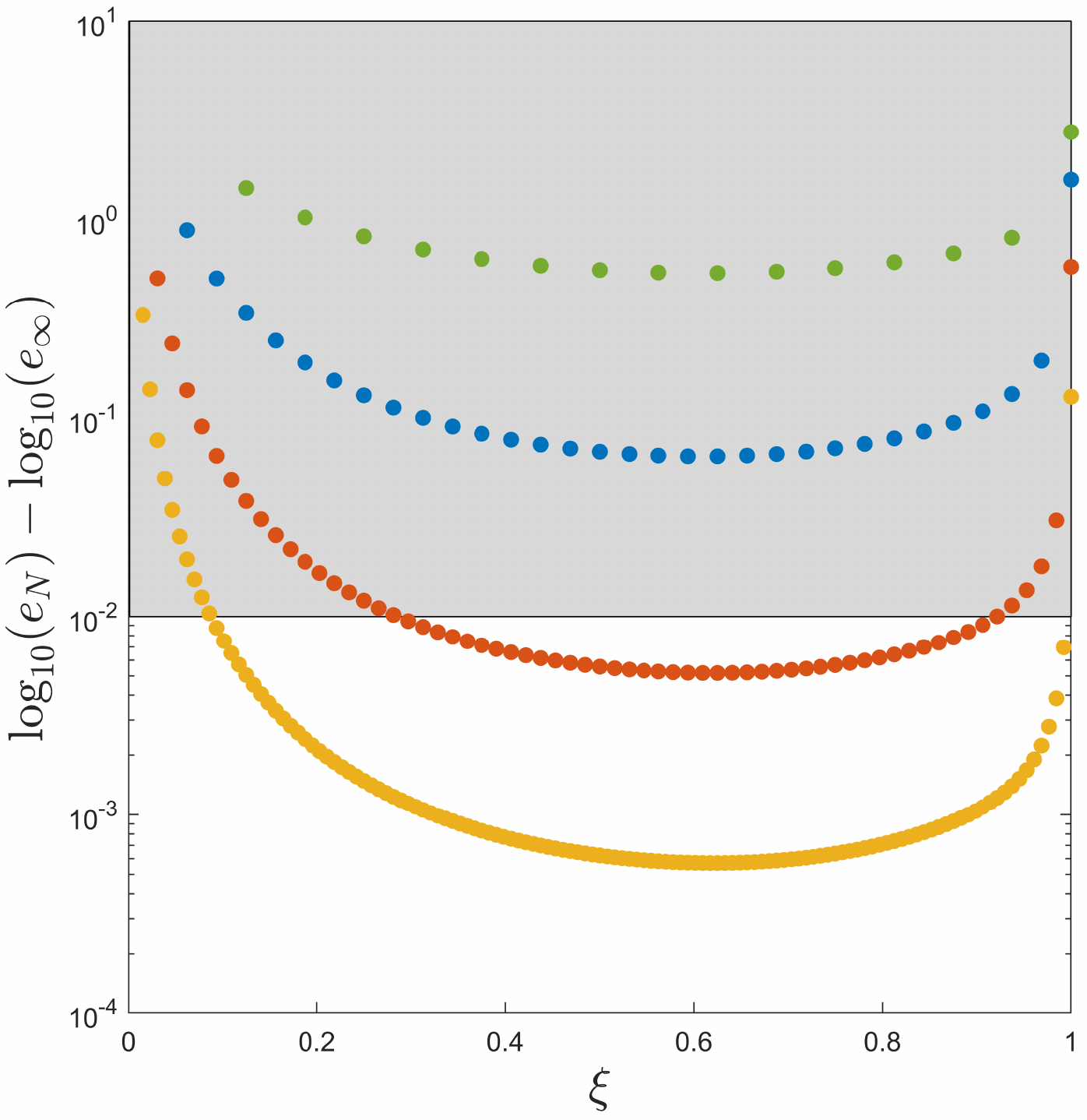}} \\
	\begin{tikzpicture}
		\filldraw[green1,line width=1pt] (2,0) circle (2pt);
		\filldraw[green1,line width=1pt] (2,0) node[right]{\scriptsize $N = 16$};	
		\filldraw[blue1,line width=1pt] (4,0) circle (2pt);
		\filldraw[blue1,line width=1pt] (4,0) node[right]{\scriptsize$N = 32$};	
		\filldraw[red1,line width=1pt] (6,0) circle (2pt);
		\filldraw[red1,line width=1pt] (6,0) node[right]{\scriptsize$N = 64$};	
		\filldraw[yellow1,line width=1pt] (8,0) circle (2pt);
		\filldraw[yellow1,line width=1pt] (8,0) node[right]{\scriptsize$N = 128$};
		\filldraw[black, line width=1pt] (9.9,0) -- (10.1,0) node[right]{\scriptsize $N \rightarrow \infty$};
	\end{tikzpicture}
	\caption{Identification of eigenvalues $\lambda^h_{-}$ that display membrane locking ($\varepsilon = 1 \cdot 10^{-2}$) in cubic and quartic discretizations.}
	\label{fig:locking_p34}
\end{figure}

Figure \ref{fig:spectra_negative_p34_t0015} and \ref{fig:locking_p34} illustrate the severity of locking of the standard formulation for cubics and quartics. Although the absolute accuracy improves dramatically with polynomial degree, the severity of membrane locking does not improve concomitantly. The degree of membrane locking is the same as for quadratics. Note also that the absolute accuracy of the mixed method again far exceeds that of the standard formulation, as may be observed from the asymptotic curves in Figure  \ref{fig:spectra_negative_p34_t0015}.

\subsection{Normalized error in mode amplitudes}
\label{sec:modeamplitudeerrors}
\begin{figure}[!h]
	\centering
	\subfloat[Standard formulation]{\includegraphics[trim = 0cm 0cm 0cm 0cm, clip,width=0.48\textwidth]{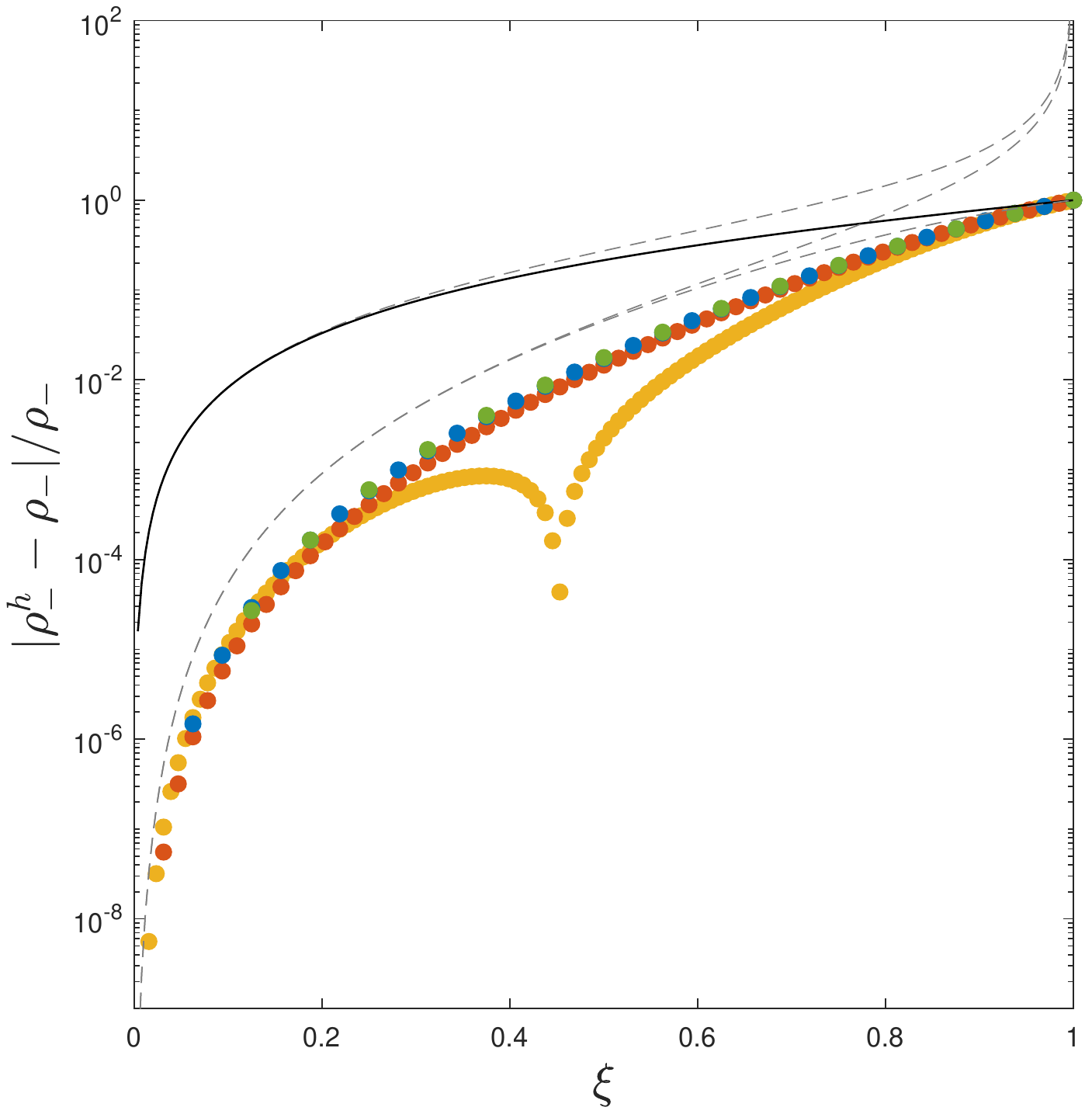}} \hspace{0.5cm}
	\subfloat[Mixed formulation]{\includegraphics[trim = 0cm 0cm 0cm 0cm, clip,width=0.48\textwidth]{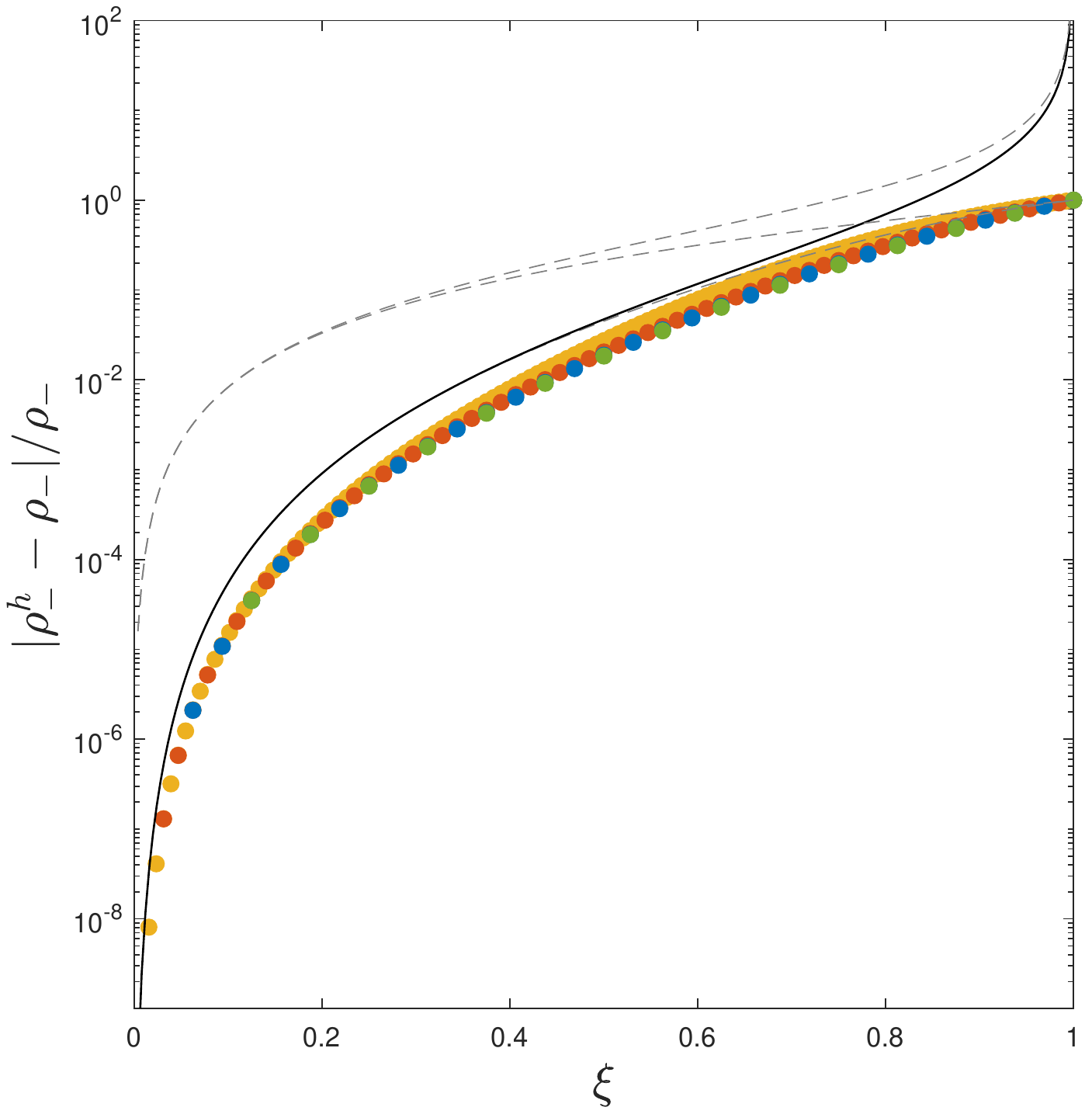}} \\
	\begin{tikzpicture}
		\filldraw[green1,line width=1pt] (2,0) circle (2pt);
		\filldraw[green1,line width=1pt] (2,0) node[right]{\scriptsize $N = 16$};	
		\filldraw[blue1,line width=1pt] (4,0) circle (2pt);
		\filldraw[blue1,line width=1pt] (4,0) node[right]{\scriptsize$N = 32$};	
		\filldraw[red1,line width=1pt] (6,0) circle (2pt);
		\filldraw[red1,line width=1pt] (6,0) node[right]{\scriptsize$N = 64$};	
		\filldraw[yellow1,line width=1pt] (8,0) circle (2pt);
		\filldraw[yellow1,line width=1pt] (8,0) node[right]{\scriptsize$N = 128$};
		\filldraw[black, line width=1pt] (9.9,0) -- (10.1,0) node[right]{\scriptsize $N \rightarrow \infty$};
	\end{tikzpicture}
	\caption{Normalized error in the mode amplitudes of the negative branch,  $\rho^h_{n-} = U_{n-}^h / W_{n-}^h$, obtained with (a) standard and (b) mixed quadratic discretizations.}
	\label{fig:normalized_amplitude_error}
\end{figure}

Figure \ref{fig:normalized_amplitude_error} illustrates the normalized error in the amplitudes of the eigenmodes corresponding to the negative branch. Both methods exhibit the same level of accuracy when applied to meshes with $2N=32, 64, 128$, and $256$ degrees of freedom, consistent with earlier findings reported in \cite{nguyen2022leveraging}. However, the asymptotic limit curves demonstrate distinct behavior. In relative terms, the approximations deteriorate with increasing refinement, especially noticeable in the standard formulation. This interesting phenomenon might be another manifestation of membrane locking, occurring specifically in the higher modes. It is important to note that the positive branch displays an identical pattern of behavior.

\subsection{Behaviour of normalized eigenvalue errors for fine discretizations}\label{sec:discretebranchswitching}
\begin{figure}[!h]
	\centering
	\subfloat[Standard formulation - $\lambda^h_{-}$]{\includegraphics[trim = 0cm 0cm 0cm 0cm, clip,width=0.48\textwidth]{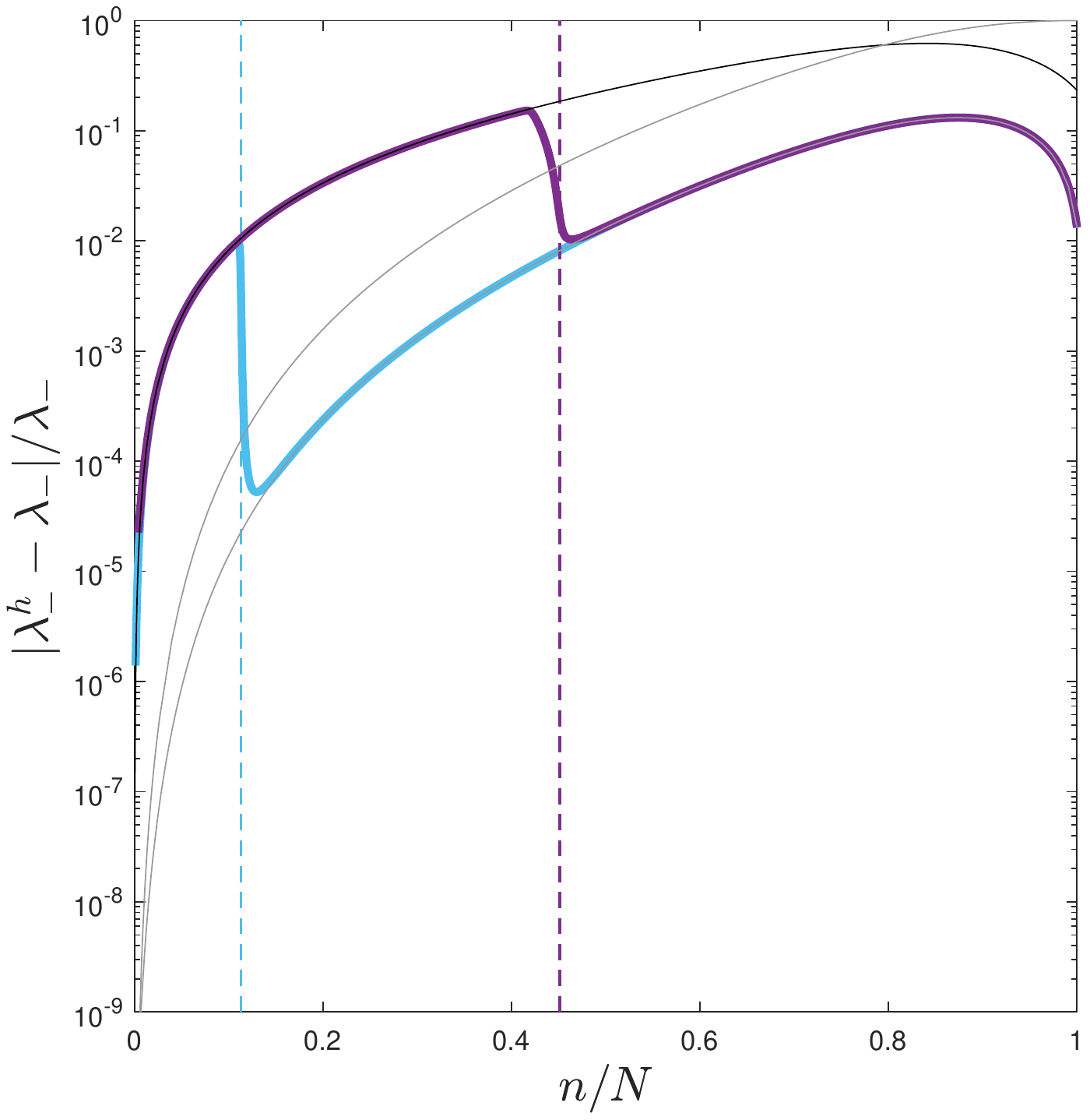}} \hspace{0.5cm}
	\subfloat[Standard formulation - $\lambda^h_{+}$]{\includegraphics[trim = 0cm 0cm 0cm 0cm, clip,width=0.48\textwidth]{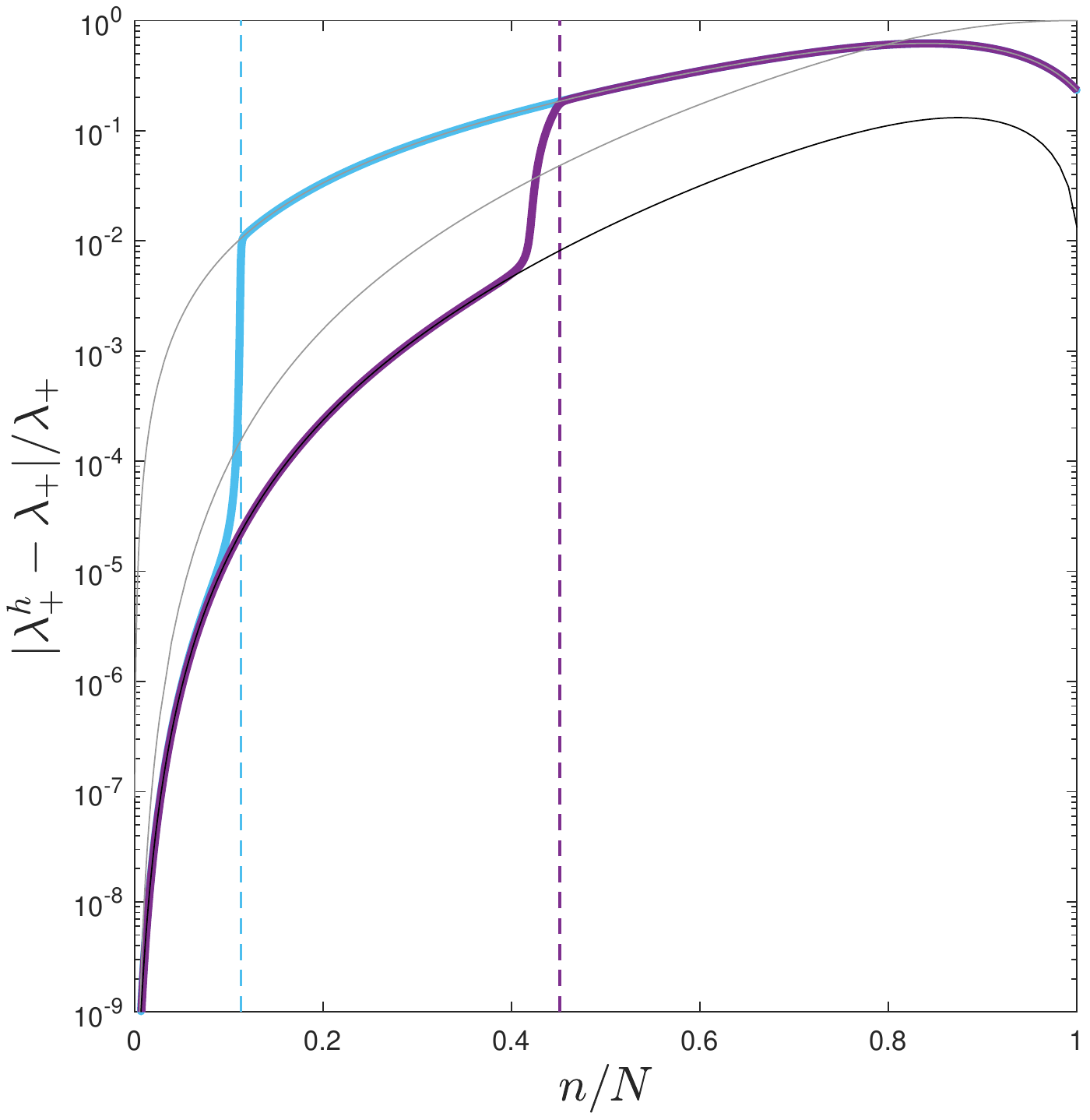}} \\
	\begin{tikzpicture}
		\filldraw[green1,line width=1pt] (2,0) circle (2pt);
		\filldraw[green1,line width=1pt] (2,0) node[right]{\scriptsize $N = 16$};	
		\filldraw[blue1,line width=1pt] (4,0) circle (2pt);
		\filldraw[blue1,line width=1pt] (4,0) node[right]{\scriptsize$N = 32$};	
		\filldraw[red1,line width=1pt] (6,0) circle (2pt);
		\filldraw[red1,line width=1pt] (6,0) node[right]{\scriptsize$N = 64$};	
		\filldraw[yellow1,line width=1pt] (8,0) circle (2pt);
		\filldraw[yellow1,line width=1pt] (8,0) node[right]{\scriptsize$N = 128$};
		\filldraw[black, line width=1pt] (9.9,0) -- (10.1,0) node[right]{\scriptsize $N \rightarrow \infty$};
	\end{tikzpicture}
	\caption{Normalized eigenvalue error obtained with quadratic, fine discretizations of the standard formulation.}
	\label{fig:normalized_spectra_error_fine_discretizations}
\end{figure}
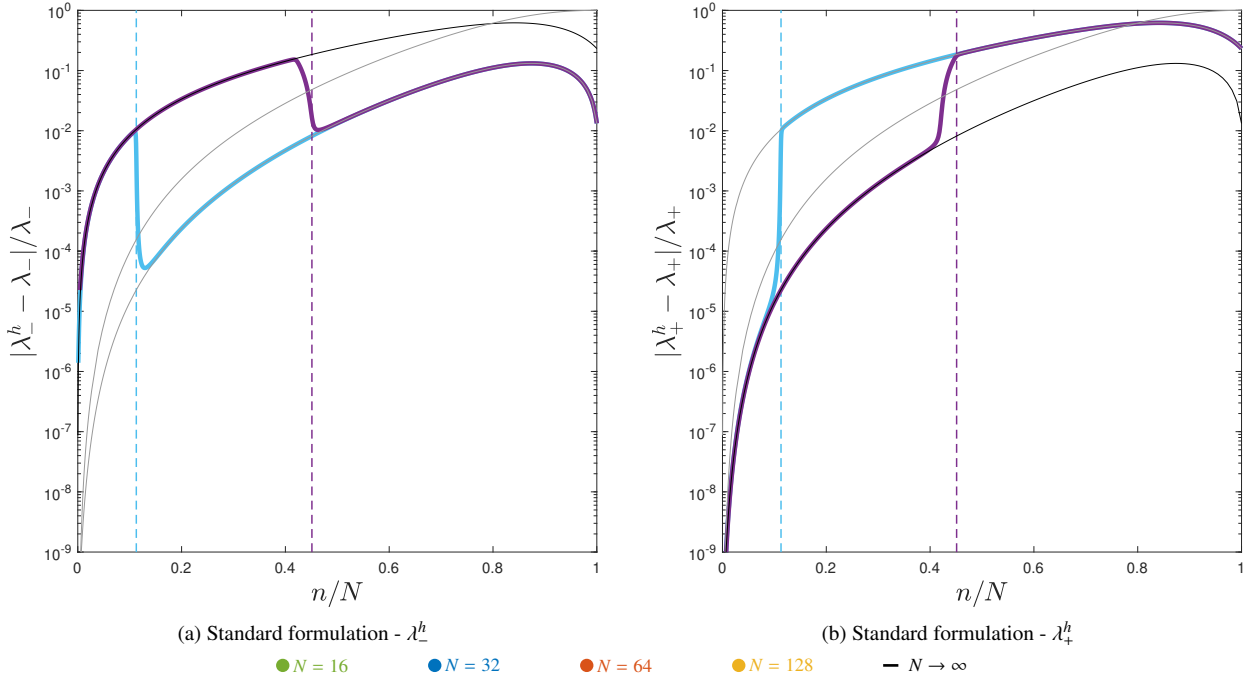

In addition to our investigation into membrane locking occurring in coarse to medium refined discretizations, we present interesting and previously unobserved findings that occur in fine discretizations of the circular Euler-Bernoulli ring, see Figure \ref{fig:normalized_spectra_error_fine_discretizations}. Notably, for discretizations where $N > \hat{n}$, the curves of normalized eigenvalue errors deviate from their asymptotic pattern. Although initially surprising, this behaviour can be anticipated upon examining the relative amplitudes presented in Figures \ref{fig:analytical_eigenmodes} and \ref{fig:discrete_eigenvalues_and_amplitudes}. We should recall that for $n < \hat{n}$, the modes corresponding to $\lambda_{n+}$ are membrane-dominated ($U_{n+}/W_{n+} > 1$), whereas the modes corresponding to $\lambda_{n-}$ are dominated by bending effects ($U_{n-} / W_{n-} < 1$). However, this behaviour is reversed for $n > \hat{n}$. While we specifically demonstrate this behavior using quadratic discretizations of the standard formulation, it holds true for polynomial degrees of any order and is also applicable to the mixed formulation.

\section{Mathematical analysis \label{sec:6}}
In this section, we present a detailed mathematical analysis of the error between the discrete eigenvalues and their exact counterparts as a function of the relative mode number, the mesh size, and the relevant physical parameters. This relative error is measured by the difference between $\lambda^h_{n\pm}$ and $\lambda_{n\pm}$ normalized by $\lambda_{n\pm}$ with $n=\xi N=\xi\pi/h$ (see \eqref{eq:relerrordef} in Section \ref{sec:lockingcriterion}). Our primary focus is on the standard formulation, where membrane locking is more pronounced. However, our analysis is comprehensive and applies to general polynomial degrees and could be adapted to the mixed method as well. In light of our characterization of membrane locking, we are particularly interested in studying the behavior of the relative error as $h$ approaches zero while keeping $\xi$ fixed. Additionally, we extend our investigation to explore the relationship between eigenvalue error and normalized thickness, thereby providing a broader understanding of the phenomenon.

\subsection{Relative eigenvalue error as a function of the normalized mode number}
\label{sec:exactrelerr}
We will focus on the expression
\begin{equation}\label{eq:relerrorgen}
	\frac{\lambda^h_{n\pm}-\lambda_{n\pm}}{\lambda_{n\pm}}
		=\frac{\lambda^h_{n\pm}}{\lambda_{n\pm}}-1
			=\frac{K_n^h}{K_n}
				\left[\frac{1\pm\sqrt{1-M_n^h/(K_n^h)^2}}{1\pm\sqrt{1-M_n/(K_n)^2}}\right]-1\,,
\end{equation}
where we take the values of $\lambda_{n\pm}$ and  $\lambda_{n\pm}^h$ from \eqref{eq:lambaanalyticKL} and \eqref{eq:lambdadiscreteKL}. 

First, note that rewriting $K_n$ in \eqref{eq:exactKML} with $n=\tfrac{\xi\pi}{h}$ yields \eqref{eq:Kexact}. On the other hand, a careful analysis of \ref{app:standard} yields an analogous expression for the discrete case \eqref{eq:Khexact},
\begin{subequations}\label{eq:KandKhexact}
\begin{align}
	K_n  	&= \frac{1}{2}\frac{1}{R^2}\frac{1}{h^4}\left(\alpha_0+\alpha_2h^2+h^4\right),  \label{eq:Kexact} \\
	K_n^h &=\frac{1}{2}\frac{1}{R^2}\frac{1}{h^4}\left(\alpha_0^h+\alpha_2^h h^2+h^4\right), \label{eq:Khexact}
\end{align}
\end{subequations}
with
\begin{subequations}
\begin{align}
	\alpha_0 &= \frac{\bar{t}^2}{12}(\xi\pi)^4,	&  \alpha_0^h &= \frac{\bar{t}^2}{12}\frac{\tilde{\mat{K}}_b^{22}}{\tilde{\mat{M}}}, & \\
	\alpha_2 &=\left(1+\frac{\bar{t}^2}{12}\right)(\xi\pi)^2, &  \alpha_2^h &=\left(1+\frac{\bar{t}^2}{12}\right)\frac{\tilde{\mat{K}}_b^{11}}{\tilde{\mat{M}}}. &
\end{align}
\end{subequations}
Here, the $\tilde{\mat{K}}_b^{ij}$ and $\tilde{\mat{K}}_m^{ij}$ represent the eigenvalues of the ``normalized'' circulant matrices, as presented in \eqref{eq:KbKmstandardform}, by which we mean that matrices have been stripped of any dependence on $\bar{t}$, $R$, and $h$. Similarly, $\tilde{\mat{M}}$ is the eigenvalue of the normalized circulant symmetric matrix described by \eqref{eq:MassMtandardform}. More precisely, these are the eigenvalues corresponding to the matrices reported in Table \ref{tab:circulant_matrices}.

As highlighted in Section \ref{sec:EVcirculant}, these eigenvalues can be analytically computed and expressed as weighted sums of sines and cosines of integer multiples of $\xi\pi=nh=n\tfrac{2\pi}{2N}$ (see equations \eqref{eq:EVcirculantsymmetric} and \eqref{eq:EVcirculantantisymmetric}). This key observation implies that the eigenvalues are solely determined by the normalized mode number $\xi$ and are independent of other parameters. This point holds significant importance in understanding the behavior of the eigenvalues as functions of the normalized mode number.

A similar approach to $M_n$ and $M_n^h$ yield the following expressions,
\begin{subequations}\label{eq:MandMhexact}
\begin{align}
	M_n &= \frac{1}{4}\frac{1}{R^4}\frac{1}{h^6}\left(a_0+a_2h^2+a_4h^4\right), & \\
M_n^h &=\frac{1}{4}\frac{1}{R^4}\frac{1}{h^6}\left(a_0^h+a_2^hh^2+a_4^hh^4\right), &
\end{align}
\end{subequations}
with
\begin{subequations}
\begin{align}
	a_0 &= \frac{4\bar{t}^2}{12} \, (\xi\pi)^6, &
	a^h_0 &= \frac{4\bar{t}^2}{12}\frac{1}{\tilde{\mat{M}}^2} \,  \Bigg[\bigg(1+\frac{\bar{t}^2}{12}\bigg)\tilde{\mat{K}}_b^{11}\tilde{\mat{K}}_b^{22}-\frac{\bar{t}^2}{12}\big|\tilde{\mat{K}}_b^{12}\big|^2\Bigg], & \\
	a_2 &= -\frac{4\bar{t}^2}{12} \, 2(\xi\pi)^4, & 
	a^h_2 &= - \frac{4\bar{t}^2}{12}\frac{1}{\tilde{\mat{M}}^2} \, 2\tilde{\mat{K}}_m^{12}(\tilde{\mat{K}}_b^{12})^{\ast}, & \\
	a_4 &= \frac{4\bar{t}^2}{12} \,  (\xi\pi)^2, & 
	a^h_4 &= \frac{4\bar{t}^2}{12}\frac{1}{\tilde{\mat{M}}^2} \, \frac{12}{\bar{t}^2}\Bigg[\Bigg(1+\frac{\bar{t}^2}{12}\Bigg)\tilde{\mat{K}}_b^{11}\tilde{\mat{K}}_m^{22}-\big|\tilde{\mat{K}}_m^{12}\big|^2\Bigg].   
\end{align}
\end{subequations}

Within this context we note that $\tilde{\mat{K}}_m^{12}(\tilde{\mat{K}}_b^{12})^{\ast}$ (with $(\cdot)^\ast$ denoting complex conjugation) is a real number in view of the fact that they are purely imaginary eigenvalues coming from skew-symmetric matrices.

Meanwhile, it follows directly from \eqref{eq:Kexact} and \eqref{eq:Khexact} that
\begin{subequations}
\begin{align}
	&\begin{aligned}
		K_n^2&=\frac{1}{4}\frac{1}{R^4}\frac{1}{h^8}\Big(\alpha_0^2+2\alpha_0\alpha_2h^2
			+\big(2\alpha_0+\alpha_2^2\big)h^4+2\alpha_2h^6+h^8\Big)\\
			&=\frac{1}{4}\frac{1}{R^4}\frac{1}{h^8}\Big(b_0+b_2h^2+b_4h^4+b_6h^6+h^8\Big)
	\end{aligned}\\
	&\!\!\!\!\!\!\!\begin{aligned}
		(K_n^h)^2&=\frac{1}{4}\frac{1}{R^4}\frac{1}{h^8}\Big((\alpha_0^h)^2+2\alpha_0^h\alpha_2^hh^2
			+\big(2\alpha_0^h+(\alpha_2^h)^2\big)h^4+2\alpha_2^hh^6+h^8\Big)\\
			&=\frac{1}{4}\frac{1}{R^4}\frac{1}{h^8}\Big(b_0^h+b_2^hh^2+b_4^hh^4+b_6^hh^6+h^8\Big)\,.
	\end{aligned}
\end{align}
\end{subequations}

Putting everything together gives,
\begin{equation}\label{eq:exactrelerror}
	\frac{\lambda^h_{n\pm}}{\lambda_{n\pm}}
		=\frac{K_n^h}{K_n}\left[\frac{1\pm\sqrt{1-\frac{M_n^h}{(K_n^h)^2}}}
			{1\pm\sqrt{1-\frac{M_n}{K_n^2}}}\right]
		=\frac{\alpha_0^h+\alpha_2^hh^2+ h^4}{\alpha_0+\alpha_2h^2+h^4}
			\left[\frac{1\pm\sqrt{1-h^2\frac{a_0^h+a_2^hh^2+a_4^hh^4}
				{b_0^h+b_2^hh^2+b_4^hh^4+b_6^hh^6+h^8}}}
					{1\pm\sqrt{1-h^2\frac{a_0+a_2h^2+a_4h^4}
						{b_0+b_2h^2+b_4h^4+b_6h^6+h^8}}}\right]\,.
\end{equation}
The main observation at this point is that there is no dependence on $R$ outside of the normalized thickness $\bar{t}=\tfrac{t}{R}$ (more generally outside of $\bar{I}=\tfrac{I}{AR^2}=\frac{t^2}{12R^2}$).

\subsection{Relative eigenvalue error in the asymptotic limit of refinement} \label{sec:relerrasymtpoticlimit}
We are particularly concerned with studying the relative error with respect to the normalized mode number $\xi\in(0,1]$ in the asymptotic refinement limit, where $h$ approaches zero. As previously mentioned, we hypothesize that membrane locking is responsible for the disparity between the spectrum computed at a fixed $h$ and the spectrum in this asymptotic limit.

The expression in \eqref{eq:exactrelerror} lends itself very well to taking this limit in both the positive and negative branches of eigenvalues. For the positive branch, taking the limit $h\to0$ gives
\begin{equation}\label{eq:aymptoticbendingdominated}
	\lim_{h\to0}\frac{\lambda^h_{n+}}{\lambda_{n+}}-1
		=\frac{\alpha_0^h}{\alpha_0}\left[\frac{1+\sqrt{1}}{1+\sqrt{1}}\right]-1
		=\frac{1}{(\xi\pi)^4}\frac{\tilde{\mat{K}}_b^{22}}{\tilde{\mat{M}}}-1\,,
\end{equation}
where the values of $\alpha_0$ and $\alpha_0^h$ were taken from \eqref{eq:Kexact} and \eqref{eq:Khexact}. It is important to observe that the normalized thickness does not affect this expression as it cancels out. In simpler terms, this expression solely relies on $\xi$ and the polynomial order $p$, which determine $\tilde{\mat{K}}_b^{22}$ and $\tilde{\mat{M}}$.

With that being said, it is crucial to further examine this expression. For the case of $p=2$, we can compute the circulant eigenvalues of $\tilde{\mat{K}}_b^{22}$ and $\tilde{\mat{M}}$ using \eqref{eq:EVcirculantsymmetric} along with the values provided in Table \ref{tab:circulant_matrices}. The resulting expression is given by:
\begin{equation}\label{eq:aymptoticbendingdominatedp2}
	\lim_{h\to0}\frac{\lambda^h_{n+}}{\lambda_{n+}}-1
		=\frac{60}{(\xi\pi)^4} \cdot \frac{6 - 8 \cos{(\xi\pi)} + 2\cos{(2 \xi\pi)}}{33 + 26 \cos{(\xi\pi)} + \cos{(2\xi\pi)} } - 1.
\end{equation}
Remarkably, this expression exactly matches \eqref{eq:EVfourthorder}, which represents the relative error resulting from solving the isolated bending response. Notably, this correspondence holds generally true for any value of $p$. In essence, it implies that in the limit of mesh refinement, the relative error of the positive branch coincides with that observed in the isolated bending response!

To analyze the negative branch of eigenvalues, we employ the Taylor series to avoid encountering a limit of the form $\frac{0}{0}$. However, in the end, the following expression holds in that limit:
\begin{equation}\label{eq:aymptoticmembranedominated}
	\lim_{h\to0}\frac{\lambda^h_{n-}}{\lambda_{n-}}-1
		=\frac{\alpha_0^h}{\alpha_0}\frac{a_0^h}{a_0}\frac{b_0}{b_0^h}-1
		=\frac{1}{(\xi\pi)^2}\left[\left(1+\frac{\bar{t}^2}{12}\right)
			\frac{\tilde{\mat{K}}_m^{11}}{\tilde{\mat{M}}}
				-\frac{\bar{t}^2}{12}\frac{\big|\tilde{\mat{K}}_b^{12}\big|^2}
					{\tilde{\mat{M}}\tilde{\mat{K}}_b^{22}}\right]-1\,.
\end{equation}
Here, we utilized the fact that $\tilde{\mat{K}}_m^{11}=\tilde{\mat{K}}_b^{11}$. It is evident that this expression is dependent on $\bar{t}$, and as $\bar{t}\to0$, the relative error becomes $\frac{1}{(\xi\pi)^2}\frac{\tilde{\mat{K}}_m^{11}}{\tilde{\mat{M}}}-1$, which precisely corresponds to the isolated membrane response, regardless of the value of $p$. Indeed, for $p=2$, it is not difficult to notice that the above expression coincides with \eqref{eq:EVsecondorder} in the limit of $\bar{t}\to0$.

\subsection{Relation to classical eigenvalue error analysis} \label{sec:standarderror}
Typically, convergence in eigenvalue problems is studied by fixing a mode number $n$ and refining the mesh. Interestingly, our analysis involving the normalized mode number $\xi=\frac{n}{N}=\frac{nh}{\pi}$ incorporates this aspect naturally. This is due to the fact that for a fixed mode number $n$, as we decrease the mesh size $h$ (equivalently, increase the number of elements $N$), the value $nh=\xi\pi$ becomes small. In other words, when we focus on small values of $\xi$, we essentially recover the behavior of the relative error at small $h$ for any fixed $n$. This implies that by examining the relative error curves in the asymptotic limit of refinement, as derived in \eqref{eq:aymptoticbendingdominated} and \eqref{eq:aymptoticmembranedominated}, we can gain insights into the behavior of the relative error for any fixed $n$.

To illustrate this point, let us consider the case of $p=2$. First,  examine the bending-dominated asymptotic limit of refinement described by \eqref{eq:aymptoticbendingdominatedp2}. When $\xi\pi$ is small (i.e., $\xi\pi \ll 1$), the expression should approach zero. However, upon closer inspection, we notice that it exhibits an \emph{inverse} power of $\xi\pi$, which may initially seem counterintuitive as it suggests a large value for small $\xi$. This apparent contradiction can be resolved by analyzing the expression more carefully using its Taylor series expansion, which allows us to account for the cancellation of these inverse powers. This technique is analogous to the methods employed in finite differences \cite[\S3.4]{leveque2007finite}. By performing the necessary calculations, we obtain the following expression:
\begin{equation}
	\lim_{h\to0}\frac{\lambda^h_{n+}}{\lambda_{n+}}-1
		=\frac{60}{(\xi\pi)^4}\frac{6-8\cos(\xi\pi)+2\cos(2\xi\pi)}
			{33+26\cos(\xi\pi)+\cos(2\xi\pi)}-1
			=\frac{1}{12}(\xi\pi)^2+\frac{1}{240}(\xi\pi)^4+\mathcal{O}\big((\xi\pi)^6\big)\,.
\end{equation}
Therefore, when considering a fixed $n$ and small enough $h$, the relative error exhibits a behavior of $\mathcal{O}\big(h^2\big)$. This result aligns with expectations, as the anticipated convergence rate for eigenvalues in a fourth-order problem is $\mathcal{O}\big(h^{2(p-1)}\big)$. It can be easily verified that this behavior holds true for other values of $p$ as well.

Moving forward, let's examine the membrane-dominated asymptotic limit of refinement for the case of $p=2$, as presented in \eqref{eq:aymptoticmembranedominated}. To analyze this expression, we need to consider the Taylor series expansions of its two components:
\begin{equation}
	\begin{aligned}
		\frac{1}{(\xi\pi)^2}\frac{\tilde{\mat{K}}_m^{11}}{\tilde{\mat{M}}}
		   &=\frac{60}{3}\frac{3-2\cos(\xi\pi)-\cos(2\xi\pi)}{33+26\cos(\xi\pi)+\cos(2\xi\pi)}\\
		 		&=1+\frac{1}{720}(\xi\pi)^4+\frac{1}{3360}(\xi\pi)^6
		 			+\mathcal{O}\big((\xi\pi)^8\big)\\
		\frac{1}{(\xi\pi)^2}\frac{\big|\tilde{\mat{K}}_b^{12}\big|^2}
			{\tilde{\mat{M}}\tilde{\mat{K}}_b^{22}}
				&=\frac{60\big(2\sin(\xi\pi)-\sin(2\xi\pi)\big)^2}
					{\big(6-8\cos(\xi\pi)+2\cos(2\xi\pi)\big)
						\big(33+26\cos(\xi\pi)+\cos(2\xi\pi)\big)}\\
			&=1-\frac{1}{12}(\xi\pi)^2-\frac{1}{180}(\xi\pi)^4
				+\mathcal{O}\big((\xi\pi)^6\big)\,.
	\end{aligned}
\end{equation}
It is important to note that the isolated membrane response, represented by $\frac{1}{(\xi\pi)^2}\frac{\tilde{\mat{K}}_m^{11}}{\tilde{\mat{M}}}$, exhibits eigenvalue convergence of $\mathcal{O}\big(h^4\big)$. This convergence rate coincides with the expected rate for second-order problems, which is $\mathcal{O}\big(h^{2p}\big)$. On the other hand, the bending component, given by $\frac{1}{(\xi\pi)^2}\frac{|\tilde{\mat{K}}_b^{12}|^2}{\tilde{\mat{M}}\tilde{\mat{K}}_b^{22}}$, converges at a slower rate of $\mathcal{O}\big(h^2\big)$. Consequently, the overall convergence rate of the relative error is determined by the bending component, resulting in a convergence rate of $\mathcal{O}\big(h^2\big)$ for small $h$. For arbitrary $p$, the membrane-dominated eigenvalues will  converge at a rate of $\mathcal{O}\big(h^{2(p-1)}\big)$.

The inclusion of the relative error curves as functions of $\xi$ provides additional insights. We can determine the value of $\xi^*$ below which the relative error behaves consistently with $\mathcal{O}\big(h^{2(p-1)}\big)$. This value serves as a measure of the preasymptotic domain's size, indicating the threshold at which the eigenvalues begin to converge at the expected rate for a given $n$. These $\xi^*$ values (or equivalently, $h^*=\tfrac{\xi^*\pi}{n}$) may vary depending on the value of $p$. From the curves of the asymptotic limit of refinement shown in Figure \ref{fig:normalized_spectra_error}, we can estimate that the convergence for both branches of eigenvalues begins at $\xi^*<0.1$. It is also noteworthy that the mixed formulation exhibits a distinct asymptotic curve for the bending-dominated eigenvalues, suggesting an earlier onset of convergence (i.e., at a larger value of $\xi$).

\subsection{Considerations for comparing relative errors}
To accurately assess the relative errors, it is crucial to consider a subtle aspect. When $h$ takes on moderate values such that $N < \hat{n}$ (i.e., the number of elements in the discretization is less than the transition mode $\hat{n}$), the eigenvalues in the positive branch correspond to a membrane-dominated response, as depicted in Figure \ref{fig:analytical_eigenmodes}. However, as $h$ approaches zero and the number of elements $N$ becomes much larger than the transition mode $\hat{n}$ (i.e., $N \gg \hat{n}$), the positive branch of eigenvalues predominantly consists of bending-dominated eigenvalues (see Figure \ref{fig:normalized_spectra_error_fine_discretizations}). This aligns with the observation made earlier regarding the curve derived in \eqref{eq:aymptoticbendingdominated}, which is associated with the isolated bending response.

To effectively analyze the discrepancy in the relative error of an eigenvalue at a finite value of $h$, it is essential to compare it against the asymptotic limit of refinement associated with the same \emph{class of physical response}, rather than focusing solely on the same eigenvalue branch. In particular, at $\bar{t}=0.01$, the first 100 eigenvalues corresponding to the positive branch exhibit membrane-dominated behavior (see Figure \ref{fig:analytical_eigenmodes}). When comparing their relative errors, we should consider the membrane-dominated asymptotic limit, which corresponds to the asymptotic limit of the negative branch of eigenvalues. Similar considerations apply to the eigenvalues in the negative branch.

Lastly, it is worth recalling that the transition mode $\hat{n}$ depends solely on the normalized thickness $\bar{t}$, as discussed in Section \ref{sec:analyticalEVEFs}. This implies that the behavior of the eigenvalues and their dominance by membrane or bending responses are determined by a combination of mesh refinement (controlled by $h$, see Figure \ref{fig:normalized_spectra_error_fine_discretizations}) and the normalized thickness $\bar{t}$. By considering the interplay between $h$ and $\bar{t}$, we can accurately interpret and compare eigenvalues with the correct limit curves. Therefore, in the locking criterion discussed in Section \ref{sec:lockingcriterion}, which we propose as a means of determining whether a particular method is affected by membrane locking, the comparison being made might span across eigenvalue branches. This detail should be taken into account when examining the definition.

\subsection{Analysis of the physical parameters}
Most physical parameters, such as Young's modulus, density, and area, have a multiplicative effect and therefore do not influence the behavior of either the analytical or discrete solutions. However, the thickness $t$ and radius of the ring $R$ do have a relative effect on certain terms. In fact, they can be consolidated into a single parameter, the normalized thickness $\bar{t} = \frac{t}{R}$ as defined previously, or alternatively, the normalized moment of inertia $\bar{I} = \frac{I}{AR^2}$, where $A$ represents the cross-sectional area of the ring. To see the equivalence, consider that for a rectangular cross-section, $\bar{I} = \frac{\bar{t}^2}{12}$, and $\bar{I}\to0$ if and only if $\bar{t}\to0$. To simplify the resulting formulas as much as possible, we will consistently use the parameter $\bar{I}$ throughout this section. However, it should be noted that $\bar{I}$ is merely proportional to $\bar{t}^2$.

The idea is to proceed similarly to Section \ref{sec:exactrelerr}, and obtain an expression analogous to \eqref{eq:exactrelerror}. Reorganizing the terms in \eqref{eq:KandKhexact} in terms of $\bar{I}$ yields
\begin{align}
	K_n=\frac{1}{2}\frac{1}{R^2}\frac{1}{h^4}\left(\beta_0+\beta_1\bar{I}\right),
	\qquad\qquad\quad
	K_n^h=\frac{1}{2}\frac{1}{R^2}\frac{1}{h^4}\left(\beta_0^h+\beta_1^h\bar{I}\right), \label{eq:K_tbar_dependence}
\end{align}
with 
\begin{subequations}
\begin{align}
	\beta_0 &= h^2(\xi\pi)^2+h^4, &
	\beta_0^h &= h^2\frac{\tilde{\mat{K}}_b^{11}}{\tilde{\mat{M}}}+h^4, & \\
	\beta_1 &= (\xi\pi)^4+h^2(\xi\pi)^2, & 
	\beta_1^h &= \frac{1}{\tilde{\mat{M}}}
		\left(\tilde{\mat{K}}_b^{22}+h^2\tilde{\mat{K}}_b^{11}\right).  &    
\end{align}
\end{subequations}

Proceeding similarly with \eqref{eq:MandMhexact} gives
\begin{align}\label{eq:MnMnhIbar}
	M_n=\frac{1}{4}\frac{1}{R^4}\frac{1}{h^8}\left(\tilde{a}_0 + \tilde{a}_1\bar{I} + \tilde{a}_2\bar{I}^2\right), 
	\qquad\qquad\quad
	M_n^h=\frac{1}{4}\frac{1}{R^4}\frac{1}{h^8}
		\left(\tilde{a}_0^h+\tilde{a}_1^h\bar{I}+\tilde{a}_2^h\bar{I}^2\right),
\end{align}
with
\begin{subequations}\label{eq:acoeffsIbar}
\begin{align}
	\tilde{a}_0   &=  0, &
	\tilde{a}_0^h &= \frac{4h^2}{\tilde{\mat{M}}^2}
		h^4\left[\tilde{\mat{K}}_m^{11}\tilde{\mat{K}}_m^{22}
			-\big|\tilde{\mat{K}}_m^{12}\big|^2\right], & \\
	\tilde{a}_1 &= 4h^2\left((\xi\pi)^6-2h^2(\xi\pi)^4+h^4(\xi\pi)^2\right), & 
	\tilde{a}_1^h &= \frac{4h^2}{\tilde{\mat{M}}^2}
		\left[\tilde{\mat{K}}_b^{11}\tilde{\mat{K}}_b^{22}
			-2h^2\tilde{\mat{K}}_m^{12}(\tilde{\mat{K}}_b^{12})^*
				+h^4\tilde{\mat{K}}_m^{11}\tilde{\mat{K}}_m^{22}\right],  &  \\
	 \tilde{a}_2 &= 0, &
	\tilde{a}_2^h &= \frac{4h^2}{\tilde{\mat{M}}^2}
		\left[\tilde{\mat{K}}_b^{11}\tilde{\mat{K}}_b^{22}
			-\big|\tilde{\mat{K}}_b^{12}\big|^2\right]. & 
\end{align}
\end{subequations}

Putting everything together yields the following exact expression for the relative magnitude of the eigenvalues,
\begin{equation}\label{eq:exactrelerrorIbar}
	\frac{\lambda^h_{n\pm}}{\lambda_{n\pm}}
		=\frac{K_n^h}{K_n}\left[\frac{1\pm\sqrt{1-\frac{M_n^h}{(K_n^h)^2}}}
			{1\pm\sqrt{1-\frac{M_n}{K_n^2}}}\right]
		=\frac{\beta_0^h+\beta_1^h\bar{I}}{\beta_0+\beta_1\bar{I}}
			\left[\frac{1\pm\sqrt{1-
				\frac{\tilde{a}_0^h+\tilde{a}_1^h\bar{I}+\tilde{a}_2^h\bar{I}^2}
					{(\beta_0^h)^2+2\beta_0^h\beta_1^h\bar{I}+(\beta_1^h)^2\bar{I}^2}}}
					{1\pm\sqrt{1-\frac{\tilde{a}_1\bar{I}}
						{(\beta_0)^2+2\beta_0\beta_1\bar{I}+(\beta_1)^2\bar{I}^2}}}\right]\,.
\end{equation}

To gain insights into the influence of $\bar{I}$ in \eqref{eq:exactrelerrorIbar}, one can employ a Taylor series expansion centered around $\bar{I}=0$. However, more interesting conclusions can be drawn by directly examining \eqref{eq:exactrelerrorIbar} as $\bar{I}$ approaches zero. Notably, the behavior of the eigenvalues in both branches exhibits distinct characteristics.

Regarding the positive eigenvalue branch, as $\bar{I}\to0$ (or equivalently, $\bar{t}\to0$), all terms in the numerators and denominators of \eqref{eq:exactrelerrorIbar} tend towards a finite value. Consequently, the limit yields the relative error.
\begin{equation}\label{eq:exactrelerrorIbar2}
	\lim_{h\to0}\frac{\lambda^h_{n+}}{\lambda_{n+}}-1
		=\frac{\beta_0^h}{\beta_0}\frac{1}{2}\left[1+\sqrt{1-\frac{\tilde{a}_0^h}{(\beta_0^h)^2}}\,\right]-1\,.
\end{equation}
Next, in the limit of refinement ($h\to0$) we obtain $\tilde{a}_0^h/(\beta_0^h)^2\to0$, and $\beta_0^h/\beta_0\to\frac{1}{(\xi\pi)^2}\frac{\tilde{\mat{K}}_m^{11}}{\tilde{\mat{M}}}$. Interestingly, this means we recover \eqref{eq:aymptoticmembranedominated} in the limit of $\bar{I}\to0$, which was an expression associated to the negative eigenvalue branch. Thus, when taking limits as $h\to0$ and $\bar{I}\to0$, the order of the limits matter.

The more interesting case arises when considering the negative eigenvalue branch and examining the limit as $\bar{I} \to 0$. In this case, the denominator involving the square root approaches zero, while the numerator converges to the finite value of $\frac{\tilde{a}_0^h}{2(\beta_0^h)^2}$. Consequently, the relative error becomes unbounded! Indeed, a proper Laurent series expansion reveals a leading term of $(\frac{\beta_0}{\beta_0^h})(\frac{\tilde{a}_0^h}{\tilde{a}_1})(\frac{1}{\bar{I}})$. This indicates that for small values of $\bar{I}$, the error can become exceedingly large. However, it is important to note that this effect diminishes as the mesh is refined since it can be observed that $\tilde{a}_0^h/\tilde{a}_1$ rapidly approaches zero as $h \to 0$. In other words, the ``blowup'' resulting from small $\bar{I}$ or $\bar{t}$ disappears as the mesh is refined. Nevertheless, it is highly likely that this phenomenon is responsible for the significant relative error observed in our results for moderate values of $h$ (as seen in the negative eigenvalue branch of the standard method in Figure \ref{fig:normalized_spectra_error}). We will discuss this further in the next section.
\section{Discussion \label{sec:7}}
In this section we discuss the qualitative and quantitative mathematical results obtained in the previous sections. Specifically, we explore the dependence of locking behavior on the physical parameters and provide insightful discussions on this topic.

\subsection{Characterization of membrane locking}
To initiate our discussion, it is worth noting that membrane locking has traditionally been described in \textit{physical} rather than \emph{mathematical} terms. However, this does preclude its mathematical description, as we aim to demonstrate in the subsequent analysis, even if it deviates from other conventional forms of locking phenomena. We reiterate the key advantages of our one-dimensional Euler-Bernoulli circular beam model: the deliberate absence of shear strains, a constant nonzero curvature, and periodic boundary conditions that eliminate other spurious behaviors due to boundaries. These unique features positions our model as highly suitable for isolating and characterizing membrane locking.

To dispense with rigid-body transformations, consider $V$ to be the span of eigenfunctions not associated with rigid-body motions (i.e., with nonzero eigenvalues). In technical terms, considering a linear functional $l: V \rightarrow \mathbb{R}$ with a given forcing within the range of the operator, the system governed by the standard variational formulation can be represented as:
\begin{align}
	\aa{\vect{u}}{\vect{\delta u}} = l(\vect{\delta u}) \qquad \forall \vect{\delta u} = (\delta u, \delta w) \in V,
\end{align}
where the bilinear operator $a$ is defined as in equation \eqref{eq:bilinearform}. This system has a solution (up to rigid-body translations) expressed by:
\begin{align}\label{eq:solutioneigenexpansion}
	\vect{u} = \sum_{n=1}^\infty \alpha_n \, \vect{u}_n 
	\qquad \text{with } \alpha_n = l(\vect{u}_n) / \lambda_n.
\end{align}
This solution involves an expansion in terms of the eigenmodes $\vect{u}_n$, with corresponding eigenvalues $\lambda_n\neq0$. The coefficients $\alpha_n$ are determined based on the eigenvalues and the linear functional $l$. This expression also holds true in the discrete case, but the sum is finite up to the total number of modes being discretized.

The total error in approximating $\vect{u}$ in \eqref{eq:solutioneigenexpansion} is a combination of the error accrued in approximating the eigenmodes $\vect{u}_n$ and its corresponding eigenvalues $\lambda_n$. The approximation error in the modes is influenced by the accuracy of each component, namely the discrete and analytical Fourier modes, as well as the approximation error in their relative amplitudes. It is possible to demonstrate that the discrete Fourier modes provide optimal approximations of the analytical modes. Additionally, in Figure \ref{fig:normalized_amplitude_error} we show that the relative amplitudes associated with the membrane- and bending-dominated responses are well approximated, even when using coarse meshes. Based on these observations, we argue that the primary source of error stems from the misapproximation of the eigenvalues.

In Section \ref{sec:eigenvalueerrors} (Figures \ref{fig:normalized_spectra_error} and \ref{fig:spectra_negative_p34_t0015}), we have observed that the discrete eigenvalues corresponding to the negative branch significantly overestimate the true eigenvalues, particularly for coarse discretizations. On the other hand, all other eigenvalues, especially those in the positive branch, are well approximated. This discrepancy in the negative eigenvalue branch, particularly for low mode numbers $n$, results in a significant underestimation of the associated coefficients $\alpha_n$ in \eqref{eq:solutioneigenexpansion} since they are inversely proportional to the eigenvalues ($\alpha_n = l(\vect{u}_n) / \lambda_n$). Figure \ref{fig:analytical_eigenmodes} further confirms that these eigenvalues correspond to a bending-dominated response, indicated by $U_{n-} / W_{n-} < 1$. Consequently, the bending-dominated response is unduly suppressed in the discrete solution, leading to a relative over-representation of the membrane-dominated contributions. In essence, the discrete solution exhibits membrane-locking. In light of the fact that the bending part is underestimated while the membrane response is accurately represented, one could argue if ``membrane locking'' is ``the right'' term for the phenomenon being observed.

This discussion serves to reinforce and validate the arguments presented throughout this article, further justifying our focus on the spectrum and the proposed criterion developed in Section \ref{sec:lockingcriterion}. Furthermore, with an explicit expression for the relative error in the eigenvalues associated with membrane locking, we are well-equipped to further study this phenomenon. In the following subsections, we will explore the implications of the observed eigenvalue errors and investigate their causes.

\subsection{Membrane locking as a function of mesh size and physical parameters}

Upon examining expressions \eqref{eq:MnMnhIbar}--\eqref{eq:exactrelerrorIbar}, we can gain insights into the phenomenon of membrane locking as it relates to the mesh size $h$ and the physical parameters, which, as we previously discussed, are captured by the normalized thickness $\bar{t}=\frac{t}{R}$ or equivalently the normalized moment of inertia $\bar{I}=\frac{\bar{t}^2}{12}$.

The occurrence of membrane locking can be attributed to a fundamental difference among the expressions for the terms $M_n$ and $M_n^h$ in \eqref{eq:MnMnhIbar} and \eqref{eq:acoeffsIbar}: $M_n$ is a multiple of $\bar{I}=\frac{\bar{t}^2}{12}$ so, as a polynomial in $\bar{I}$, it lacks constant and quadratic terms; however, in the discrete case, $M_n^h$ is a polynomial that does include constant and quadratic terms, specifically denoted by $\tilde{a}_0^h$ and $\tilde{a}_2^h$ in \eqref{eq:acoeffsIbar}. The nonzero values of $\tilde{a}_0^h$ and $\tilde{a}_2^h$ have a profound impact on the behavior of the discrete solution.

\begin{figure}[!h]
	\centering
	\includegraphics[trim = 0cm 0cm 0cm 0cm, clip,width=0.48\textwidth]{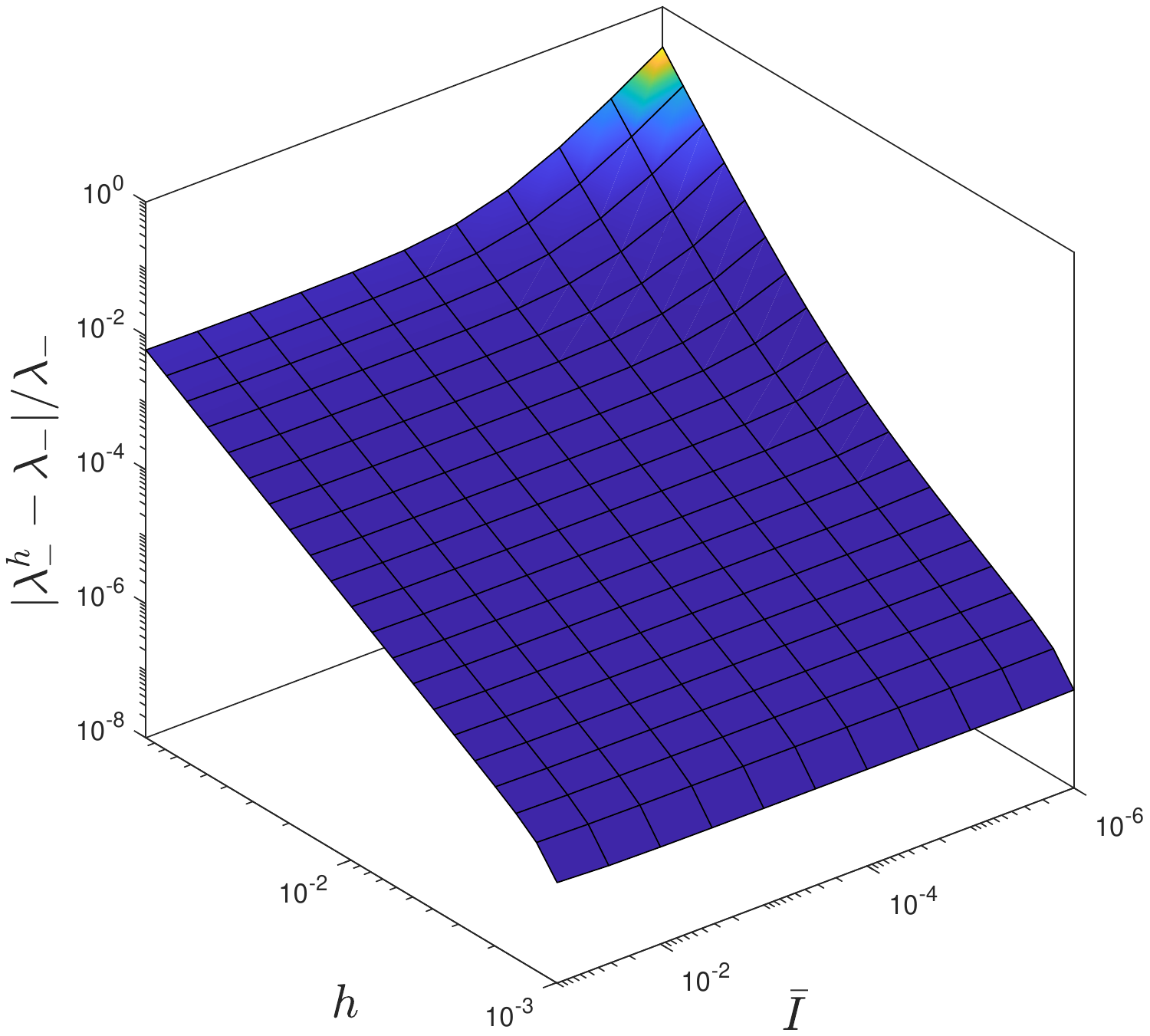}
	\caption{Normalized eigenvalue error of the second eigenvalue, $n=2$, of the negative branch as a function of mesh size $h = \pi / N$ and the normalized moment of inertia $\bar{I} = \bar{t}^2 / 12$. Note the blowup that occurs on coarse meshes when $\bar{I} \to 0$. This is a manifestation of membrane locking.}
	\label{fig:behavior_as_a_fun_of_h_and_I}
\end{figure}

\begin{figure}[!h]
	\centering
	\includegraphics[trim = 0cm 0cm 0cm 0cm, clip,width=0.7\textwidth]{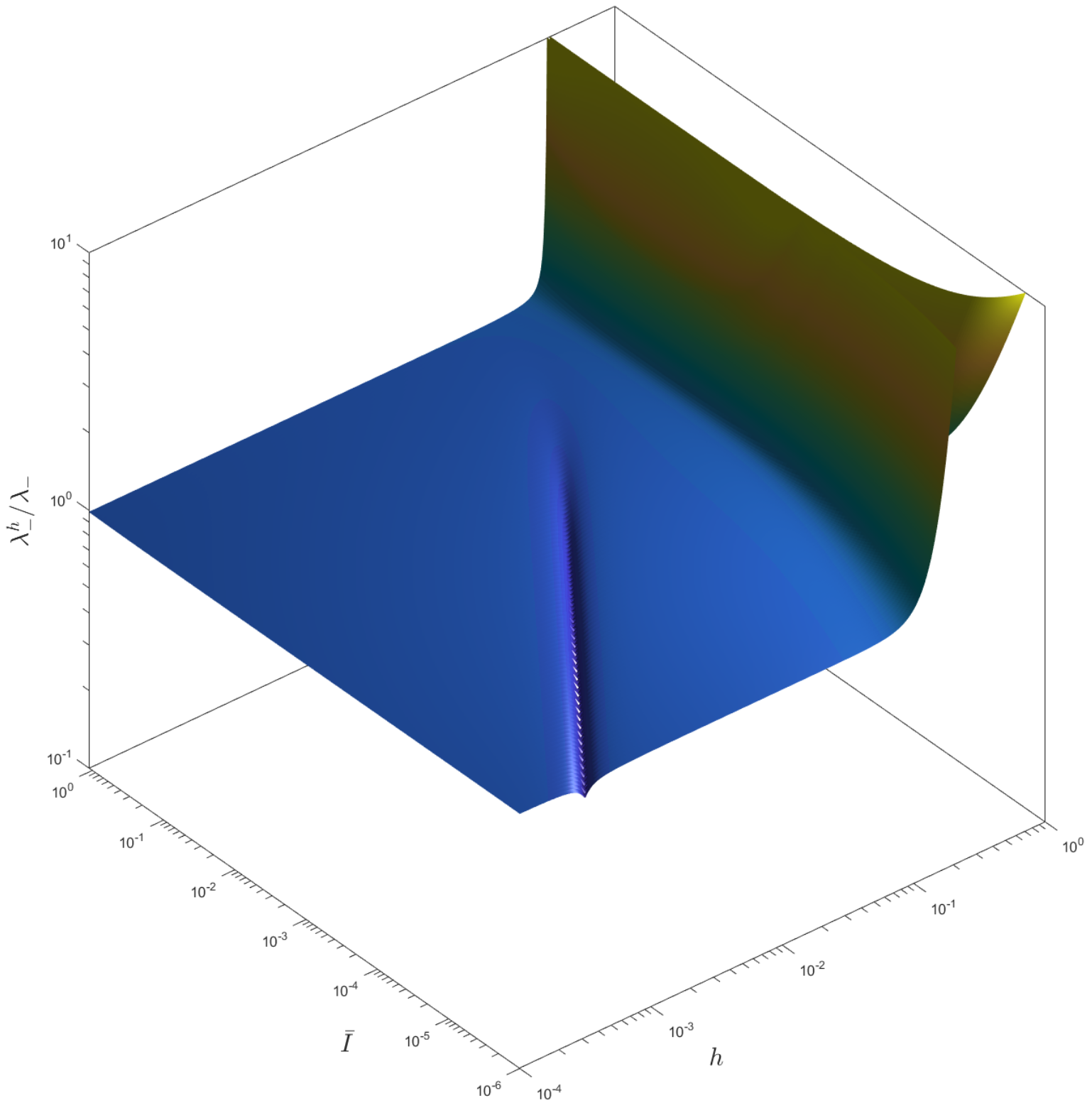}
	\caption{Landscape of the normalized eigenvalues at fixed $\xi = 0.1$ as a function of (continuous) mesh size $h$ and the normalized moment of inertia $\bar{I} = \bar{t}^2 / 12$. Note the blowup that occurs on coarse meshes, which gets more severe when $\bar{I} \to 0$. This is a manifestation of membrane locking. Note the peculiar behavior that runs diagonally across the plot, which occurs approximately for $h\approx\xi\pi\sqrt{\bar{I}}$, and corresponds to when this negative eigenvalue branch becomes membrane-dominated (recall \eqref{eq:nhatformula} which translates to $\hat{n}(\bar{I})\approx1/\sqrt{\bar{I}}$ and Section \ref{sec:discretebranchswitching}).}
	\label{fig:behavior_as_a_fun_of_h_and_I_fixed_xi}
\end{figure}

First, the presence of $\tilde{a}_0^h$ leads to a singularity in the relative error as $\bar{I}$ approaches zero with the leading relevant coefficient being $\tilde{a}_0^h/\tilde{a}_1$. As the mesh is refined, this effect quickly subsides because this coefficient converges rapidly to zero, specifically at a rate of $\mathcal{O}(h^4)$. This is shown in Figure \ref{fig:behavior_as_a_fun_of_h_and_I} where the $n=2$ eigenvalue of the negative branch is seen to be very large at large values of $h$ and small values of $\bar{I}$. At a fixed value of $\xi$, Figure \ref{fig:behavior_as_a_fun_of_h_and_I_fixed_xi} also shows how the negative branch of eigenvalues is severely over-approximated at coarse $h$, but quickly recovers accuracy for small $h$. 

On the other hand, the non-zero coefficient $\tilde{a}_2^h$ introduces a term that is linear in $\bar{I}$ (quadratic in $\bar{t}$). The relevant coefficient for this term is $\frac{\tilde{a}_2^h}{\tilde{a}_1}$. Unlike $\tilde{a}_0^h$, this coefficient does not rapidly diminish as the mesh is refined for a fixed value of $\xi$. However, it may decrease with $p$-refinement, as discussed further below. In fact, according to \eqref{eq:aymptoticmembranedominated}, this term eventually becomes part of the \emph{membrane}-dominated response in the limit of mesh refinement (and this is indeed what happens to the left of the diagonal dip observed in Figure \ref{fig:behavior_as_a_fun_of_h_and_I_fixed_xi}).

Thus, the non-trivial interplay between $\tilde{a}_0^h$ and $\tilde{a}_2^h$ is crucial in understanding the occurrence of membrane locking. When dealing with coarse meshes, the magnitude of $\tilde{a}_0^h$ is significant, resulting in the $\frac{1}{\bar{I}}$ term to dominate the expression of relative error. As a consequence, the relative error becomes large as $\bar{I}$ decreases, leading to membrane locking. Conversely, in the case of fine meshes with a fixed $\xi$ value, $\tilde{a}_0^h$ approaches zero, allowing the linear term in $\bar{I}$ to dominate. Consequently, the relative error becomes large as $\bar{I}$ increases.

Finally, it is worth noting that there are mathematical similarities between membrane locking and shear locking. Specifically, membrane locking occurs when $\bar{I}$ is small, corresponding to a small thickness $t$. This resemblance is not unexpected, as shear locking also introduces an artificial bending stiffness, particularly when the thickness $t$ approaches zero (\cite[Chapter 9.4, Page 368]{wriggers2008nonlinear}.).

\subsection{Membrane locking as a function of polynomial degree}
As pointed out in Section \ref{sec:6} the dependence on polynomial degree $p$ lies hidden in the circulant eigenvalues $\tilde{\mat{K}}_b^{ij}$, $\tilde{\mat{K}}_m^{ij}$ and $\tilde{\mat{M}}$. This dependence implies that the asymptotic limits of mesh refinement are influenced by the choice of $p$. Understanding these limits is crucial, as they not only determine the convergence rate but also affect the absolute value of the relative error.

While increasing the value of $p$ significantly improves the convergence rate and the overall accuracy of the method, it does not eliminate the parasitic behavior we identify as membrane locking as defined in Section \ref{sec:lockingcriterion}. This is because the relative error of the eigenvalues still deviates significantly from the asymptotic curve, particularly at low values of $\xi$. However, it is worth mentioning that the terms $\tilde{\mat{K}}_m^{11}\tilde{\mat{K}}_m^{22}-|\tilde{\mat{K}}_m^{12}|^2$ and $\tilde{\mat{K}}_b^{11}\tilde{\mat{K}}_b^{22}-|\tilde{\mat{K}}_b^{12}|^2$ appearing in $\tilde{a}_0^h$ and $\tilde{a}_2^h$ from   \eqref{eq:acoeffsIbar}, tend to decrease as the value of $p$ increases. This reduction may potentially help in mitigating the effects of membrane locking. However, it is important to note that in the case of cubic and quartic discretizations, we have not observed a significant reduction in membrane locking despite the presence of these terms. Membrane locking remains prominent in these cases, as may be observed in Figures \ref{fig:spectra_negative_p34_t0015} and \ref{fig:locking_p34}.

Lastly, as discussed in Section \ref{sec:standarderror}, the preasymptotic domain can be effectively described by the relative error in the asymptotic limit of refinement, which is $p$-dependent. For instance, Figure \ref{fig:spectra_negative_p34_t0015} illustrates that for $p=4$, the onset of ``good" convergence initiates at an earlier stage (for coarser values of $h$ at a fixed $n$). This occurs because the accelerated decay of the error occurs at a larger value of $\xi$ compared to $p=2$ and $p=3$.

\subsection{Comparison of the standard and mixed formulation}
The mixed Galerkin discretization based on the Hellinger-Reissner principle was shown to be largely free of membrane locking (Figure \ref{fig:locking_p2}). Furthermore, the asymptotic accuracy of eigenvalues for bending-dominated modes was significantly better in the mixed formulation compared to the standard formulation, see Figure \ref{fig:normalized_spectra_error} and \ref{fig:spectra_negative_p34_t0015}. A more detailed study of the relative error in the case of this variational formulation is viable, but is left for a future analysis.

\section{Conclusion \label{sec:8}}
This paper addresses the issue of membrane locking which affects finite element approximations of thin beams and shells. Membrane locking, along with transverse shear locking and volumetric locking, negatively impacts the accuracy and reliability of finite element discretizations. However, membrane locking is still not well understood, and the absence of a rigorous assessment methodology presents challenges in the development of effective beam and shell elements. Our research aims to achieve two main goals: firstly, to gain insights into the mechanisms underlying membrane locking, and secondly, to establish a rigorous criterion for accurately identifying and evaluating the extent of membrane locking.

Our study involved a comprehensive Fourier analysis of the standard and mixed Galerkin formulations of the Euler-Bernoulli beam model for a circular ring, employing smooth splines for the discretization. Importantly, this model selection enabled the effective isolation and investigation of membrane locking due to several key factors: the model inherently lacks shear strains, maintains a constant nonzero curvature, and incorporates periodic boundary conditions to eliminate other undesired phenomena, including e.g. boundary layers. Furthermore, our choice of model facilitated exact formulas for the \emph{discrete} eigenvalues and corresponding eigenfunctions by exploiting properties of circulant matrices. The exact expressions were presented in terms of the relevant discretization and physical parameters, and compared to the exact analytical solutions.

Our analysis revealed two distinct types of eigenmodes: bending-dominated and membrane-dominated. Remarkably, we observed that all eigenmodes were well approximated even with coarse to medium discretizations. However, the eigenvalues associated with bending modes in the standard formulation exhibited significant deviations, especially in coarse discretizations. This motivated the development of a criterion to accurately identify and quantify the occurrence of membrane locking based on accuracy of eigenvalues. This criterion effectively illuminated the presence and severity of membrane locking in the standard formulation for quadratics, cubics and quartics. Despite the substantial improvement in absolute accuracy with increasing polynomial order, we found that it did not mitigate the occurrence of locking behavior. In contrast, when utilizing a mixed Galerkin discretization to numerically solve the same problem, we observed minimal effects of membrane locking on the solutions across all polynomial degrees and meshes. Moreover, we discovered that the absolute accuracy of the mixed method was significantly superior, particularly in relation to the eigenvalues associated with bending-dominated modes.

Furthermore, we mathematically derived and analyzed exact expressions for the relative error. We established the expected convergence rates for a fixed mode number and determined that the normalized spectra in the asymptotic refinement limit exhibited a relationship with the spectra associated with isolated bending and membrane responses. Significantly, our investigation revealed that in coarse discretizations of the standard formulation, the eigenvalues associated with the bending-dominated response were grossly overestimated, particularly for a small thickness and large radius of curvature. The result is a suppression of the bending-dominated component of the discrete solution. Hence, our careful eigenvalue analysis demonstrates that the phenomenon observed is not simply an overestimation of the membrane response but rather a suppression of bending, prompting a question regarding the suitability of the term ``membrane locking" for this phenomenon. We hypothesize that this effect shares mathematical similarities with shear locking, despite our model not incorporating shear strains. 

Looking ahead, we plan to investigate membrane locking in doubly curved shells and explore the extension of our analytical techniques to periodic, uniform spline discretizations of a torus. Additionally, we would like to explore other types of shell elements, such as non-conforming techniques based on $C^0$ finite elements or DG. Despite the resulting submatrices not being circulant, they exhibit a circulant block structure with block size proportional to the polynomial order, which may enable analytical solution techniques. Furthermore, we note that our techniques may provide insights into other forms of locking, specifically shear locking in thin beams.

\section*{Acknowledgments}
R.R. Hiemstra and D. Schillinger gratefully acknowledge financial support from the German Research Foundation (Deutsche Forschungsgemeinschaft) through the DFG Emmy Noether Grant SCH 1249/2-1 and the standard DFG grant SCH 1249/5-1. F. Fuentes' work was partially supported by the National Center for Artificial Intelligence CENIA FB210017, Basal ANID. The authors thank T.J.R Hughes and T.-H. Nguyen for helpful discussions.

\appendix
\section{Eigenvalue decomposition of matrices with circulant blocks \label{sec:EVcirculant}}
The discrete eigenvalue problems that arise in this paper involve two-by-two block matrices with circulant blocks, sometimes referred to as circulant block matrices (as opposed to block-circulant matrices). First, we describe the special structure of circulant matrices and their eigenvalue decomposition based on the Discrete Fourier Transform (DFT). Then we use these properties to diagonalize two-by-two block matrices, which are the focus of this paper.

\subsection{Eigenvalue decomposition of a circulant matrix \label{sec:circulent-matrix}}
A circulant $\N \times \N$ real matrix is uniquely defined by $\N$ real numbers, $c_0, c_1, \ldots , c_{\N-1}$, and has the following special structure
\begin{align}
	\mat{C} = 
	\begin{bmatrix} 
		c_0 			& c_{1} 		& \cdots 		& c_{\N-2} 	& c_{\N-1}						\\
		c_{\N-1} 	& c_0 			& c_{1} 		& \cdots 		& c_{\N-2}						\\
		\vdots 		&	\ddots		&	\ddots		&	\ddots		&	\vdots							\\
		c_2				&	\cdots		&	c_{\N-1}	& c_0 			& c_1								\\
		c_1				& c_2				& \cdots 		& c_{\N-1}	& c_0
	\end{bmatrix}
\end{align}
Every circulant matrix has an eigenvalue decomposition 
\begin{align}
	\mat{C} = \mat{F} \mat{\Lambda} \mat{F}^{-1}
\end{align}
where $\mat{F} \in \mathbb{C}^{N\times N}$, called the Fourier matrix, denotes the collection of eigenvectors and $\mat{\Lambda}$ the corresponding eigenvalues. The eigenvectors, the columns of the Fourier matrix, are always the same, namely the Fourier modes,
\begin{align}
	& \mat{v}_j = 
	\frac{1}{\sqrt{\N}} 
	\begin{pmatrix}
		1 & \omega^j & \omega^{2j} & \ldots & \omega^{(\N-1)j}
	\end{pmatrix}^T, \; \; j=0,1,\ldots, \N-1,& 
	&\text{with}& 
	&\omega = \exp{(2\pi \cdot i / \N)}&
\end{align}
The corresponding eigenvalues are defined via the following simple relationship
\begin{align}\label{eq:EVcirculant}
	\lambda_j = c_0 + c_{1} \cdot \omega^j + c_{2} \cdot \omega^{2j} +  \cdots  + c_{\N-1} \cdot \omega^{(\N-1)j}, \qquad j=0, 1 , \ldots , \N-1.
\end{align}
If $\mat{C}$ is symmetric then the eigenvalues are all real and are given by
\begin{align}\label{eq:EVcirculantsymmetric}
	\lambda_j = c_0 + 2c_{1} \cdot \cos(2\pi \cdot j /\N) + 2c_{2} \cdot \cos(4\pi \cdot j /\N) + \ldots + 2c_{\lfloor \N/2 \rfloor} \cdot \cos(\lfloor \N/2 \rfloor \pi \cdot j / \N).
\end{align}
On the other hand, if $\mat{C}$ is skew-symmetric then the eigenvalues are complex and are given by
\begin{align}\label{eq:EVcirculantantisymmetric}
	\lambda_j = 2i c_{1} \cdot \sin(2\pi \cdot j /\N) + 2i c_{2} \cdot \sin(4\pi \cdot j /\N) + \ldots + 2i c_{\lfloor \N/2 \rfloor} \cdot \sin(\lfloor \N/2 \rfloor \pi \cdot j / \N).
\end{align}

\subsection{Similarity transformations of matrices with circulant blocks \label{app:similarity_transforms}}
Let $\mat{A}, \mat{B}, \mat{C}, \mat{D} \in \mathbb{R}^{\N \times \N}$ be circulant matrices. We are interested in the eigenvalue decomposition of the $2 \times 2$ block matrix
\begin{align}
	\mat{K} = \begin{bmatrix} \mat{A} & \mat{B}  \\ \mat{C} & \mat{D}  \end{bmatrix}
\end{align}
Using the Fourier matrix $\mat{F}$ we can bring matrix $\mat{K}$ to the following simpler form
\begin{align}
\tilde{\mat{K}} = 
\begin{bmatrix} 
	\tilde{\mat{A}} 	& \tilde{\mat{B}}  \\ 
	\tilde{\mat{C}}	& \tilde{\mat{D}}  
\end{bmatrix} =
\begin{bmatrix} 
	\mat{F}^{-1} 	&   \\ 
	 						& \mat{F}^{-1}  
\end{bmatrix}
\begin{bmatrix} 
	\mat{A} & \mat{B}  \\ 
	\mat{C} & \mat{D}  
\end{bmatrix}
\begin{bmatrix} 
	\mat{F} 	&   \\ 
	 						& \mat{F}  
\end{bmatrix}
\end{align}
where $\tilde{\mat{A}}, \tilde{\mat{B}}, \tilde{\mat{C}}$ and $\tilde{\mat{D}}$ are diagonal with entries the eigenvalues of the respective matrices. This matrix is easier to diagonalize than the original matrix $\mat{K}$.

\subsection{Eigenvalue decomposition of a $2\times2$ block matrix with diagonal blocks \label{sec:circulent-block-matrix}}
Let $\tilde{\mat{K}}$ denote a $2\times2$ block matrix with diagonal blocks. Since diagonal matrices commute, $\tilde{\mat{C}} \tilde{\mat{D}} = \tilde{\mat{D}} \tilde{\mat{C}}$, we have that the characteristic equation is similar to that of a $2 \times 2$ matrix (see \cite{silvester2000determinants}),
\begin{align}
\det{ \left(\tilde{\mat{K}} - \lambda \mat{I} \right)} &= 
\det
\left(
\begin{bmatrix} 
	\tilde{\mat{A}}-\lambda_n \mat{I} 			& \tilde{\mat{B}}  \\ 
	\tilde{\mat{C}}											& \tilde{\mat{D}}-\lambda_n \mat{I}  
\end{bmatrix}\right)\\
&= \det{\left(\left( \tilde{\mat{A}}-\lambda_n \mat{I} \right)\left( \tilde{\mat{D}}-\lambda_n \mat{I} \right) - \tilde{\mat{B}} \tilde{\mat{C}} \right)}	\\
&= \det{\left(\tilde{\mat{A}} \tilde{\mat{D}} - \tilde{\mat{B}} \tilde{\mat{C}} - \lambda_n \left(\tilde{\mat{A}} + \tilde{\mat{D}} \right) + \lambda_n^2 \mat{I}  \right)}	\\
&=\prod_{n=1}^\N \left(\tilde{A}_{n} \tilde{D}_{n} - \tilde{B}_{n} \tilde{C}_{n} - \lambda_n \left(\tilde{A}_{n} + \tilde{D}_{n} \right) + \lambda_n^2  \right)
\end{align}
The characteristic equation of the $2N \times 2N$ system is simply the product of $N$ characteristic equations of $2\times2$ matrices. Hence, instead of the $2N \times 2N$ matrix eigenvalue problem, we may consider, for each $n$, the $2\times 2$ matrix eigenvalue problem 
\begin{align}
\begin{bmatrix} 
	\tilde{A}_n - \lambda_n 	& \tilde{B}_n  \\ 
	\tilde{C}_n	& \tilde{D}_n - \lambda_n  
\end{bmatrix}
\begin{bmatrix} 
	U_n  \\ 
	W_n  
\end{bmatrix}
= 
\begin{bmatrix} 
	0  \\ 
	0  
\end{bmatrix}.
\label{eq:two_by_two_problem_abstract}
\end{align}
The eigenvalues follow as the roots of the characteristic equation
\begin{subequations}
\begin{align}
	\lambda_{n\pm} 	&= K_n \left(1 \pm L_n\right)
\end{align}
where $K_n$ and $L_n$ are constants that depend on $n$
\begin{align}	
	K_n 	= \frac{\tilde{A}_{n} + \tilde{D}_{n}}{2}, \qquad
	M_n 	= \tilde{A}_{n} \cdot \tilde{D}_{n} - \tilde{C}_{n} \cdot \tilde{B}_{n}, \qquad
	L_n 	=\sqrt{1 - M_n / K^2_n}.	
\end{align}
\end{subequations}
The eigenfunctions are defined up to a constant. However, the ratio of the amplitudes the eigenfunctions, $U_n / W_n$, satisfy
\begin{align}
	\rho_{n\pm} = \frac{U_{n\pm}}{W_{n\pm}} = \frac{\tilde{B}_n}{\lambda_{n\pm} - \tilde{A}_n}
\end{align}
\section{Matrix eigenvalue problem \label{app:splines}}
We compute exact expressions for the entries of the block-circulant stiffness and mass matrices that arise in the standard and mixed formulation of the discrete eigenvalue problem based on approximation with uniform spline basis functions of degree $p=2,3,4$. The expressions are given in terms of rational numbers and symbolic variables, such as the mesh size $h$, the radius $R$ and the normalized beam thickness $\bar{t}=\tfrac{t}{R}$.

\subsection{Standard formulation  \label{app:standard}}
The matrix eigenvalue problem is
\begin{align}
	\mmat{K} \, \mmat{u} = \lambda^h \mmat{M} \, \mmat{u}.
\end{align}
The mass and stiffness matrices are $2 \times 2$ block-circulant matrices that have the structure
\begin{align}
	&\mathbf{K}  = \mathbf{K}_m + \mathbf{K}_b
			 = 	\begin{bmatrix} 
						\mat{K}^{11}_m 		& \mat{K}^{12}_m 			\\ 
						\mat{K}^{21}_m 		& \mat{K}^{22}_m 
					\end{bmatrix} + 
					\begin{bmatrix} 
						\mat{K}^{11}_b 		& \mat{K}^{12}_b 			\\ 
						\mat{K}^{21}_b 		& \mat{K}^{22}_b 
					\end{bmatrix},&
	& \mathbf{M} = 
	\begin{bmatrix}
		\mat{M} 	&  				\\  
						& \mat{M}
	\end{bmatrix}.&
\end{align}
Here $ \mathbf{K}_m$ denote the contribution to the membrane stiffness and  $\mathbf{K}_b$ the portion related to the bending. The blocks $\mat{M}$ are circulant, symmetric positive definite matrices, and are computed as
\begin{align}\label{eq:MassMtandardform}
	& [\mat{M}]_{ij} = R \int_{0}^{2 \pi} B_{i,p}(\theta) B_{j,p}(\theta) \, \mathrm{d} \theta.
\end{align}
In the stiffness matrix we make a distinction between the action of membrane and bending terms, respectively. The submatrices on the diagonal are circulant and symmetric. The off-diagonal submatrices are skew-symmetric circulant matrices. These are computed according to the relations   
\begin{equation}\label{eq:KbKmstandardform}
\begin{aligned}
	[\mat{K}^{11}_{m}]_{ij} 	&= \frac{1}{R} \int_{0}^{2\pi} B^\prime_{i,p}(\theta) B^\prime_{j,p}(\theta) \, \mathrm{d} \theta,						\qquad\quad &
	[\mat{K}^{11}_{b}]_{ij} 	&= \frac{1}{R}\frac{\bar{t}^2}{12} \int_{0}^{2\pi} B^\prime_{i,p}(\theta) B^\prime_{j,p}(\theta) \, \mathrm{d}\theta,	\\
	[\mat{K}^{12}_{m}]_{ij}	&= \frac{1}{R}\int_{0}^{2\pi} B^\prime_{i,p}(\theta) B_{j,p}(\theta) \, \mathrm{d}\theta, 										\qquad\quad &
	[\mat{K}^{12}_{b}]_{ij}	&= \frac{1}{R}\frac{\bar{t}^2}{12} \int_{0}^{2\pi} -B^\prime_{i,p}(\theta) B^{\prime \prime}_{j,p}(\theta) \, \mathrm{d}\theta 	\\
	[\mat{K}^{22}_{m}]_{ij}	&= \frac{1}{R} \int_{0}^{2\pi} B_{i,p}(\theta) B_{j,p}(\theta) \, \mathrm{d}\theta, 														\qquad\quad &
	[\mat{K}^{22}_{b}]_{ij}	&= \frac{1}{R}\frac{\bar{t}^2}{12} \int_{0}^{2\pi} B^{\prime \prime}_{i,p}(\theta) B^{\prime \prime}_{j,p}(\theta) \, \mathrm{d}\theta
\end{aligned}
\end{equation}
Table \ref{tab:circulant_matrices} lists the non-zero values of the circulant blocks of the mass and stiffness matrix in the case of $\p=2, 3$ and $4$. Inverses, transposes and products of circulant matrices are circulant. Hence, the blocks in the matrix $\tilde{\mmat{K}} = \mmat{M}^{-1} \mmat{K}$ remain circulant and the technique in \ref{sec:EVcirculant} can be used to compute the discrete eigenvalues and eigenvectors of the system.

\begin{table}[h!]
\begin{center}
\caption{Circulant submatrices of the mass and stiffness matrix of the standard formulation}
\label{tab:circulant_matrices}
\begin{tabular}{C{8em} C{4.8em} C{4.8em} C{4.8em} C{4.8em} C{4.8em}}
		\toprule
		$\frac{1}{R}\frac{1}{h}\mat{M}=R\frac{1}{h}\mat{K}^{22}_{m}$	& $c_0$						& $c_1 = c_{N-1}$			& $c_2 = c_{N-2}$				& $c_3 = c_{N-3}$			& $c_4 = c_{N-4}$			\\
		\midrule
		$\p=2$	& $\frac{11}{20}$				& $\frac{13}{60}$ 				& $\frac{1}{120}$ 			& 									& 								\\ 
		$\p=3$ 	& $\frac{151}{315}$			& $\frac{397}{1680}$ 		& $\frac{1}{42}$ 			& $\frac{1}{5040}$		& 								\\
		$\p=4$ 	& $\frac{15619}{36288}$	& $\frac{44117}{181440}$	& $\frac{913}{22680}$ 	& $\frac{251}{181440}$	& $\frac{1}{362880}$\\
		\midrule
		 $Rh\mat{K}^{11}_{m}=Rh\frac{12}{\bar{t}^2}\mat{K}^{11}_{b}$	& $c_0$						& $c_1 = c_{N-1}$			& $c_2 = c_{N-2}$				& $c_3 = c_{N-3}$			& $c_4 = c_{N-4}$			\\
		\midrule
		$\p=2$ 	& $1$								& $-\frac{1}{3}$ 				& $-\frac{1}{6}$ 			& 								& 								\\
		$\p=3$ 	& $\frac{2}{3}$					& $-\frac{1}{8}$ 				& $-\frac{1}{5}$			& $-\frac{1}{120}$	&								\\
		$\p=4$ 	& $\frac{35}{72}$ 			& $-\frac{11}{360}$ 			& $-\frac{17}{90}$ 		& $-\frac{59}{2520}$ & $-\frac{1}{5040}$	\\
		\midrule
		 $R\mat{K}^{12}_{m}$	& $c_0$						& $c_1 = -c_{N-1}$		& $c_2 = -c_{N-2}$		& $c_3 = -c_{N-3}$		& $c_4 = -c_{N-4}$			\\
		\midrule
		$\p=2$ 	& $0$ 	& $\frac{5}{12}$ 			& $\frac{1}{24}$			& 									&								\\
		$\p=3$ 	& $0$	& $\frac{49}{144}$ 		& $\frac{7}{90}$			& $\frac{1}{720}$			&								\\
		$\p=4$ 	& $0$	& $\frac{809}{2880}$ 	& $\frac{289}{2880}$ 	& $\frac{41}{6720}$ 		& $\frac{1}{40320}$ 	\\
		\midrule
		$Rh^2\frac{12}{\bar{t}^2}\mat{K}^{12}_{b}$	& $c_0$							& $c_1 = -c_{N-1}$		& $c_2 = -c_{N-2}$	& $c_3 =-c_{N-3}$		& $c_4 = -c_{N-4}$			\\
		\midrule
		$\p=2$		& 	$0$ 		& $1$ 						& $-\frac{1}{2}$				&									&										\\	
		$\p=3$		&	$0$ 		& $\frac{19}{24}$ 		& $-\frac{1}{3}$ 				& $-\frac{1}{24}$			& 										\\
		$\p=4$		&	$0$ 		& $\frac{217}{360}$ 	& $-\frac{67}{360}$ 			& $-\frac{3}{40}$ 			& $-\frac{1}{720}$			\\
		\midrule
		$Rh^3\frac{12}{\bar{t}^2}\mat{K}^{22}_{b}$ & $c_0$		& $c_1 = c_{N-1}$	& $c_2 = c_{N-2}$			& $c_3 = c_{N-3}$		& $c_4 = c_{N-4}$			\\
		\midrule
		$\p=2$		& 	$6$ 				& $-4$ 			& $1$				&					&					\\	
		$\p=3$		&	$\frac{8}{3}$ 		& $-\frac{3}{2}$ 		& $0$ 			& $\frac{1}{6}$			& 					\\
		$\p=4$		&	$\frac{19}{12}$ 		& $-\frac{43}{60}$ 	& $-\frac{4}{15}$ 		& $\frac{11}{60}$ 			& $\frac{1}{120}$			\\
		\bottomrule
\end{tabular}
\end{center}
\end{table}

\subsection{Mixed formulation \label{app:mixed}}
The generalized matrix eigenvalue problem is
\begin{align}
	\begin{bmatrix} 
		-\mathbf{K}_{\mmat{\varepsilon}\mmat{\varepsilon}} 		& \mathbf{K}_{\mmat{\varepsilon}\mathbf{u}} 			\\ 
		\mathbf{K}_{\mmat{\varepsilon}\mathbf{u}}^T 		& \mathbf{0} 
	\end{bmatrix}
	\begin{bmatrix} 
		\mmat{\varepsilon} \\ 
		\mathbf{u} 
	\end{bmatrix}
 = 
	\begin{bmatrix} 
		\mathbf{0} 		& \mathbf{0} 			\\ 
		\mathbf{0} 		& \lambda^h \mathbf{M} 
	\end{bmatrix}
	\begin{bmatrix} 
		\mmat{\varepsilon} \\ 
		\mathbf{u} 
	\end{bmatrix}.
\end{align}
The matrices have the following structure
\begin{align}
\mathbf{K}_{\mmat{\varepsilon}\mmat{\varepsilon}} &= 					
	\begin{bmatrix} 
		\mat{K}^{11}_{\mmat{\varepsilon}\mmat{\varepsilon}}  		& 0 			\\ 
		0		& \mat{K}^{22}_{\mmat{\varepsilon}\mmat{\varepsilon}} 
	\end{bmatrix}, &
\mathbf{K}_{\mmat{\varepsilon}\mathbf{u}} &= 					
	\begin{bmatrix} 
		\mat{K}^{11}_{\mmat{\varepsilon}\mathbf{u}}  	& \mat{K}^{12}_{\mmat{\varepsilon}\mathbf{u}} 			\\ 
		\mat{K}^{21}_{\mmat{\varepsilon}\mathbf{u}} 		& \mat{K}^{22}_{\mmat{\varepsilon}\mathbf{u}} 
	\end{bmatrix},
	& \mathbf{M} = 
	\begin{bmatrix}
		\mat{M} 	&  	0			\\  
				0		& \mat{M}
	\end{bmatrix}.&
\end{align}
The non-zero blocks of the mass and stiffness matrix are all circulant. The entries of $\mat{M}$ are computed as
\begin{align}
	& [\mat{M}]_{ij} = R \int_{0}^{2 \pi} B_{i,p}(\theta) B_{j,p}(\theta) \, \mathrm{d} \theta,
\end{align}
and the entries of the submatrices of $\mathbf{K}_{\mmat{\varepsilon}\mmat{\varepsilon}}$ and $\mathbf{K}_{\mmat{\varepsilon}\mathbf{u}}$ are given by
\begin{align}
[\mat{K}^{11}_{\mmat{\varepsilon}\mmat{\varepsilon}}]_{ij}  &= R \int_{0}^{2\pi} B_{i,p-1}(\theta) B_{j,p-1}(\theta) \, \mathrm{d}\theta, 	 &
[\mat{K}^{22}_{\mmat{\varepsilon}\mmat{\varepsilon}}]_{ij}  &= R^3\frac{\bar{t}^2}{12} \int_{0}^{2\pi} B_{i,p-1}(\theta) B_{j,p-1}(\theta) \, \mathrm{d}\theta \\
[\mat{K}^{11}_{\mmat{\varepsilon}\mathbf{u}}]_{ij}  &= \int_{0}^{2\pi} B_{i,p-1}(\theta) B^\prime_{j,p}(\theta) \, \mathrm{d}\theta, &
[\mat{K}^{12}_{\mmat{\varepsilon}\mathbf{u}}]_{ij}  &= \int_{0}^{2\pi} B_{i,p-1}(\theta) B_{j,p}(\theta) \, \mathrm{d}\theta	\\
[\mat{K}^{21}_{\mmat{\varepsilon}\mathbf{u}}]_{ij}  &= R\frac{\bar{t}^2}{12} \int_{0}^{2\pi} B_{i,p-1}(\theta) B^\prime_{j,p}(\theta) \, \mathrm{d}\theta, &
[\mat{K}^{22}_{\mmat{\varepsilon}\mathbf{u}}]_{ij}  &=  R\frac{\bar{t}^2}{12}  \int_{0}^{2\pi} -B_{i,p-1}(\theta) B^{\prime \prime}_{j,p}(\theta) \, \mathrm{d}\theta
\end{align}

The above generalized eigenvalue problem can be transformed to standard form using static condensation, that is, by solving for $\mathbf{\varepsilon}$ in terms of $\mathbf{u}$, leading to
\begin{align}
	&\mathbf{K} \, \mathbf{u} = \lambda^h \mathbf{M} \mathbf{u}& &\text{with} &
	&\mathbf{K} =  \mathbf{K}_{\mmat{\varepsilon}\mathbf{u}}^T \mathbf{K}^{-1}_{\mmat{\varepsilon}\mmat{\varepsilon}} \, \mathbf{K}_{\mmat{\varepsilon}\mathbf{u}}.&
\end{align}
As before, the contribution to the stiffness matrix can be factored into a contribution related to membrane strains and a contribution related to bending
\begin{align}
	&\mathbf{K}  = \mathbf{K}_m + \mathbf{K}_b
			 = 	\begin{bmatrix} 
						\mat{K}^{11}_m 		& \mat{K}^{12}_m 			\\ 
						\mat{K}^{21}_m 		& \mat{K}^{22}_m 
					\end{bmatrix} + 
					\begin{bmatrix} 
						\mat{K}^{11}_b 		& \mat{K}^{12}_b 			\\ 
						\mat{K}^{21}_b 		& \mat{K}^{22}_b 
					\end{bmatrix}.
\end{align}
where,
\begin{align}
	\mat{K}^{11}_m &= (\mat{K}^{11}_{\mmat{\varepsilon}\mathbf{u}})^T \, (\mat{K}^{11}_{\mmat{\varepsilon}\mmat{\varepsilon}})^{-1} \, \mat{K}^{11}_{\mmat{\varepsilon}\mathbf{u}}, &
	\mat{K}^{11}_b  &= (\mat{K}^{21}_{\mmat{\varepsilon}\mathbf{u}})^T \, (\mat{K}^{22}_{\mmat{\varepsilon}\mmat{\varepsilon}})^{-1} \, \mat{K}^{21}_{\mmat{\varepsilon}\mathbf{u}} \\
	\mat{K}^{12}_m &= (\mat{K}^{11}_{\mmat{\varepsilon}\mathbf{u}})^T \, (\mat{K}^{11}_{\mmat{\varepsilon}\mmat{\varepsilon}})^{-1} \, \mat{K}^{12}_{\mmat{\varepsilon}\mathbf{u}}, &
	\mat{K}^{12}_b  &= (\mat{K}^{21}_{\mmat{\varepsilon}\mathbf{u}})^T \, (\mat{K}^{22}_{\mmat{\varepsilon}\mmat{\varepsilon}})^{-1} \, \mat{K}^{22}_{\mmat{\varepsilon}\mathbf{u}} \\
	\mat{K}^{22}_m &= (\mat{K}^{12}_{\mmat{\varepsilon}\mathbf{u}})^T \, (\mat{K}^{11}_{\mmat{\varepsilon}\mmat{\varepsilon}})^{-1} \, \mat{K}^{12}_{\mmat{\varepsilon}\mathbf{u}}, &
	\mat{K}^{22}_b  &= (\mat{K}^{22}_{\mmat{\varepsilon}\mathbf{u}})^T \, (\mat{K}^{22}_{\mmat{\varepsilon}\mmat{\varepsilon}})^{-1} \, \mat{K}^{22}_{\mmat{\varepsilon}\mathbf{u}}.
\end{align}
Table \ref{tab:circulant_matrices_mixed} lists the non-zero values of the circulant blocks of the mass and stiffness matrix in the case of $\p=2, 3$ and $4$. Inverses, transposes and products of circulant matrices are circulant. Hence, the blocks in the matrix $\tilde{\mmat{K}} = \mmat{M}^{-1} \mmat{K}$ remain circulant and the technique in \ref{sec:EVcirculant} can be used to compute the discrete eigenvalues and eigenvectors of the system.

\begin{table}[h!]
\begin{center}
\caption{Circulant submatrices of the mass and stiffness matrix of the mixed formulation}
\label{tab:circulant_matrices_mixed}
\begin{tabular}{C{8em} C{4.8em} C{4.8em} C{4.8em} C{4.8em} C{4.8em}}
		\toprule
		$\frac{1}{R}\frac{1}{h}\mat{M}$	& $c_0$						& $c_1 = c_{N-1}$			& $c_2 = c_{N-2}$				& $c_3 = c_{N-3}$			& $c_4 = c_{N-4}$			\\
		\midrule
		$\p=2$	& $\frac{11}{20}$				& $\frac{13}{60}$ 				& $\frac{1}{120}$ 			& 									& 								\\ 
		$\p=3$ 	& $\frac{151}{315}$			& $\frac{397}{1680}$ 		& $\frac{1}{42}$ 			& $\frac{1}{5040}$		& 								\\
		$\p=4$ 	& $\frac{15619}{36288}$	& $\frac{44117}{181440}$	& $\frac{913}{22680}$ 	& $\frac{251}{181440}$	& $\frac{1}{362880}$\\
		\midrule
		$\frac{1}{R}\frac{1}{h}\mat{K}^{11}_{\mmat{\varepsilon}\mmat{\varepsilon}} = \frac{1}{R^3}\frac{1}{h}\frac{12}{\bar{t}^2}\mat{K}^{22}_{\mmat{\varepsilon}\mmat{\varepsilon}}$	& $c_0$						& $c_1 = c_{N-1}$			& $c_2 = c_{N-2}$				& $c_3 = c_{N-3}$			& $c_4 = c_{N-4}$			\\
		\midrule
		$\p=2$ & $\frac{2}{3}$				& $\frac{1}{6}$ 				&  			& 									& 	\\
		$\p=3$	& $\frac{11}{20}$				& $\frac{13}{60}$ 				& $\frac{1}{120}$ 			& 									& 								\\ 
		$\p=4$ 	& $\frac{151}{315}$			& $\frac{397}{1680}$ 		& $\frac{1}{42}$ 			& $\frac{1}{5040}$		& 								\\
		\midrule
		 $\mat{K}^{11}_{\mmat{\varepsilon}\mathbf{u}} = \frac{1}{R}\frac{12}{\bar{t}^2} \, \mat{K}^{21}_{\mmat{\varepsilon}\mathbf{u}}$	& $c_0 = -c_{N-1}$						& $c_1 = -c_{N-2}$		& $c_2 = -c_{N-3}$		& $c_3 = -c_{N-4}$		& $c_4 = -c_{N-5}$			\\
		\midrule
		$\p=2$ 	& $-\frac{1}{2}$ 	& $-\frac{1}{6}$ 			&			& 									&								\\
		$\p=3$ 	& $-\frac{1}{3}$	& $-\frac{5}{24}$ 		& $-\frac{1}{120}$			&			&								\\
		$\p=4$ 	& $-\frac{35}{144}$	&  $-\frac{17}{80}$ 	& $-\frac{17}{720}$ 		& $-\frac{1}{5040}$ &	\\
		\midrule
		$\frac{1}{h}\mat{K}^{12}_{\mmat{\varepsilon}\mathbf{u}}$	& $c_0=c_{N-1}$	& $c_1 = c_{N-2}$		& $c_2 = c_{N-3}$	& $c_3 =c_{N-4}$		& $c_4 = c_{N-5}$			\\
		\midrule
		$\p=2$		& 	$\frac{11}{24}$ 		& $\frac{1}{24}$ 						&				&									&										\\	
		$\p=3$		&	$\frac{151}{360}$ 		& $\frac{19}{240}$ 		& $\frac{1}{720}$ 				&			& 										\\
		$\p=4$		&	$\frac{15619}{40320}$ 		& $\frac{477}{4480}$ 	& $\frac{247}{40320}$ 			& $\frac{1}{40320}$ 			&			\\
		\midrule
		$\frac{1}{R}h\frac{12}{\bar{t}^2}\mat{K}^{22}_{\mmat{\varepsilon}\mathbf{u}}$	& $c_0 = c_{N-1}$		& $c_1 = c_{N-2}$			& $c_2 = c_{N-3}$				& $c_3 = c_{N-4}$			& $c_4 = c_{N-5}$			\\
		\midrule
		$\p=2$ 	& $\frac{1}{2}$		& $-\frac{1}{2}$ 				&  			& 								& 								\\
		$\p=3$ 	& $\frac{5}{12}$	& $-\frac{3}{8}$ 				& $-\frac{1}{24}$			& 	&								\\
		$\p=4$ 	& $\frac{49}{144}$ 			& $-\frac{21}{80}$ 			& $-\frac{11}{144}$ 		& $-\frac{1}{720}$ & 	\\
		\bottomrule
\end{tabular}
\end{center}
\end{table}

\section{Analytical eigenvalues and eigenfunctions \label{app:fourier_series}}
We represent the circumferential and radial displacement, $u$ and $w$ respectively, by the following Fourier series  \cite[Chapter 5, Section 5.3, page 82]{soedel_vibrations_2004}
\begin{align}
	u(\theta)  &= \sum_{n=1}^N U_n \sin(n \theta), \\
	w(\theta) &= \sum_{n=1}^N W_n \cos(n \theta).
\end{align}
The resulting stiffness and mass matrix have the structure
\begin{align}
		& \mmat{K} = 
			\begin{bmatrix} 
					\mat{A} & \mat{B}  \\ 
					\mat{C} & \mat{D} 
			\end{bmatrix}, &
	& \mathbf{M} = 
			\begin{bmatrix}
					\mat{M} 	&  											\\  
							     & \mat{M}
			\end{bmatrix}.&
\end{align}
Due to the orthogonality of the Fourier basis the submatrices of $\mmat{K}$ and $\mmat{M}$ are diagonal, with the diagonal entries of the mass matrix being $\pi R$. The eigenvalues and eigenmodes of the system $\mmat{K}\mmat{u}=\lambda \mmat{M}\mmat{u}$ are computed by performing an eigenvalue decomposition of the $2\times2$ block matrix 
\begin{align}
	\tilde{\mmat{K}} = \mmat{M}^{-1} \mmat{K} = 
\begin{bmatrix} 
	\tilde{\mat{A}} 	& \tilde{\mat{B}}  \\ 
	\tilde{\mat{C}}		& \tilde{\mat{D}} 
\end{bmatrix}.
\end{align}
where each of the blocks remains diagonal by virtue of $\mmat{M}$ being diagonal. The values at the diagonals are
\begin{subequations}
\begin{align}
	&  \tilde{A}_{n} = \frac{n^2}{R^2} \left( 1 + \frac{t^2}{12 R^2} \right), \\
	&  \tilde{B}_{n} = \tilde{C}_{n} = \frac{n}{R^2} \left( 1 + \frac{t^2 n^2}{12 R^2} \right), \\
	& \tilde{D}_{n} = \frac{1}{R^2} \left( 1 + \frac{t^2 n^4}{12 R^2} \right).
\end{align}
\end{subequations}
The eigenvalues and amplitudes of the Fourier modes may now be determined using the technique in  \ref{sec:circulent-block-matrix}.

\bibliographystyle{elsarticle-num}
\bibliography{paper-locking}

\begin{thebibliography}{10}
\expandafter\ifx\csname url\endcsname\relax
  \def\url#1{\texttt{#1}}\fi
\expandafter\ifx\csname urlprefix\endcsname\relax\def\urlprefix{URL }\fi
\expandafter\ifx\csname href\endcsname\relax
  \def\href#1#2{#2} \def\path#1{#1}\fi

\bibitem{bischoff2004models}
M.~Bischoff, W.~A. Wall, K.-U. Bletzinger, E.~Ramm, {M}odels and {F}inite
  {E}lements for {T}hin-walled {S}tructures, in: E.~Stein, R.~de~Borst,
  T.~J.~R. Hughes (Eds.), {E}ncyclopedia of {C}omputational {M}echanics,
  Vol.~2, John Wiley \& Sons, 2004, Ch.~3, pp. 59--137.

\bibitem{stolarski1982membrane}
H.~Stolarski, T.~Belytschko, Membrane locking and reduced integration for
  curved elements, Journal of Applied Mechanics 49 (1982) 172.

\bibitem{suri1995locking}
M.~Suri, I.~Babu{\v{s}}ka, C.~Schwab, Locking effects in the finite element
  approximation of plate models, Mathematics of Computation 64~(210) (1995)
  461--482.

\bibitem{babuvska1992locking2}
I.~Babu{\v{s}}ka, M.~Suri, Locking effects in the finite element approximation
  of elasticity problems, Numerische Mathematik 62~(1) (1992) 439--463.

\bibitem{hughes2012finite}
T.~J.~R. Hughes, The finite element method: linear static and dynamic finite
  element analysis, Courier Corporation, 2012.

\bibitem{wriggers2008nonlinear}
P.~Wriggers, Nonlinear finite element methods, Springer, 2008.

\bibitem{braess2007finite}
D.~Braess, Finite elements: Theory, fast solvers, and applications in solid
  mechanics, Cambridge University Press, 2007.

\bibitem{ambroziak2013locking}
A.~Ambroziak, Locking effects in the finite element method, Shell structures:
  Theory and application (2013) 369.

\bibitem{nguyen2022leveraging}
T.-H. Nguyen, R.~R. Hiemstra, D.~Schillinger, Leveraging spectral analysis to
  elucidate membrane locking and unlocking in isogeometric finite element
  formulations of the curved {E}uler--{B}ernoulli beam, Computer Methods in
  Applied Mechanics and Engineering 388 (2022) 114240.

\bibitem{arnold1981discretization}
D.~N. Arnold, Discretization by finite elements of a model parameter dependent
  problem, Numerische Mathematik 37 (1981) 405--421.

\bibitem{babuvska1992locking1}
I.~Babu{\v{s}}ka, M.~Suri, On locking and robustness in the finite element
  method, SIAM Journal on Numerical Analysis 29~(5) (1992) 1261--1293.

\bibitem{belytschko1985stress}
T.~Belytschko, H.~Stolarski, W.~K. Liu, N.~Carpenter, J.~S. Ong, Stress
  projection for membrane and shear locking in shell finite elements, Computer
  Methods in Applied Mechanics and Engineering 51~(1-3) (1985) 221--258.

\bibitem{macneal1985proposed}
R.~H. Macneal, R.~L. Harder, A proposed standard set of problems to test finite
  element accuracy, Finite Elements in Analysis and Design 1~(1) (1985) 3--20.

\bibitem{hughes_finite_2014}
T.~J.~R. Hughes, J.~A. Evans, A.~Reali, Finite element and {NURBS}
  approximations of eigenvalue, boundary-value, and initial-value problems,
  Computer Methods in Applied Mechanics and Engineering 272 (2014) 290--320.

\bibitem{strang_analysis_2008}
G.~Strang, G.~Fix, An {Analysis} of the {Finite} {Element} {Method}, 2nd
  Edition, Wellesley-Cambridge Press, 2008.

\bibitem{lee2015modal}
Y.~Lee, H.-M. Jeon, P.-S. Lee, K.-J. Bathe, The modal behavior of the mitc3+
  triangular shell element, Computers \& Structures 153 (2015) 148--164.

\bibitem{zou2021galerkin}
Z.~Zou, T.~J.~R. Hughes, M.~A. Scott, R.~A. Sauer, E.~J. Savitha, Galerkin
  formulations of isogeometric shell analysis: Alleviating locking with
  greville quadratures and higher-order elements, Computer Methods in Applied
  Mechanics and Engineering 380 (2021) 113757.

\bibitem{braess1988multigrid}
D.~Braess, A multigrid method for the membrane problem, Computational Mechanics
  3~(5) (1988) 321--329.

\bibitem{park1984fourier}
K.~C. Park, D.~L. Flaggs, A {F}ourier analysis of spurious mechanisms and
  locking in the finite element method, Computer Methods in Applied Mechanics
  and Engineering 46~(1) (1984) 65--81.

\bibitem{soedel_vibrations_2004}
W.~Soedel, {Vibrations of Shells and Plates}, Marcel Dekker, Inc., New York,
  2004.

\bibitem{boor_practical_2001}
C.~de~Boor, A practical guide to splines, Springer, 1978.

\bibitem{leveque2007finite}
R.~J. LeVeque, Finite difference methods for ordinary and partial differential
  equations: steady-state and time-dependent problems, SIAM, 2007.

\bibitem{silvester2000determinants}
J.~R. Silvester, Determinants of block matrices, The Mathematical Gazette
  84~(501) (2000) 460--467.

\end{thebibliography}

\end{document}